\RequirePackage{fix-cm}
\documentclass[smallextended]{svjour3}
\usepackage{url}

\usepackage{graphicx, algorithmic, algorithm2e,  psfrag, amssymb, lineno, amsfonts, amsmath, pstricks-add, bm, float, comment, color, hyperref, placeins}
\usepackage{subcaption}
\usepackage{caption}
\usepackage[]{natbib}
\usepackage[misc,geometry]{ifsym} 
\RestyleAlgo{boxruled}
\hypersetup{colorlinks,breaklinks,
           linkcolor=blue,urlcolor=blue,
           anchorcolor=blue,citecolor=blue}

\newcommand{\am}{\mathbf{a}}
\newcommand{\Cm}{\mathbf{C}}
\newcommand{\dm}{\mathbf{d}}
\newcommand{\Gm}{\mathbf{G}}
\newcommand{\gm}{\mathbf{g}}

\newcommand{\hb}{\mathbf{h}}

\newcommand{\mm}{\mathbf{m}}
\newcommand{\ngg}{n_{\mathbf{g}}}
\newcommand{\np}{n_{\mathbf{\ppm}}}
\newcommand{\ny}{n_{\mathbf{\Ym}}}
\newcommand{\Ns}{N_\mathrm{s}}
\newcommand{\ppm}{{\pmb\theta}}

\newcommand{\vm}{\mathbf{v}}
\newcommand{\xm}{\mathbf{x}}
\newcommand{\Ym}{\pmb{\xi}}

\newcommand{\ie}{\textit{i.e.},}
\newcommand{\eg}{\textit{e.g.},}
\newcommand{\tpos}[1]{#1^{\mathrm{T}}}
\newcommand{\Exp}{\mathbb{E}}
\newcommand{\Var}{\mathrm{Var}}

\newcommand{\kappaa}{\pmb\kappa}
\newcommand{\thetaa}{\pmb\theta}

\usepackage{pgfplots}
\usepackage{mathptmx}

\usepackage{tikz}
\usetikzlibrary{matrix} 
\usetikzlibrary{arrows} 
\usetikzlibrary{arrows.meta}
\usetikzlibrary{calc} 
\usetikzlibrary{shapes}
\usetikzlibrary{fit}
\usetikzlibrary{positioning}
\usetikzlibrary{intersections,patterns,pgfplots.fillbetween}
\usetikzlibrary{calc}
\usetikzlibrary{decorations.pathreplacing}
\usepackage{ifpdf}
\ifpdf
  \usepackage{epstopdf}
\fi

\tikzstyle{block} = [draw,rectangle,thick,minimum height=2em,minimum width=2em]
\tikzstyle{sum} = [draw,circle,inner sep=0mm,minimum size=2mm]
\tikzstyle{connector} = [->,thick]
\tikzstyle{line} = [thick]
\tikzstyle{branch} = [circle,inner sep=0pt,minimum size=1mm,fill=black,draw=black]
\tikzstyle{guide} = []

\begin{document}



\title{Topology Optimization under Uncertainty {using a} Stochastic Gradient-based {Approach}
}


\author{Subhayan De        \and
	Jerrad Hampton \and  Kurt Maute \and Alireza Doostan
}


\institute{
	S. De \at
	Smead Aerospace Engineering Sciences Department, University of Colorado, Boulder, CO 80309, USA\\
	\email{Subhayan.De@colorado.edu}         
\and
 J. Hampton\\ \email{Jerrad.Hampton@colorado.edu}\\
\and
 K. Maute\\ \email{Kurt.Maute@colorado.edu}\\
\and 
A. Doostan (\Letter)\\ \email{Alireza.Doostan@colorado.edu}
}

\date{Received: date / Accepted: date}

\maketitle



\begin{abstract}
Topology optimization under uncertainty (TOuU) often defines objectives and constraints by statistical moments of geometric and physical quantities of interest. Most traditional TOuU methods use gradient-based optimization algorithms and rely on accurate estimates of the statistical moments and their gradients, \eg~ via adjoint calculations. When the number of uncertain inputs is large or the quantities of interest exhibit large variability, a large number of adjoint (and/or forward) solves may be required to ensure the accuracy of these gradients. The optimization procedure itself often requires a large number of iterations, which may render TOuU computationally expensive, if not infeasible. To tackle this difficulty, we here propose an optimization approach that generates a stochastic approximation of the objective, constraints, and their gradients via a small number of adjoint (and/or forward) solves, per optimization iteration. A statistically independent (stochastic) approximation of these quantities is generated at each optimization iteration. The total cost of this approach is only a small factor larger than that of the corresponding deterministic TO problem. We incorporate the stochastic approximation of objective, constraints and their design sensitivities into two classes of optimization algorithms. First, we investigate the stochastic gradient descent (SGD) method and a number of its variants, which have been successfully applied to large-scale optimization problems for machine learning. Second, we study the use of the proposed stochastic approximation approach within conventional nonlinear programming methods, focusing on the Globally Convergent Method of Moving Asymptotes (GCMMA). The performance of these algorithms is investigated with structural design optimization problems utilizing a Solid Isotropic Material with Penalization (SIMP), as well as an explicit level set method. These investigations, conducted on both two- and three-dimensional structures, illustrate the efficacy of the proposed stochastic gradient approach for TOuU applications. 
\end{abstract}




\section{\texorpdfstring{Introduction}{Introduction}}
\label{sec:intro}

Finding the optimum geometry and arrangement of materials for the design of engineered materials and components is a key challenge across many applications. This challenge becomes even more complicated in the presence of uncertainty resulting from, \eg~ manufacturing imprecision or incomplete/inaccurate measurements. In structural topology optimization {(TO)}~\citep{Sigmund1998,sigmund2013topology}, the arrangement of one or multiple materials within a design domain is investigated to optimize the mechanical performance of a structure while accounting for design constraints. In recent years, TO has been used in several other fields, such as fluid flows, acoustics, optics, and multi-physics applications. The reader is referred to the survey papers by \cite{sigmund2013topology, deaton2014survey}.

In the presence of uncertainty, we often seek \textit{robust} designs that provide a compromise between mean performance and insensitivity to variation in geometry, material properties, such as elastic modulus or heat conductivity of the material, and operating conditions. 
A poor understanding of uncertainties in these system inputs can lead to poor design choices \citep{schueller2008computational}. 
In this paper, we consider TO problems with uncertainties in geometry, loading conditions, and material properties, which we model as random variables. 
We refer to this TO in the presence of uncertainty as Topology Optimization under Uncertainty (TOuU), similar to ~\cite{conti2009shape,tootkaboni2012topology,maute2014touu,keshavarzzadeh2017topology}. 

The effects of uncertainty in TO were first considered by adding a probabilistic constraint, mostly based on the probability of failure \citep{bae2002reliability,maute2003reliability,kharmanda2004reliability,jung2004reliability,moon2004reliability,kim2006reliability,mogami2006reliability,eom2011reliability}. This approach is known as reliability-based topology optimization (RBTO) and mostly uses first-order or second-order Taylor series expansion to approximate the moments of the limit state functions with respect to the uncertain parameters \citep{haldar2000probability}. On the other hand, \cite{eldred2011design} and \cite{keshavarzzadeh2016gradient} used polynomial chaos expansion method \citep{roger2003stochastic} to solve RBTO. In these studies, the uncertainty is assumed to be present in the loads, in the material properties, such as the modulus of elasticity and yield stress, or in the thickness of the structural components. 

Robust TO, on the other hand, considers the effects of uncertainty in the objective and constraints through higher order statistics, \textit{e.g.}, variance. \cite{alvarez2005minimization}, \cite{dunning2011introducing}, and \cite{dunning2013robust} considered uncertainties in the load magnitude and direction and minimized the expected and/or variance of the compliance of the structure. 
TO under geometric uncertainty has also been investigated in \cite{guest2008structural}, \cite{chen2010level}, and \cite{chen2011new}. These authors used a Karhunen-Loeve expansion to represent the random field resulting from material, geometry and loading uncertainties. Perturbations in the stiffness matrix are used in \cite{asadpoure2011robust},  \cite{zhao2014robust}, \cite{jansen2015robust}, and \cite{kriegesmann2019robust} to estimate the statistical moments of the compliance by
Monte Carlo approximations. \cite{tootkaboni2012topology}, \cite{lazarov2012topology}, \cite{keshavarzzadeh2017topology}, and \cite{zhang2017robust} used polynomial chaos expansions while considering uncertainties in the loading and geometry of the structure and estimated the mean and variance of the compliance. A combined reduced-order modeling  and polynomial chaos expansion approach for shape optimization has been proposed in \cite{maute2009reduced}.

Monte Carlo simulation approaches use a large number of gradient evaluations at independent samples of inputs to form sample average estimates of the gradients. The number of required evaluations is large when the variances of the gradients are large. This, in turn, will lead to a high computational cost. Taylor series based perturbations are computationally efficient; however, they lose accuracy when the objective and constraints depend on the input uncertainties in a highly nonlinear manner. While methods based on polynomial chaos expansion are frequently used, they suffer from the \textit{curse of dimensionality}, \textit{i.e.}, the number of expansion coefficients increases rapidly with the number of uncertain parameters. Sparse polynomial chaos expansions \citep{Doostan09b,doostan2011non,blatman2010adaptive,hampton2016compressive} can be used to reduce the computational cost, but for problems with a large number of uncertain inputs, such as TO with multiple uncertainty sources, the computational cost associated with this approach also becomes unbearable.

To alleviate the computational burden and motivated by the recent progress in stochastic gradient descent (SGD) methods,~\citep{bottou2018optimization}, we advocate for a new approach for TOuU that uses a stochastic approximation of the gradients. SGD methods have recently seen significant use for high-dimensional problems arising in machine learning~\citep{bottou2010SGDML, bottou2012SGDML, sutskever2013SGDML}
In detail, we generate an unbiased approximation of the gradients in a manner similar to the standard Monte Carlo approach but with a small, \eg~ 4 or 10, number of Monte Carlo samples of the gradient. Unlike in the standard Monte Carlo approach, we generate statistically independent estimates of the gradient in each optimization iteration. The small sample size gradient estimates are essentially stochastic approximations of those computed via standard Monte Carlo approach, hence the term {\it stochastic gradients}. Using these stochastic gradients, we will illustrate in Section \ref{sec:sgd} that TOuU problems may be solved via gradient descent methods. 

Recently, \cite{martin2018analysis} analyzed a standard SGD method with variable step size to solve an optimal control problem, where the underlying partial-differential equation is solved using the finite element method and used as a constraint. 
{For a desired tolerance value in the optimal solution Martin et al. (2018) showed that when compared to the standard Monte Carlo method with fixed random samples the use of SGD can reduce the computational cost by a factor equal to the natural logarithm of the desired tolerance.}
In contrast to these investigations, in the present study, we {consider} the use of improved variants of SGD for non-convex TO problems, and to our best knowledge, this is the first such investigation. We explore AdaGrad~\citep{duchi2011adagrad} and Adam \citep{kingma2014adam}, which are variants of the standard SGD method that retard the movement in {directions} with historically large gradient magnitudes and show promise for problems that lack convexity. To a lesser extent, we investigate methods like Stochastic Average Gradient (SAG)~\citep{rouxetalSAG2012}, which performs well on objectives that are strongly-convex, stochastic variance reduced gradient (SVRG) \citep{johnson2013accelerating}, and Adadelta~\citep{zeiler2012adadelta}. Second, we study the performance of a standard {nonlinear} programming method when supplying stochastic gradients. We focus on the Globally Convergent Method of Moving Assympotes (GCMMA) by ~\cite{svanberg1987}, owing to its popularity for TO. 

In this paper, we describe the geometry and material distribution by either an explicit level set {or a density based method. }
The explicit level set method follows the work of \cite{kreissl2012levelset}, \cite{villanueva2014density}, and \cite{sharma2017sensitivities}. For the density approach, we build on the Solid Isotropic Material with Penalization (SIMP) method~\citep{bendsoe1989optimal,ZhouSIMP91,sigmund2013topology}. 
We illustrate the use of stochastic gradients in SGD methods and the GCMMA with two static linear elastic problems. Uncertainty in geometry, material properties, and loading conditions is considered. 
%
The results show that Adam and GCMMA supplied with stochastic gradients produce designs that are robust in presence of uncertainty. Further, only a handful of random samples are needed at every iteration, thus reducing the cost of TOuU to a factor that is only slightly larger than the one of deterministic TO.


The rest of the paper is organized as follows. Section \ref{sec:methodology} gives a brief background of the TO approaches used in this paper and provides a short outline on different {SGD methods and GCMMA}. Then, we formulate the TOuU problem. In Section \ref{sec:examples},  we illustrate the performance of the stochastic gradient-based methods with two TO problems under uncertainty. { Finally, in Section \ref{sec:conc}, we conclude the paper with a discussion on our observations.}


\section{Methodology}
\label{sec:methodology}


\subsection{Design Model}\label{sec:des_model}

In this paper, we describe the geometry and material layout within a specified design domain by an explicit level set method \citep{van2013level} {or a density method, namely the SIMP approach. Two different numerical examples are solved herein using these two methods separately to show the versatile applicability of the proposed TOuU approach}. The parameters of the discretized level set functions are defined explicitly in terms of the optimization variables. To this end, we follow closely the work of \cite{kreissl2012levelset,villanueva2014density}. 
{For the SIMP method, we use the settings described in \cite{bendsoe1989optimal} and \cite{ZhouSIMP91} as well as projection \cite{sigmund2013topology} as described in this section.} 

To account specifically for shape imperfections, we consider structures made of an assembly of geometric primitives \citep{guo2015morphable,norato2015geometry,zhang2016new,liu2018efficient} {in one of our examples herein}. Each primitive is defined by a small number of parameters that are treated as optimization and/or random variables. {In that example,
a 2D structure is composed of $n_r$ rectangles or bars.} The level set field associated with {$r$th rectangle}, $\phi^{\mathrm{rec}}_r$, is given in a local coordinate system, $\left({\tilde{x},\tilde{y}}\right)$, as follows:
%
\begin{equation}\label{eq:rect_def}
    \phi^{\mathrm{rec}}_r=\left({\left({\frac{\tilde{x}}{a_r}}\right)^\mu+\left({\frac{\tilde{y}}{b_r}}\right)^\mu }\right)^{1/\mu}-1,
\end{equation}
where the parameters $a_r$ and $b_r$ denote the dimensions of the {$r$th rectangle}, and $\mu$ controls the sharpness of its corners. Note that the value of the level set function within the rectangle is negative. The local and the global coordinates are related by the following coordinate transformation:
\begin{equation}
    \left[{
    \begin{array}{c}
    \tilde{x}\\
    \tilde{y}
    \end{array}
    }\right] 
    =
    \left[{
    \begin{array}{cc}
      \cos(\alpha_r)    &  \sin(\alpha_r) \\
     -\sin(\alpha_r)   &   \cos(\alpha_r)
    \end{array}
    }\right] 
    \
    \left[{
    \begin{array}{c}
    x-x_r^c\\
    y-y_r^c
    \end{array}
    }\right] ,
\end{equation}
where the center for the rectangle is defined by $\left({x_r^c,y_r^c}\right)$ and $\alpha_r$ denotes the angle between the coordinate axes. The assembly of bars within the design domain is defined by the minimum of all level set values {$\{\phi^{\mathrm{rec}}_r\}_{r=1}^{n_r}$} at a point $\xm$. To obtain a differentiable function $\phi(\xm)$ we approximate the minimum level set value at Node $i$ by the Kreisselmeier-Steinhauser function:
\begin{equation}\label{eq:ks}
    \phi_i=\frac{1}{\beta_\mathrm{KS}} \ln \left({\sum\limits_{r=1}^{n_r}{e^{\beta_{KS} \ \phi^{\mathrm{rec}}_r(\xm_i)}}}\right) \qquad \beta_{KS}<0,
\end{equation}
where $\xm_i$ are the coordinates at Node $i$, and the parameter $\beta_\mathrm{KS}$ controls the accuracy of the approximation. Assuming that all bars are made of a single material, the geometry of the structure is defined by $n_r$ sets of the following five parameters: $\left\{{a_r,b_r,x_r^c,y_r^c,\alpha_r}\right\}$. These parameters will be treated as design and random parameters in the numerical examples of Section \ref{sec:mbb}. 

{
For our second example, we use SIMP to formulate the optimization problem. In SIMP, the material properties are interpolated by a power-law model in terms of the density $\rho$ of a fictitious porous material, \eg
\begin{equation}\label{eq:simp}
E(\rho) = \rho^{\beta_\mathrm{SIMP}} E_0,
\end{equation}
where $\beta_\mathrm{SIMP}$ is the penalization parameter and $E_0$ is the bulk material's elastic modulus. For $\beta_\mathrm{SIMP}>1$, intermediate densities are penalized and thus favors a $0-1$ design for $\rho$. 
Here, we assign each node $i, \ i=1, \dots, n_N$, in the finite element mesh an optimization variable $\theta_i$.
To avoid a checker-board design, we treat the elemental volume fractions, or normalized densities $\rho$, of the fictitious porous material as a weighted combination of the optimization variables $\ppm$ using the following linear filter \citep{bruns2001topology,bourdin2001filters,sigmund2007morphology,andreassen2011efficient}: 
\begin{equation}
\label{eq:linfilter}
{\rho}_i = \frac{\sum\limits_{j=1}^{n_N} {w_{ij} \theta_j}}{\sum\limits_{j=1}^{n_N} {w_{ij}}}
\quad \textrm{with} \ w_{ij}=\max \left({0,r_f - \left|{ \xm_i - \xm_j }\right|}\right) \ ; \ 0 \leq \rho_i \leq 1 \ ; \ i=1 \ldots N_e,
\end{equation}
where $r_f$ is a filter radius and the design domain is divided in $N_e$ number of elements using a finite element mesh. 
To counteract the emergence of intermediate densities introduced by the filter \eqref{eq:linfilter}, we project the filtered densities by a smoothed Heaviside function \citep{sigmund2013topology,wang2011projection} given by 
\begin{equation}\label{eq:projection}
    \bar\rho_i = \frac{\tanh(\beta_{\mathrm{pr}}(\rho_i-\nu_\mathrm{pr}))+\tanh(\beta_\mathrm{pr}\nu_\mathrm{pr})}{\tanh(\beta_\mathrm{pr}(1-\nu_\mathrm{pr}))+\tanh(\beta_\mathrm{pr}\nu_\mathrm{pr})}, \quad i=1 \ldots N_e,
\end{equation}
where $\bar\rho_i$ is the projected density; $\beta_{\mathrm{pr}}$ is a projection strength parameter; and $\nu_\mathrm{pr}$ is a projection threshold parameter. 
}

\subsection{Analysis Model}\label{sec:ana_model}

This study considers structural optimization problems where the physical response is described by a static linear elastic model. When describing the geometry of the primitives by a level set method, we follow the work of \citep{villanueva2014density,sharma2017sensitivities} discretize the governing equations by the eXtended Finite Element Method (XFEM). A generalized Heaviside enrichment strategy is used to consistently approximate the displacement field in the solid domain. 
To mitigate ill-conditioning caused by particular intersection configurations, we stabilize the XFEM formulation by the face-oriented Ghost Penalty Method \citep{burman2014fictitious,shott2014face}. 
For a specific realization of the design and boundary conditions, we compute the design sensitivities, \ie~ the gradients of the objective and gradients measures 
by the adjoint method. Details of evaluating the sensitivities using the adjoint method can be found in ~\cite{sharma2017sensitivities}. 
When describing the material layout in the design design domain by the SIMP approach, we adopt a standard finite element scheme. The design sensitivities are computed by the adjoint method, accounting for the smoothing and projection operators. The reader is referred to Sigmund and Maute \cite{sigmund2013topology}, and Deaton and Grandhi \cite{deaton2014survey} for details. 
Depending on the optimization algorithm, the gradients of objective and constrains with respect to the optimization variables are used directly or modified to determine the search direction in the optimization process. Details for the optimization algorithms considered in this study are presented in the next subsection. 
\subsection{Topology Optimization under Uncertainty (TOuU) Problem Formulation}\label{sec:method_details}

Let $\ppm\in \mathbb{R}^{\np}$ denote the vector of optimization variables, and $\Ym\in\mathbb{R}^{\ny}$ the vector of random variables associated with the problem. Let $f(\ppm;\Ym): \mathbb{R}^{\np} \times \mathbb{R}^{\ny} \rightarrow \mathbb{R}$ denote the deterministic performance measure given $\ppm$, which also depends on the realized values of $\Ym$. Similarly, let $\gm(\ppm;\Ym): \mathbb{R}^{\np} \times \mathbb{R}^{\ny} \rightarrow \mathbb{R}^{\ngg}$ be $\ngg$ real-valued deterministic measure of the constraints, which potentially depend on realized values of $\Ym$. We say that $(\ppm;\Ym)$ satisfies the constraints if $\gm(\ppm;\Ym) \le 0$ and refer to any positive value of $\gm(\ppm;\Ym)$ as a constraint violation. In structural TO, the function $f$ is often the strain energy and the constraint functions, $\gm$, often include the volume or mass of the structure.

In the present study, for TOuU, we define the objective as a combination of its expected value and variance as follows
\begin{equation}\label{eq:exp_risk}
R(\ppm) = \Exp[f(\ppm;\Ym)] + \lambda \Var(f(\ppm;\Ym)),\\
\end{equation}
%
%
where $\Exp[\cdot]$ and $\Var(\cdot)$ denote the mathematical expectation and variance of their arguments, respectively. Here we simply add the variance to the objective {following \cite{beyer2007robust}} so that the calculation of the gradients will be straightforward. 
Similarly, {we use the constraint violation
\begin{equation}
\label{eq:glam_def}
\Cm_j(\ppm) = \Exp[\Gm_j(\ppm;\Ym)] + \lambda \Var(\Gm_j(\ppm;\Ym)), \quad j=1,\dots,\ngg,
\end{equation}  
where $\Gm_j(\ppm;\Ym) = \left(\max\left(0,\gm_j(\ppm;\Ym)\right)\right)^2$, for $j=1,\dots,\ngg$. We refer the interested reader to, e.g., \cite{griva2009linear}[Chapter 6] for a discussion on this choice of $\Gm_j$.
}
%
We note that while $f$ and $\gm$ are a random variable {and a random vector, respectively}, $R$ and $\Cm$ are not as they are associated with moments taken with respect to the probability measure of $\Ym$. Here, the parameter $\lambda \ge 0$ denotes the importance of variations in $f$ or $\gm$, relative to their means. {Note that, $\lambda$ in \eqref{eq:exp_risk} has the same unit as ${1}/{f(\ppm,\Ym)}$. Similarly, in \eqref{eq:glam_def}, the unit of $\lambda$ is same as the unit of {$1/\gm_j^2(\ppm;\Ym)$}. Here, we keep the same notation for brevity and use the same value of $\lambda$ in \eqref{eq:exp_risk} and \eqref{eq:glam_def} while acknowledging different values of $\lambda$ in \eqref{eq:exp_risk} and \eqref{eq:glam_def} may be used.

} The effect of adding the variance to the objective and constraints is to aim for a design that shows little variations in the two responses, namely, $f$ and $\gm$ even in the presence of uncertainty. 
Hence, we denote this formulation as \textit{robust optimization} formulation for $\lambda>0$.
Although one may use different values of $\lambda$ in (\ref{eq:exp_risk}) and (\ref{eq:glam_def}) here we restrict our investigation to utilizing the same $\lambda$ for both functions.
Therefore, we are interested here in solving a problem of the type,
\begin{align}
\label{eq:opt_def}
\mathop{\min~}\limits_{\ppm}R(\ppm) \mbox{   subject to } \Cm(\ppm) = \mathbf{0},
\end{align}
or the closely related unconstrained formulation,
\begin{align}
\label{eq:opt_def2}
\mathop{\min~}\limits_{\ppm}R(\ppm) + \tpos{\kappaa} {\Cm}(\ppm), 
\end{align}
{where 
$\kappaa$ is a user-specified vector of parameters that penalizes against constraint violation. 
}  
%
%
For TO, an extensive search for $\kappaa$ is not computationally feasible. Hence, based on a few preliminary runs the values of $\kappaa$ were determined in the numerical examples of this paper.

In standard Monte Carlo approach, we approximate $R(\ppm)$ and $\Cm(\ppm)$, from (\ref{eq:exp_risk}) and (\ref{eq:glam_def}) utilizing $\Ns$ forward solves of the model for specific values of the design variables $\ppm$ and $\Ns$ realizations of $\Ym$. However, evaluating the functions $f(\ppm,\Ym)$ and $\gm(\ppm;\Ym)$ at each iteration $\Ns$ times may be computationally expensive for structures with many degrees of freedom. 
Herein, we propose to use $n$ number of random samples to estimate (\ref{eq:exp_risk}) and (\ref{eq:glam_def}) where $n\ll \Ns$ (\eg~$n=4,10,$ etc.).  
We use a hat notation to denote an estimate of a quantity, and define the following
\begin{align}
\widehat{\mathbb{E}}[f(\ppm)] &:= \frac{1}{n}\mathop{\sum}\limits_{i=1}^{n} f(\ppm;\Ym_i);\\
\widehat{\mathbb{E}}[\Gm(\ppm)] &:= \frac{1}{n}\mathop{\sum}\limits_{i=1}^{n} \Gm(\ppm;\Ym_i);\\
\widehat{\Var}(f(\ppm)) &:= \frac{1}{n-1}\mathop{\sum}\limits_{i=1}^{n} \bigg( f(\ppm;\Ym_i) - \widehat{\mathbb{E}}[f(\ppm)]\bigg)^2;\\
\widehat{\Var}(\Gm_j(\ppm)) &:= \frac{1}{n-1}\mathop{\sum}\limits_{i=1}^{n} \bigg( \Gm_j(\ppm;\Ym_i) - \widehat{\mathbb{E}}[\Gm_j(\ppm)]\bigg)^2;\quad j=1,\dots,\ngg;\\
\label{eq:Fhatdef}
\widehat{R}(\ppm) &= \widehat{\mathbb{E}}[f(\ppm)] + \lambda \widehat{\Var}(f(\ppm));\\\
\label{eq:Ghatdef}
\widehat{\Cm}_j(\ppm) &= \widehat{\mathbb{E}}[\Gm_j(\ppm)] + \lambda \widehat{\Var}(\Gm_j(\ppm));\quad j=1,\dots,\ngg.
\end{align}
%
For brevity, these definitions do not explicitly state the dependence on $\Ym$.  

We then approximate the gradients of the objective and constraints with respect to the design parameters via the gradients of (\ref{eq:Fhatdef}) and (\ref{eq:Ghatdef}), which are computed via chain rule to be
\begin{align}
\widehat{\mathbb{E}}[\nabla f(\ppm)] &:= \frac{1}{n}\mathop{\sum}\limits_{i=1}^{n} \nabla f(\ppm;\Ym_i);\label{eq:EdF}\\ 
\widehat{\mathbb{E}}[\nabla \Gm(\ppm)] &:= \frac{1}{n}\mathop{\sum}\limits_{i=1}^{n} \nabla \Gm(\ppm;\Ym_i);\\
\widehat{\mathbb{E}}[f(\ppm)\nabla f(\ppm)] &:= \frac{1}{n}\mathop{\sum}\limits_{i=1}^{n} f(\ppm;\Ym_i)\nabla f(\ppm;\Ym_i));\\
\widehat{\mathbb{E}}[\Gm(\ppm)\nabla \Gm(\ppm)] &:= \frac{1}{n}\mathop{\sum}\limits_{i=1}^{n} \Gm(\ppm;\Ym_i)\nabla \Gm(\ppm;\Ym_i));\\
\widehat{\Var}(\nabla f(\ppm)) &:= \frac{2n}{n-1} {\left[\widehat{\mathbb{E}}[f(\ppm)\nabla f(\ppm)] - \widehat{\mathbb{E}}[f(\ppm)]\widehat{\mathbb{E}}[\nabla f(\ppm)]\right]};\\
\widehat{\Var}(\nabla \Gm(\ppm)) &:= \frac{2n}{n-1}\left[\widehat{\mathbb{E}}[\Gm(\ppm)\nabla \Gm(\ppm)] - \widehat{\mathbb{E}}[\Gm(\ppm)]\widehat{\mathbb{E}}[\nabla \Gm(\ppm)]\right];\label{eq:Fhatgraddef} \\
\nabla\widehat{R}(\ppm) &= \widehat{\mathbb{E}}[\nabla f(\ppm)] + \lambda \widehat{\Var}(\nabla f(\ppm)); \label{eq:Ghatgraddef}\\
\nabla\widehat{\Cm}_j(\ppm) &= \widehat{\mathbb{E}}[\nabla \Gm(\ppm)] + \lambda \widehat{\Var}(\nabla \Gm(\ppm));\quad j =1,\dots,\ngg, \label{eq:EdC}
\end{align}
where the normalization terms $n/(n-1)$ are to ensure unbiasedness. 
Notice that (\ref{eq:EdF})-(\ref{eq:EdC}) hold as expectation and differentiation are linear operators. We note that in this case, the expectation and variance are computed component-wise; specifically, there is no covariance information relying on correlations between different components. Finally, the gradient of the objective formulation \eqref{eq:opt_def2}, which combines gradients of the objective and the constraint functions, is given by
\begin{align}
\label{eq:Hhatgrad}
\widehat{\hb}(\ppm;\kappaa;\lambda) &:= \nabla\widehat{R}(\ppm) + \tpos{\kappaa}{\nabla{\widehat\Cm}}(\ppm).
\end{align}
%

\begin{figure}[htb!]
    \centering
    \begin{subfigure}[t]{\textwidth}
    \begin{tikzpicture}
    \node[draw=black,fill = yellow!20,minimum width = 11.5cm,minimum height = 9.5cm] () at (3,-2) {};
    \node[draw=none] () at (-1,0) {Random samples};
    
    \node[draw=none] () at (-1,-6) {Optimization iteration};
    \node[draw=none] () at (-1,-4.25) {Evaluate $\widehat{R}(\ppm)\in \mathbb{R}$,};
    \node[draw=none] () at (-1,-4.75) {$\widehat{\Cm}(\ppm)\in \mathbb{R}^{\ngg\times 1}$,};
    \node[draw=none] () at (-1,-5.25) {and $\widehat{\hb}(\ppm)\in \mathbb{R}^{n_{\ppm}\times 1}$};
\scriptsize
\foreach \j in {1,...,4}
{
\foreach \i in {1,...,3}
{
        \pgfmathtruncatemacro{\y}{(\i - 1)};
        \pgfmathtruncatemacro{\x}{\j}; 
        \pgfmathtruncatemacro{\label}{3-\y};
        \node[circle,thick,draw=black,fill=white!60!green,minimum size=22]
        (\label) at (1*\x,1*\y) {$\Ym_{\label}$};
}
}
\foreach \j in {1,...,4}
{
\foreach \i in {1,...,3}
{
        \pgfmathtruncatemacro{\y}{(\i - 6)};
        \pgfmathtruncatemacro{\x}{\j}; 
        \pgfmathtruncatemacro{\label}{5-\y};
        \node[circle,draw=black,fill=white!50!black,minimum size=0.25,inner sep=1pt]
        (\label) at (1*\x,0.25*\y) {};
}
}

\foreach \i in {1,...,3}
{
        \pgfmathtruncatemacro{\y}{(4)};
        \pgfmathtruncatemacro{\x}{20+\i}; 
        \node[circle,draw=black,fill=white!50!black,minimum size=0.25,inner sep=1pt]
        () at (0.25*\x,0.25*\y) {};
}

\foreach \j in {1,...,2}
{
\foreach \i in {1,...,3}
{
        \pgfmathtruncatemacro{\y}{(\i - 1)};
        \pgfmathtruncatemacro{\x}{6+\j}; 
        \pgfmathtruncatemacro{\label}{3-\y};
        \node[circle,thick,draw=black,fill=white!60!green,minimum size=22]
        (\label) at (1*\x,1*\y) {$\Ym_{\label}$};
}
}
\foreach \j in {1,...,2}
{
\foreach \i in {1,...,3}
{
        \pgfmathtruncatemacro{\y}{(\i - 6)};
        \pgfmathtruncatemacro{\x}{6+\j}; 
        \pgfmathtruncatemacro{\label}{5-\y};
        \node[circle,draw=black,fill=white!50!black,minimum size=0.25,inner sep=1pt]
        (\label) at (1*\x,0.25*\y) {};
}
}
\foreach \j in {1,...,2}
{
\foreach \i in {1,...,2}
{
        \pgfmathtruncatemacro{\y}{(\i - 4)};
        \pgfmathtruncatemacro{\x}{6+\j}; 
        \pgfmathtruncatemacro{\label}{3+\y};
        \ifthenelse{\label=0}{
        \node[circle,thick,draw=black,fill=white!60!green,minimum size=22]
        (\label) at (1.0*\x,1.0*\y) {$\Ym_{\Ns}$};}{
        \node[circle,thick,draw=black,fill=white!60!green,minimum size=22]
        (\label) at (1.0*\x,1.0*\y) {$\Ym_{\Ns-\label}$};
        }
}
}

\foreach \j in {1,...,4}
{
\foreach \i in {1,...,2}
{
        \pgfmathtruncatemacro{\y}{(\i - 4)};
        \pgfmathtruncatemacro{\x}{\j}; 
        \pgfmathtruncatemacro{\label}{3+\y};
        \ifthenelse{\label=0}{
        \node[circle,thick,draw=black,fill=white!60!green,minimum size=22]
        (\label) at (1.0*\x,1.0*\y) {$\Ym_{\Ns}$};}{
        \node[circle,thick,draw=black,fill=white!60!green,minimum size=22]
        (\label) at (1.0*\x,1.0*\y) {$\Ym_{\Ns-\label}$};
        }
}
}

\foreach \j in {1,...,4}
{
        \pgfmathtruncatemacro{\y}{( -6)};
        \pgfmathtruncatemacro{\x}{\j}; 
        \pgfmathtruncatemacro{\optiter}{\x};

        \node[thick,draw=black,fill=white!80!blue,minimum size=15]
        (\optiter) at (1.0*\x,1.0*\y) {$\ppm_{\optiter}$};
        
        \draw[-latex, line width=1mm] (\x,\y+2.5) -- (\x,\y+0.5);
        \draw[-latex,thick] (\x+0.25,\y) -- (\x+0.75,\y);
}

\foreach \j in {1,...,2}
{

        \pgfmathtruncatemacro{\y}{( - 6)};
        \pgfmathtruncatemacro{\x}{6+\j}; 
        \pgfmathtruncatemacro{\optiter}{8-\x};
        \ifthenelse{\optiter=0}{
        \node[thick,draw=black,fill=white!80!blue,minimum size=15]
        (\optiter) at (1.0*\x,1.0*\y) {$\ppm_{M_f}$};}{
        \node[thick,draw=black,fill=white!80!blue,minimum size=15]
        (\optiter) at (1.0*\x,1.0*\y) {$\ppm_{M_f-\optiter}$};
        }
        \draw[-latex, line width=1mm] (\x,\y+2.5) -- (\x,\y+0.5);
        \draw[-latex,thick] (\x-0.675,\y) -- (\x-0.3,\y);
}

\node[circle,thick,draw=black,fill=white!60!green,minimum size=22]
        () at (10,-1) {$\Ym_{i}$};
\node[circle,thick,draw=black,fill=white!40!red,minimum size=22]
        () at (10,0) {$\Ym_{i}$};
\node[draw=none,right]
        () at (10.5,0) {: Not used};
\node[draw=none,right]
        () at (10.5,-1) {: Used};

    \end{tikzpicture}
    		\caption{Schematic of the full-gradient approach} \label{fig:FGschematic}
	\end{subfigure}\\
	\vspace{10pt}
	    \begin{subfigure}[t]{\textwidth}
    \begin{tikzpicture}
    \node[draw=black,fill = yellow!20,minimum width = 11.5cm,minimum height = 9.5cm] () at (3,-2) {};
    \node[draw=none] () at (-1,0) {Random samples};
    
    \node[draw=none] () at (-1,-6) {Optimization iteration};
    \node[draw=none] () at (-1,-4.25) {Evaluate $\widehat{R}(\ppm)\in \mathbb{R}$,};
    \node[draw=none] () at (-1,-4.75) {$\widehat{\Cm}(\ppm)\in \mathbb{R}^{\ngg\times 1}$,};
    \node[draw=none] () at (-1,-5.25) {and $\widehat{\hb}(\ppm)\in \mathbb{R}^{n_{\ppm}\times 1}$};
\scriptsize
\foreach \j in {1,...,4}
{
\foreach \i in {1,...,3}
{
        \pgfmathtruncatemacro{\y}{(\i - 1)};
        \pgfmathtruncatemacro{\x}{\j}; 
        \pgfmathtruncatemacro{\label}{3-\y};
        \node[circle,thick,draw=black,fill=white!40!red,minimum size=22]
        (\label) at (1*\x,1*\y) {$\Ym_{\label}$};
}
}

\foreach \j in {1,...,4}
{
\foreach \i in {1,...,3}
{
        \pgfmathtruncatemacro{\y}{(\i - 6)};
        \pgfmathtruncatemacro{\x}{\j}; 
        \pgfmathtruncatemacro{\label}{5-\y};
        \node[circle,draw=black,fill=white!50!black,minimum size=0.25,inner sep=1pt]
        (\label) at (1*\x,0.25*\y) {};
}
}

\foreach \i in {1,...,3}
{
        \pgfmathtruncatemacro{\y}{(4)};
        \pgfmathtruncatemacro{\x}{20+\i}; 
        \node[circle,draw=black,fill=white!50!black,minimum size=0.25,inner sep=1pt]
        () at (0.25*\x,0.25*\y) {};
}

\foreach \j in {1,...,2}
{
\foreach \i in {1,...,3}
{
        \pgfmathtruncatemacro{\y}{(\i - 1)};
        \pgfmathtruncatemacro{\x}{6+\j}; 
        \pgfmathtruncatemacro{\label}{3-\y};
        \node[circle,thick,draw=black,fill=white!40!red,minimum size=22]
        (\label) at (1*\x,1*\y) {$\Ym_{\label}$};
}
}
\foreach \j in {1,...,2}
{
\foreach \i in {1,...,3}
{
        \pgfmathtruncatemacro{\y}{(\i - 6)};
        \pgfmathtruncatemacro{\x}{6+\j}; 
        \pgfmathtruncatemacro{\label}{5-\y};
        \node[circle,draw=black,fill=white!50!black,minimum size=0.25,inner sep=1pt]
        (\label) at (1*\x,0.25*\y) {};
}
}
\foreach \j in {1,...,2}
{
\foreach \i in {1,...,2}
{
        \pgfmathtruncatemacro{\y}{(\i - 4)};
        \pgfmathtruncatemacro{\x}{6+\j}; 
        \pgfmathtruncatemacro{\label}{3+\y};
        \ifthenelse{\label=0}{
        \node[circle,thick,draw=black,fill=white!40!red,minimum size=22]
        (\label) at (1.0*\x,1.0*\y) {$\Ym_{\Ns}$};}{
        \node[circle,thick,draw=black,fill=white!40!red,minimum size=22]
        (\label) at (1.0*\x,1.0*\y) {$\Ym_{\Ns-\label}$};
        }
}
}
\node[circle,thick,draw=black,fill=white!60!green,minimum size=22]
        () at (7,0) {$\Ym_{3}$};
\node[circle,thick,draw=black,fill=white!60!green,minimum size=22]
        () at (8,2) {$\Ym_{1}$};
\node[circle,thick,draw=black,fill=white!60!green,minimum size=22]
        () at (8,-2) {$\Ym_{\Ns-1}$};
        
\node[circle,thick,draw=black,fill=white!60!green,minimum size=22]
        () at (3,2) {$\Ym_{1}$};
        
\node[circle,thick,draw=black,fill=white!60!green,minimum size=22]
        () at (2,0) {$\Ym_{3}$};
\node[circle,thick,draw=black,fill=white!60!green,minimum size=22]
        () at (1,2) {$\Ym_{1}$};
\node[circle,thick,draw=black,fill=white!60!green,minimum size=22]
        () at (4,1) {$\Ym_{2}$};
\node[circle,thick,draw=black,fill=white!60!green,minimum size=22]
        () at (3,0) {$\Ym_{3}$};

\foreach \j in {1,...,4}
{
\foreach \i in {1,...,2}
{
        \pgfmathtruncatemacro{\y}{(\i - 4)};
        \pgfmathtruncatemacro{\x}{\j}; 
        \pgfmathtruncatemacro{\label}{3+\y};
        \ifthenelse{\label=0}{
        \node[circle,thick,draw=black,fill=white!40!red,minimum size=22]
        (\label) at (1.0*\x,1.0*\y) {$\Ym_{\Ns}$};}{
        \node[circle,thick,draw=black,fill=white!40!red,minimum size=22]
        (\label) at (1.0*\x,1.0*\y) {$\Ym_{\Ns-\label}$};
        }
}
}

\foreach \j in {1,...,4}
{
        \pgfmathtruncatemacro{\y}{( -6)};
        \pgfmathtruncatemacro{\x}{\j}; 
        \pgfmathtruncatemacro{\optiter}{\x};

        \node[thick,draw=black,fill=white!80!blue,minimum size=15]
        (\optiter) at (1.0*\x,1.0*\y) {$\ppm_{\optiter}$};
        
        \draw[-latex, line width=1mm] (\x,\y+2.5) -- (\x,\y+0.5);
        \draw[-latex,thick] (\x+0.25,\y) -- (\x+0.75,\y);
}

\foreach \j in {1,...,2}
{
        \pgfmathtruncatemacro{\y}{( - 6)};
        \pgfmathtruncatemacro{\x}{6+\j}; 
        \pgfmathtruncatemacro{\optiter}{8-\x};
        \ifthenelse{\optiter=0}{
        \node[thick,draw=black,fill=white!80!blue,minimum size=15]
        (\optiter) at (1.0*\x,1.0*\y) {$\ppm_{M_b}$};}{
        \node[thick,draw=black,fill=white!80!blue,minimum size=15]
        (\optiter) at (1.0*\x,1.0*\y) {$\ppm_{M_b-\optiter}$};
        }
        \draw[-latex, line width=1mm] (\x,\y+2.5) -- (\x,\y+0.5);
        \draw[-latex,thick] (\x-0.675,\y) -- (\x-0.3,\y);
}

\node[circle,thick,draw=black,fill=white!60!green,minimum size=22]
        () at (1,-3) {$\Ym_{\Ns}$};
\node[circle,thick,draw=black,fill=white!60!green,minimum size=22]
        () at (4,-2) {$\Ym_{\Ns-1}$};
\node[circle,thick,draw=black,fill=white!60!green,minimum size=22]
        () at (2,-2) {$\Ym_{\Ns-1}$};
\node[circle,thick,draw=black,fill=white!60!green,minimum size=22]
        () at (10,-1) {$\Ym_{i}$};
\node[circle,thick,draw=black,fill=white!40!red,minimum size=22]
        () at (10,0) {$\Ym_{i}$};
\node[draw=none,right]
        () at (10.5,0) {: Not used};
\node[draw=none,right]
        () at (10.5,-1) {: Used};
    \end{tikzpicture}
    		\caption{Schematic of the proposed approach} \label{fig:proposedschematic}
	\end{subfigure}
	
    \caption{Comparison of the proposed approach for TOuU with a typical Monte Carlo approach utilizing $\Ns$ random samples for $M_b$ and $M_f$ optimization iterations, respectively. Note that, though total number of random samples $\Ns$ is finite in this schematic in reality it can be infinite.}
    \label{fig:schem}
\end{figure}
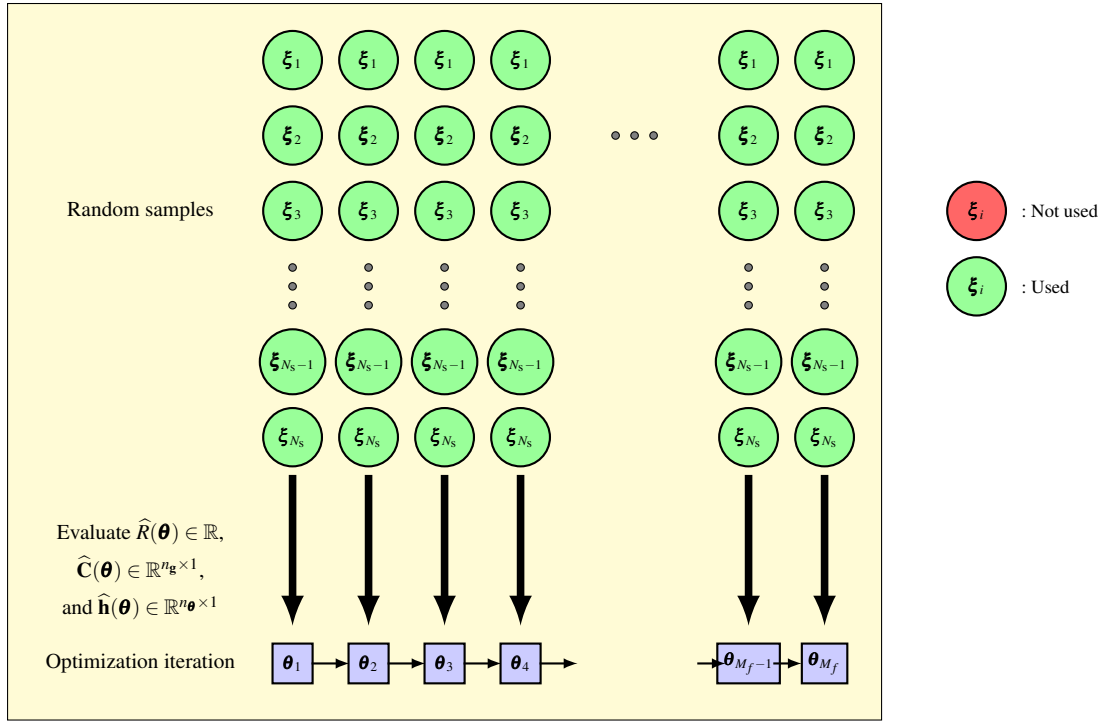
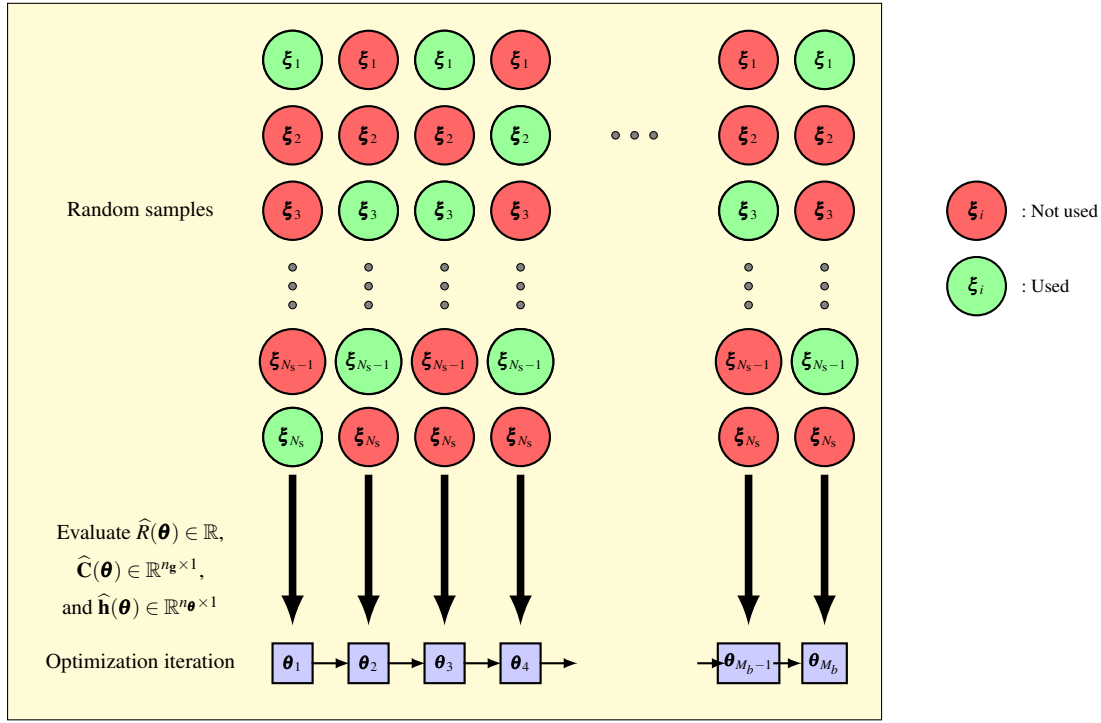

In this paper, the objective, constraints and their gradients are estimated at every design iteration with a small number of random samples $\{\Ym_i\}_{i=1}^n$ instead of using all $\Ns$ samples with $n\ll \Ns$ resulting in a significant reduction in computational cost. 
These quantities thus calculated, however, are stochastic in nature. Hence, the gradients estimated in \eqref{eq:Ghatgraddef} and \eqref{eq:EdC} using $n$ random samples can be thought of as \textit{ stochastic gradients} calculated using a small batch of random samples \citep{ruder2016overview,bottou2018optimization}. 
%
%
Note that, we compute only one mean gradient vector $\widehat{\hb}(\ppm;\kappaa;\lambda)$ at every iteration and not multiple small batches or \textit{mini-batches}. 
%
%
Also, a single random sample $\Ym$ can be used for mean design with variance weighting $\lambda = 0$ but not with $\lambda \neq 0$ as more than one sample is required to estimate the variance. For brevity, we omit the dependence of $\hb$ on $\kappaa$ and $\lambda$. 
However, to successfully implement the stochastic gradient methods for TOuU we need to use a different set of random samples for each iteration. 
The schematic of this approach is illustrated in Figure \ref{fig:schem} and is compared with a standard Monte Carlo approach with $\Ns$ gradients. 

In this proposed approach, the  samples in each mini-batch are statistically independent, so we analyze the samples in parallel, \ie~we compute the objective and constraint measures as well as their gradients in parallel.
The parallel implementation is straight forward and constitutes a trivial modification to (\ref{eq:Hhatgrad}). In the next two subsections, we discuss the SGD method, its variants, and GCMMA in brief, where we use these stochastic gradients.

\FloatBarrier

\subsection{Stochastic Gradient Descent (SGD) Method and its Variants} 
\label{sec:sgd}
The SGD method uses a single realization of $\Ym$ to update the design at the $k$th iteration utilizing the gradient $\hb_k$ \citep{bottou2018optimization} that is defined by
\begin{equation}\label{eq:sgd}
\begin{split}
&\hb_k=\nabla{f}(\ppm_k,\Ym_i) + \tpos{\kappaa}{\nabla{\Gm}}(\ppm_k,\Ym_i);\\
&\ppm_{k+1}=\ppm_k-\eta \hb_{k},\\
\end{split}
\end{equation}
where $\eta $ is the step size, also known as the \textit{learning rate}. 
Hence, the computational cost of stochastic gradient descent is very low. However, the convergence of the basic SGD method can be very slow since we might not follow the descent direction at every iteration. Instead, the descent is achieved in expectation as the expectation of the stochastic gradient is same as the gradient of the objective in \eqref{eq:exp_risk}. Note that, in literature, the use of a positive hyperparameter term $\kappaa$ with SGD and its variants was mainly to introduce regularization \citep{tsuruoka2009stochastic,lavergne2010practical,collins2014memory}. In this study, we propose the use of $\kappaa$ to enforce the design constraints. 
%

A straightforward extension of the basic SGD method is to use a small batch of random samples to estimate the gradient $\hb_k$. This version is known as \textit{mini-batch gradient descent} \citep{ruder2016overview,bottou2018optimization}. 
\begin{algorithm}[htb]
	\begin{algorithmic}
		\STATE Given $\eta$.
		\STATE Initialize $\ppm_1$.
		\FOR {$k=1,2,\dots$}
		\STATE Compute $\hb_k:=\hb(\ppm_k)$.
		\STATE Set $\ppm_{k+1} \leftarrow \ppm_{k} - \eta\hb_k$. {[see Eqn. (\ref{eq:sgd})]}
		\ENDFOR
	\end{algorithmic}
	\caption{\textit{Basic stochastic gradient descent method.} \citep{bottou2018optimization}}
	\label{alg:sgd}
\end{algorithm}
In the past few years, different modifications to the above algorithm have been proposed to improve the convergence of SGD methods. Five such algorithms that we utilize here are presented next.


\subsubsection{\texorpdfstring{Adaptive Subgradient Methods (AdaGrad)}{Stochastic Gradient Descent: Adagrad}}
\label{subsec:adagrad}

%
The first algorithm that we present is AdaGrad~\citep{duchi2011adagrad}, where historical information about the gradients is used to modify the update of optimization variables. It dampens the movements along directions with historically large gradients, thus adapting the learning rate and facilitating a faster convergence.
%
%
In this algorithm, at iteration $k$, we compute the following auxiliary variable
\begin{equation}
\am_{k,j}=\sum_{i=1}^k{\hb^2_{i,j}},\qquad j=1,\dots,\np,
\end{equation}
%
which is then used in the following update rule
\begin{equation}\label{eq:adagrad}
\ppm_{k+1}=\ppm_k-\eta\am_k^{-1/2}\hb_k,
\end{equation}
where the vector multiplication is performed component-wise. To avoid division by zero, a small number $\epsilon$ ($=10^{-8}$ used herein) is incorporated in the denominator of the update as shown in Algorithm \ref{alg:adagrad}.
\begin{algorithm}[htb]
	\begin{algorithmic}
		\STATE Given $\eta$.
		\STATE Initialize $\ppm_1$.
		\STATE Initialize $\am := \bm{0}$, having the same dimensions as $\ppm_1$.
		\FOR {$k=1,2,\dots,$}
		\STATE Compute $\hb_k := \hb(\ppm_{k})$.
		\STATE Set $\am_j := \am_j + \hb_{k,j}^2,\qquad j=1,\dots,\np$.
		\STATE Set $\ppm_{k+1,j} \leftarrow \ppm_{k,j} - \eta\frac{\hb(\ppm_{k,j})}{\sqrt{\am_j}+\sqrt{\epsilon}},\qquad j=1,\dots,n_p$. {[see Eqn. (\ref{eq:adagrad})]}
		\ENDFOR
	\end{algorithmic}
	\caption{\textit{AdaGrad} \citep{duchi2011adagrad}}
	\label{alg:adagrad}
\end{algorithm}


\subsubsection{\texorpdfstring{Adadelta}{Adadelta}}
\label{subsec:adadelta}
The second algorithm that we study in this paper is Adadelta~\citep{zeiler2012adadelta}, which is similar to AdaGrad but reduces the window of historical gradients using an exponential decay rate $\zeta$, \textit{i.e.}, at the $k$th iteration the gradient accumulation is estimated using
\begin{equation}\label{eq:gradaccu}
\am_{\hb,j} =  \zeta \hb_{k-1,j}^2 + (1-\zeta) \hb_{k,j}^2, \qquad j=1,\dots,\np.
\end{equation}
\cite{zeiler2012adadelta} suggested a value of 0.95 for $\zeta$ that we use here as well. A root-mean-square (RMS) of the gradient history accumulator is then defined as
\begin{equation}\label{eq:rmsgrad}
{\mathrm{RMS}[\hb_{k,j}]=\sqrt{\am_{\hb,j}+\epsilon}.}
\end{equation}
Again, a very small number $\epsilon$ ($=10^{-8}$ used herein) is used to avoid division by zero. To be consistent with the units an accumulation of parameter update history is used,
\begin{equation}\label{eq:paramaccu}
\am_{\ppm,j} = \zeta \am_{\ppm,j} + (1-\zeta) \Delta\ppm_{k,j}^2, \qquad j=1,\dots,\np,
\end{equation}
where $\Delta\ppm_k$ is the update of the parameter vector $\ppm$ at iteration $k$. The RMS of the parameter update history becomes
\begin{equation}\label{eq:rmsparam}
{\mathrm{RMS}[\Delta\ppm_{k,j}]=\sqrt{\am_{\ppm,j}+\epsilon}.}
\end{equation}
The parameter update is then defined as
\begin{equation}\label{eq:adadelta}
\begin{split}
    \Delta\ppm_k &= \frac{\mathrm{RMS[\Delta \ppm_{k}]}}{\mathrm{RMS[\hb_k]}}\hb_k,\\
\ppm_{k+1} &= \ppm_k - \Delta\ppm_k. 
\end{split}
\end{equation}
The main steps of the Adadelta method are summarized in Algorithm \ref{alg:adadelta}. Note that, in this algorithm, we do not need the learning rate $\eta$ explicitly.
\begin{algorithm}[htb]
	\begin{algorithmic}
		\STATE Given $\rho$, $\epsilon$.
		\STATE Initialize $\ppm_1$.
		\STATE Initialize $\am_{\hb} := \bm{0}$ and $\am_{\ppm} := \bm{0}$, having the same dimensions as $\bm{p}_1$.
		\FOR {$k=1,2,\dots,$}
		\STATE Compute $\hb_k := \hb(\ppm_{k})$.
		\STATE Set $\am_{\hb,j} \leftarrow \zeta\am_{\hb,j} + (1-\zeta)\hb_{k,j}^2,\qquad j=1,\dots,\np$. {[see Eqn. (\ref{eq:gradaccu})]}
		\STATE  Evaluate $\mathrm{RMS}[\hb_{k,j}] = \sqrt{\am_{\hb_{k,j}}+\epsilon},\qquad j=1,2,\dots,\np$. {[see Eqn. (\ref{eq:rmsgrad})]}
		\STATE  Evaluate $\mathrm{RMS}[\Delta\ppm_{k,j}] = \sqrt{\am_{\ppm_{k,j}}+\epsilon},\qquad j=1,2,\dots,\np$. {[see Eqn. (\ref{eq:rmsparam})]}
        \STATE Set $\Delta \ppm_{k,j} = \frac{\mathrm{RMS[\Delta \ppm_{k,j}]}}{\mathrm{RMS[\hb_{k,j}]}}\hb_{k,j} \qquad j = 1, 2, \dots,\np$. 
		\STATE Set $\ppm_{k+1,j} \leftarrow \ppm_{k,j} - \Delta \ppm_{k,j} \qquad j = 1, 2, \dots,\np$. {[see Eqn. (\ref{eq:adadelta})]}
        \STATE Set $\am_{\ppm,j} \leftarrow \zeta\am_{\ppm,j} + (1-\zeta)(\Delta \ppm_{k,j})^2,\qquad j=1,\dots,\np$. {[see Eqn. (\ref{eq:paramaccu})]}
		\ENDFOR
	\end{algorithmic}
	\caption{\textit{Adadelta} \citep{zeiler2012adadelta}}
	\label{alg:adadelta}
\end{algorithm}


\subsubsection{Adaptive Moment Estimation (Adam)}
\label{subsec:Adam}

Another algorithm that we consider in this paper for TO is Adam~\citep{kingma2014adam}. This algorithm is similar to Adagrad, but it intends to additionally smooth the variability in $\hb$ from iteration to iteration through the accumulation of historical gradient and squared gradient information. Two exponential decay rates $\beta_m$ and $\beta_v$ are used as follows
\begin{equation}\label{eq:adamdecay}
\begin{split}
\mm_k &= \beta_m \mm_{k-1} + (1-\beta_m)\hb_k;\\
\vm_{k,j} &= \beta_v \vm_{k-1,j} + (1-\beta_v)\hb_{k,j}^2,\qquad j=1,2,\dots,\np.\\
\end{split}
\end{equation}
The authors of Adam~\citep{kingma2014adam} recommend $\beta_m = 0.9$ and $\beta_v = 0.999$, which we also use here. The following initialization bias correction is applied to the above-accumulated quantities
\begin{equation}\label{eq:adambias}
\begin{split}
\widehat{\mm}_k &= \frac{\mm_k}{1-\beta_m^k};\\
\widehat{\vm}_k &= \frac{\vm_k}{1-\beta_v^k}.\\
\end{split}
\end{equation}
Finally, the parameters are updated using
\begin{equation}\label{eq:adam}
\ppm_{k+1,j} = \ppm_{k,j} - \eta\frac{\widehat{\mm}_j}{\sqrt{\widehat{\vm}_j}+{\epsilon}}\qquad j=1,2,\dots,\np.
\end{equation}
The main steps of the Adam algorithm are summarized in Algorithm \ref{alg:adam}.
\begin{algorithm}[htb]
	\begin{algorithmic}
		\STATE Given $\eta$, $\beta_m$, $\beta_t$, and $\epsilon$.
		\STATE Initialize $\ppm_1$.
		\STATE Initialize $\mm = \bm{0}$.
		\STATE Initialize $\vm = \bm{0}$.
		\FOR {$k=1,2,\dots,$}
		\STATE Compute $\hb_k := \hb(\ppm_{k})$.
		\STATE Set $\mm \leftarrow \beta_m\mm + (1-\beta_m)\hb_k$. {[see Eqn. (\ref{eq:adamdecay})]}
		\STATE Set $\vm_j \leftarrow \beta_v\vm_j + (1-\beta_v)\hb_{k,j}^2\qquad j=1,2,\dots,\np$. {[see Eqn. (\ref{eq:adamdecay})]}
		\STATE Set $\widehat\mm \leftarrow \mm/(1-\beta_m^k)$. {[see Eqn. (\ref{eq:adambias})]}
		\STATE Set $\widehat\vm \leftarrow \vm/(1-\beta_v^k)$. {[see Eqn. (\ref{eq:adambias})]}
		\STATE Set $\ppm_{k+1,j} \leftarrow \ppm_{k,j} - \eta\frac{\widehat{\mm}_j}{\sqrt{\widehat{\vm}_j}+{\epsilon}}\qquad j=1,2,\dots,\np$. {[see Eqn. (\ref{eq:adam})]}
		\ENDFOR
	\end{algorithmic}
	\caption{\textit{Adam} \citep{kingma2014adam}}
	\label{alg:adam}
\end{algorithm}


\subsubsection{\texorpdfstring{Stochastic Average Gradient (SAG)}{Stochastic Gradient Descent: SAG}}
\label{subsec:SAG}

The next algorithm that we evaluate in this paper is the Stochastic Average Gradient or SAG~\citep{rouxetalSAG2012}, which updates the gradient information for one random sample at every iteration. The parameters are updated using
\begin{equation}\label{eq:sag}
\begin{split}
\ppm_{k+1}&= \ppm_k - \frac{\eta}{\Ns}\sum_{i=1}^{\Ns}\dm_{k,i};\\
\dm_{k,i} &= \begin{cases}
\hb(\ppm_k;\Ym_{i})\qquad \text{if $i=t\in \{1,2,\dots,{\Ns}\}$};\\
\dm_{k-1,i}\qquad\qquad \text{otherwise},\\
\end{cases}
\end{split}
\end{equation}
where $t$ is selected randomly from $\{1,2,\dots,\Ns\}$. This method minimizes a batch of $\Ns$ input uncertainty samples that is necessarily finite, as opposed to the typically infinite potential input uncertainty samples that are assumed by the random variable model. 
\begin{algorithm}[htb]
	\begin{algorithmic}
		\STATE Given $\eta$.
		\STATE Initialize $\ppm_1$.
		\STATE Initialize $\dm = \bm{0}$.
		\FOR {$k=1,2,\dots,$}
		\STATE Draw $t$ randomly from $\{1,\cdots,\Ns\}$.
		\IF{$i=t$}
		\STATE Compute $\dm_{i} := {\hb}(\ppm_{k};\Ym_{i})$.
		\ENDIF
		\STATE $\ppm_{k+1} \leftarrow \ppm_k - \frac{\eta}{\Ns}\sum_{i=1}^{\Ns}\dm_i$. [see Eqn. (\ref{eq:sag})]
		\ENDFOR
	\end{algorithmic}
	\caption{\textit{SAG} \citep{rouxetalSAG2012}}
	\label{alg:sag}
\end{algorithm}


\subsubsection{Stochastic Variance Reduced Gradient (SVRG)}

The final algorithm that we consider in this work is the stochastic variance reduced gradient (SVRG) method \citep{johnson2013accelerating}. In this algorithm, a variance reduction method is introduced in SGD by maintaining a  parameter estimate ${\ppm}_{\mathrm{best}}$ that is updated only at every $m$ iterations. Using this parameter estimate and the following estimate of the gradient
\begin{equation}
\hb_\mathrm{best} = \frac{1}{\Ns}\sum_{i=1}^{\Ns} \hb(\ppm_\mathrm{best},\Ym_i),
\end{equation}
an update rule similar to the SGD is applied. However, the gradient term is replaced with 
\begin{equation}\label{eq:svrg}
\hb_k=\hb(\ppm_k,\Ym_{t})-\hb(\ppm_\mathrm{best},\Ym_{t})+\hb_\mathrm{best}
\end{equation}
for randomly chosen $t\in\{1,2,\dots,\Ns\}$, \textit{i.e.}, $\hb(\ppm_\mathrm{best})$ is used here as a control variate \citep{ross2013simulation}. These steps are illustrated in Algorithm \ref{alg:svrg}.
\begin{algorithm}[htb]
	\begin{algorithmic}
		\STATE Given $\eta$ and $m$.
		\STATE Initialize $\widetilde\ppm_1$.
		\FOR {$j=1,2,\dots,$}
		\STATE Set ${\ppm}_{\mathrm{best}} = \widetilde{\ppm}_{j}$.
		\STATE Set $\hb_\mathrm{best} = \frac{1}{\Ns}\sum_{i=1}^{\Ns} \hb(\ppm_\mathrm{best},\Ym_i)$.
		\STATE Set $\ppm_1 = \ppm_{\mathrm{best}}$.
		\FOR{$k=1,2,\dots,m$}
		\STATE Randomly choose $t \in \{1,2,\dots,\Ns\}$.
		\STATE Set $\ppm_{k+1}\rightarrow\ppm_k-\eta \big[\hb(\ppm_k,\Ym_{t})-\hb(\ppm_\mathrm{best},\Ym_{t})+\hb_\mathrm{best}\big]$. [see Eqn. (\ref{eq:svrg})]
		\ENDFOR
		\STATE $\widetilde\ppm_{j} \leftarrow \ppm_{m+1}$ or $\widetilde\ppm_{j} \leftarrow \ppm_{t}$ for randomly chosen $t \in \{1,2,\dots,m\}$.
		\ENDFOR
	\end{algorithmic}
	\caption{\textit{SVRG} \citep{johnson2013accelerating}}
	\label{alg:svrg}
\end{algorithm}

We note that most of the SGD methods considered in this study depend on the learning rate parameter $\eta$, which controls the size of the step along the modified gradients. The tuning of this parameter is crucial to the success of these methods. For tuning, either a learning schedule can be followed or some pilot runs are needed. Since the overall cost of one iteration of the optimization is not significant given only a few samples are used at every iteration, these pilot runs do not add substantial cost to the overall optimization procedure. This aspect is further illustrated in Section \ref{sec:mbb}. Note that a line search in the descent direction of the current iteration can also be performed removing the need for tuning of the learning rate parameter $\eta$ \citep{mahsereci2015probabilistic}; however, this is beyond the scope of the current paper.

\subsection{Globally Convergent Method of Moving Asymptotes (GCMMA) with Stochastic Gradients}\label{sec:gcmma}

Owing to its popularity in TO, we also study the performance of Globally Convergent Method of Moving Asymptotes (GCMMA) using stochastic gradients. In GCMMA, at every iteration, a convex subproblem is constructed using conservative approximations of the objective and constraint functions around the current optimization variable values. The solution of the subproblem is determined by solving its Karush-Kuhn-Tucker (KKT) conditions. GCMMA is well suited for large scale optimization problems as the subproblem is separable, \ie~ it can be decomposed in $n_\theta$ single-variable problems. In this paper, we supply GCMMA with stochastic gradients as outlined in the next section. 

\FloatBarrier
\section{Numerical Examples}
\label{sec:examples}

In this section, we investigate the characteristics and performance of using stochastic gradients with SGD methods and the GCMMA for TOuU by two numerical examples. In both examples, we seek to minimize the expected compliance of the structure subject to a volume constraint. The structural response is described by a static linear elastic model. First, a 2D design problem is studied with the explicit level set approach of Section \ref{sec:des_model}, using a primitive geometry parameterization of the design. We consider uncertainty in shape and topology. Second, we apply the SIMP method to a 3D problem with uncertainty in loading conditions and support stiffness. These examples have been chosen because they feature different influences of optimization variables and uncertain parameters on the objective and constraint. Since the TO problems considered here are non-convex, there may be several local minima. Therefore, different optimization algorithms may converge to different designs. To compare several algorithms in this work, we consider the value of the objective at the optimized solution to assess the performance of these methods, independently of the particular design obtained. 
 
For the level set method, each realization of a design candidate is analyzed by the XFEM using bi-linear elements in 2D. To suppress rigid body motion of structural members, we apply soft distributed fictitious springs to the solid domain; see \cite{villanueva2014density}. 
For the SIMP example, the structural response is predicted by a standard finite element analysis with  and tri-linear elements in 3D. 
Further details of the XFEM formulation and sensitivity analysis used in these examples are outlined in Section \ref{sec:ana_model}. The gradients of the performance and constraint measures are computed by the adjoint method. The linear systems of the forward and sensitivity analyses are solved by a direct serial solver for the 2D problems and by an iterative parallel solver for the 3D problems. The algorithmic parameters controlling the behavior of the different optimization algorithms studied in the following are given with each example.


\subsection{Example I: Beam Design Problem}
\label{sec:mbb}

\begin{figure}[htb!]
	\centering
	\begin{tikzpicture}[scale=0.75,every node/.style={minimum size=1cm},on grid]
	\draw [ultra thick,fill=blue!20,draw=none] (0,-3.5) rectangle (10,0);
	\draw [ultra thick] (0,0) -- (10,0);
	\draw [ultra thick] (10,0) -- (10,-3.5);
	\draw [ultra thick] (10,-3.5) -- (0,-3.5);
	\draw [ultra thick,dashed] (0,-3.5) -- (0,0);
	
	\draw[ultra thick,fill=gray!30] (9.85,-4.15) circle (0.125);
	\draw[ultra thick,fill=gray!30] (10.15,-4.15) circle (0.125);
	\draw[fill=gray!30]    (10,-3.5) -- ++(0.3,-0.5) -- ++(-0.6,0) -- ++(0.3,0.5);
	\draw [ultra thick] (10,-3.5) -- (10.3,-4);
	\draw [ultra thick] (10.3,-4) -- (9.7,-4);
	\draw [ultra thick] (9.7,-4) -- (10,-3.5);
	
	\node[draw=none] at (5, -1.5)  (c)     {\Large Minimize compliance};
	\node[draw=none] at (5, -2)  (c)     {\Large subject to mass constraint};
	
	\draw[-latex, line width=1mm] (0,1.25) -- (0,0);
	\node[draw = none] at (0,1.5) () {\Large $P$};
	
	\draw[ultra thick,fill=gray!30] (-0.65,-0.15) circle (0.125);
	\draw[ultra thick,fill=gray!30] (-0.65,0.15) circle (0.125);
	\draw[fill=gray!30]    (0,0) -- ++(-0.5,0.3) -- ++(0,-0.6) -- ++(0.5,0.3);
	\draw [ultra thick] (0,0) -- (-0.5,-0.3);
	\draw [ultra thick] (-0.5,-0.3) -- (-0.5,0.3);
	\draw [ultra thick] (-0.5,0.3) -- (0,0);
	\draw[ultra thick,fill=gray!30] (-0.65,-3.35) circle (0.125);
	\draw[ultra thick,fill=gray!30] (-0.65,-3.65) circle (0.125);
	\draw[fill=gray!30]    (0,-3.5) -- ++(-0.5,0.3) -- ++(0,-0.6) -- ++(0.5,0.3);
	\draw [ultra thick] (0,-3.5) -- (-0.5,-3.8);
	\draw [ultra thick] (-0.5,-3.8) -- (-0.5,-3.2);
	\draw [ultra thick] (-0.5,-3.2) -- (0,-3.5);
	
	\draw[ultra thick,fill=gray!30] (-0.65,-1.90) circle (0.125);
	\draw[ultra thick,fill=gray!30] (-0.65,-1.60) circle (0.125);
	\draw[fill=gray!30]    (0,-1.75) -- ++(-0.5,0.3) -- ++(0,-0.6) -- ++(0.5,0.3);
	\draw [ultra thick] (0,-1.75) -- (-0.5,-2.05);
	\draw [ultra thick] (-0.5,-2.05) -- (-0.5,-1.45);
	\draw [ultra thick] (-0.5,-1.45) -- (0,-1.75);
	\end{tikzpicture}
	\caption{Schematic for the beam design problem (design domain is shown as the shaded region).}
	\label{fig:prims_schematic}
\end{figure}
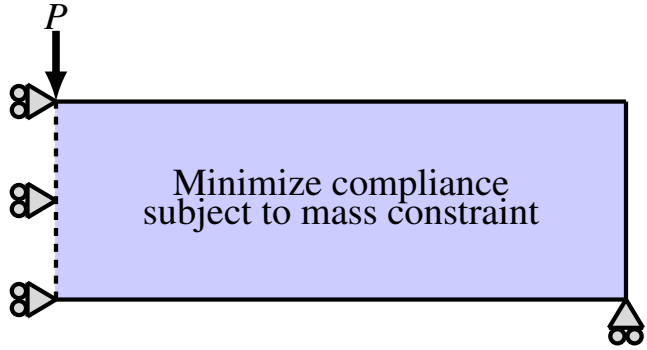

{In this problem, we are interested in optimizing the geometry of a 2D beam with a point load $2P$ placed at the center of it span as shown in \ref{fig:prims_schematic}. The design domain is formed by a $6.0 \times 1.0$ rectangle. The beam is simply supported on the bottom left and right corners. By considering symmetry, we restrict our computational domain to one-half the design domain's length, where the load is placed at the top left corner, and the structure is supported at the bottom right corner.} The objective function of the optimization is the compliance of the system, and the constraint is to ensure the mass-ratio of the design area is no more than 50\% of the maximum design mass (when entire design domain is filled with the material). 
%

We describe the geometry of the structure by the explicit level set method with a geometric primitive parameterization of the level set function.  We restrict the design space to structures that can be assembled with twelve rectangular bars considering one half of the design domain. Each of these twelve bars is assumed to have five geometric parameters, which are allowed to vary; see Section \ref{sec:des_model}. The total number of optimization variables is $\np=60$.

We introduce geometric uncertainties in this problem by assuming that the geometry parameters of the {\it manufactured} bars denoted by $a_r,b_r,x_r^c,y_r^c,\alpha_r, \ r=1, \ldots, 12$, are uncertain variations of the corresponding design geometry parameters $\bm\theta_j$, ${j=5(r-1)}$. We subject the design geometry parameters to uniform perturbations as follows,
%
%
\begin{equation}
    \begin{split}
        x_r^c    &= \theta_{j+1} \ +0.005\xi_{j+1};\\
        y_r^c    &= \theta_{j+2} \ +0.005\xi_{j+2};\\
        a_r      &= \theta_{j+3} \ +0.005\xi_{j+3};\\
        b_r      &= \theta_{j+4} \ +0.005\xi_{j+4};\\
        \alpha_r &= \theta_{j+5} \ +0.005\xi_{j+5};\\
    \end{split}
    \quad \textrm{with} \ r=1, \ldots, 12 \ ; \ {j= 5 \ (r-1)}
\end{equation}
where $\xi_i$'s are independent uniform random variables on $[-1,1]$. The lower and upper bounds of the optimization variables $\theta_i$ are given in Table \ref{tab:barparameters}. {These box constraints are accounted for by the optimization algorithms, \ie~either GCMMA or SGDs. Specifically, with different variants of SGD, if the perturbation of an optimization variable $\theta_i$ via $\xi_i${, for $i=1,\dots,60$,}  places the primitive bar outside of the shaded region in Figure \ref{fig:prims_schematic},
%
%
then that variable is simply replaced by its maximum or minimum allowed value.} 
 

%
%
\begin{table}[]
    \centering
    \begin{tabular}{c|c|c| c}
    \hline 
         Bar parameter & Optimization Variable & Lower Bound & Upper Bound \\\hline 
         $x_r^c $      & $\theta_{j+1}$            & $0.0$       & $3.0$       \\ 
         $y_r^c $      & $\theta_{j+2}$        & $0.0$       & $1.0$       \\ 
         $a_r$         & $\theta_{j+3}$        & $0.0$       & $3.15$      \\ 
         $b_r$         & $\theta_{j+4}$        & $0.0$       & $0.1$       \\
         $\alpha_r$    & $\theta_{j+5}$        & $0.0$       & $\pi$       \\ \hline
    \end{tabular}
    \caption{Lower and upper bounds of bar parameters and optimization variables for $r=1, \ldots, 12$ and {$j=5 \ (r-1)$}.}
    \label{tab:barparameters}
\end{table}

This uncertainty model mimics a scenario where, due to manufacturing imperfections, the built structure's dimensions may be slightly different from their designed values. Hence, in this example, the optimization and uncertainty parameters are closely connected. Note that, this problem is high-dimensional, \ie~ $\ny = 60$, and methods based on PCE or sparse grids will not be efficient. We use Young's modulus $E_0$, the material density $\rho_0$, and the load as unity. The Poisson's ratio is assumed as 0.3. To avoid rigid body motion of disconnected material fictitious springs with spring constants $10^{-6}$ are used. Finite elements with edge length 0.02 are used to solve the forward problem. We use $\mu = 10$ in \eqref{eq:rect_def} to define the sharpness of the corners of the primitive rectangles and $\beta_\mathrm{KS} = -20$ in the Kreisselmeier-Steinhauser function (see \eqref{eq:ks}).


\begin{figure}[htb!]
\centering
	\begin{subfigure}[t]{\textwidth}
		\centering
        \begin{tikzpicture}
        \fill [gray!40] (0,-2.5) rectangle (8,0.35);
        \fill [black] (0,0) rectangle (8,0.35);
        \fill [black] (0,-2.5) rectangle (8,-2.15);
        \fill [black] (0,-2.5) rectangle (0.35,0.35);
        \fill [black] (7.65,-2.5) rectangle (8,0.35);
		\fill [black] (0,-1.25) rectangle (8,-0.9);  
        \fill [black] (3.825,-2.5) rectangle (4.175,0.35);
        \node[draw=none,white] at (2, 0.175)   (5) {5};
        \node[draw=none,white] at (6, 0.175)   (6) {6};
        \node[draw=none,white] at (0.175,-0.5)   (8) {8};
        \node[draw=none,white] at (4,-0.5)   (10) {10};
        \node[draw=none,white] at (7.775,-0.5)   (12) {12};
        \node[draw=none,white] at (2,-1.1)   (3) {3};
        \node[draw=none,white] at (6,-1.1)   (4) {4};
        \node[draw=none,white] at (0.175,-1.7)   (7) {7};
        \node[draw=none,white] at (4,-1.7)   (9) {9};
        \node[draw=none,white] at (7.775,-1.7)   (11) {11};
        \node[draw=none,white] at (2,-2.3)   (1) {1};
        \node[draw=none,white] at (6,-2.3)   (2) {2};
        \end{tikzpicture}
		\caption{Initial structural design with bar numbers shown} \label{fig:prims_designsa}
	\end{subfigure}%
    \\
	\begin{subfigure}[t]{\textwidth}
		\centering
		\begin{tikzpicture}
		\node[inner sep=0pt] () at (0,0){\includegraphics[scale = 0.295]{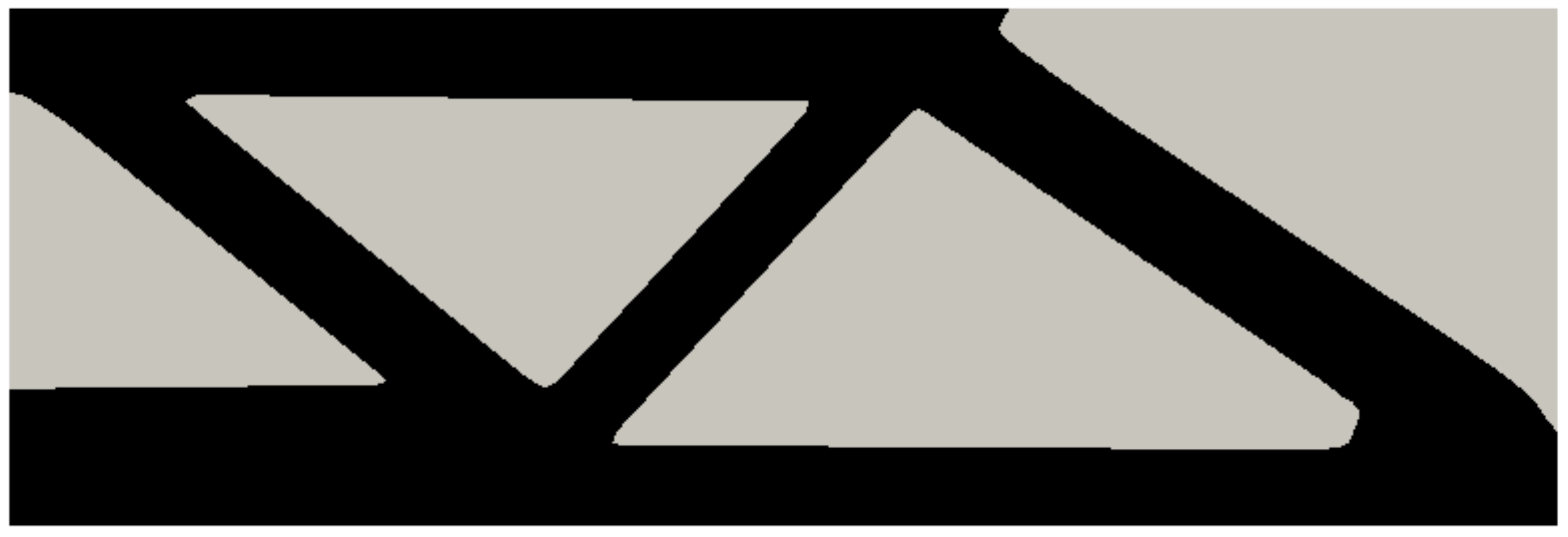}};
		\end{tikzpicture}
		\caption{Final structural design} \label{fig:prims_designsb}
	\end{subfigure}
	\caption{Initial and final design for the beam design problem obtained using GCMMA in a deterministic setting.}
	\label{fig:prims_designs}
\end{figure}
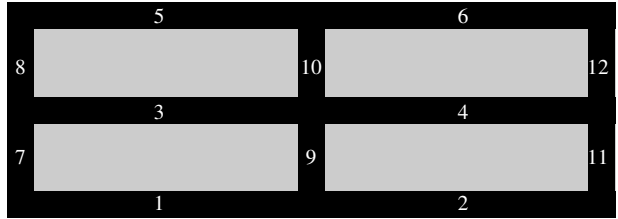
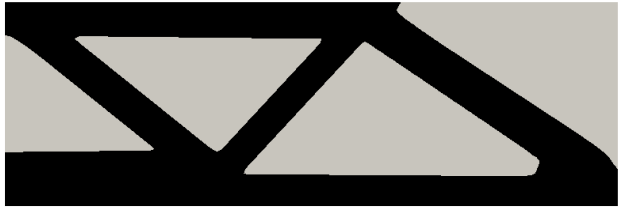
%
%
This design problem is sensitive to the level of randomness in the geometry parameters and the desired robustness, \ie~ the weighting of the variance, $\lambda$, in the formulation of the objective \eqref{eq:exp_risk}. This sensitivity is primarily caused by the fact that due to small design changes and/or stochastic perturbations, either thin structural members disconnect or bars connect in an unfavorable fashion. Either case may lead to large changes in the objective function and dramatically reduced step sizes may be needed to obtain a stable convergence of the optimization process. {Further, to avoid the emergence of thin members we add a regularization term $\lambda_\mathrm{reg}\left(\sum_{r=1}^{12}(a_r^2+b_r^2)\right)$ with $\lambda_\mathrm{reg}=0.001$ to the objective \citep{sharma2017advances,wein2019review}, which encourages small $a_r$ and $b_r$ and large $\phi_r^\mathrm{rec}$ (see \eqref{eq:rect_def}).} 
%
%
%

%
%
Here, we start with an initial design as shown in Figure \ref{fig:prims_designs}. A final design obtained using GCMMA for this problem without the uncertainty is shown in Figure \ref{fig:prims_designsb}. 
For this example, we have chosen $\kappa=1000$ in (\ref{eq:opt_def2}), which is identified by preliminary numerical experimentations to be an appropriate level of $\kappa$.


\subsubsection{Average Design}
\label{subsec:prims_mean_opt}

First, we consider solving the optimization problem (\ref{eq:exp_risk}) and (\ref{eq:glam_def}), omitting the variance term in the objective, \ie~ $\lambda = 0$. We use a sample size of $n=4$ and learning rate $\eta=0.05$ for the algorithms in Section~\ref{sec:method_details} (except for Adadelta and Adam). For Adam, we use $\eta=0.01$. The learning rates are chosen based on some pilot runs. 

\begin{figure}[htb!]
	\centering
	\includegraphics[scale = 0.25]{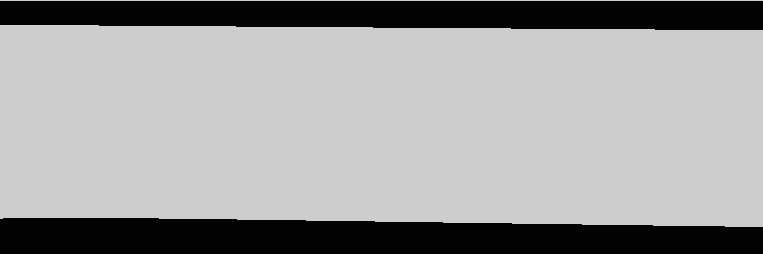}
	\caption{An example of failed design obtained using SVRG with $m=10$ and $\eta=0.05$ (see Algorithm \ref{alg:svrg}). Note that the two structural members are not connected; thus the compliance is infinite.}
	\label{fig:prims_fail}
\end{figure}

To confine the presentation to the most relevant findings, we do not explicitly discuss results from methods  that failed to converge and do not produce any meaningful designs, \textit{e.g.}, the design shown in Figure~\ref{fig:prims_fail}. Specifically, the methods that are not shown in the legends of the plots in this paper have failed. These failures are the result of some designs having disconnected members which leads to poor designs with dramatically larger compliance. Similarly, we exclude methods that do not result in a converged solution with a significantly reduced objective. 
%
%

%
\begin{figure}[htb!]
	\centering
	\begin{subfigure}[t]{\textwidth}
		\centering
		\includegraphics[scale =0.35]{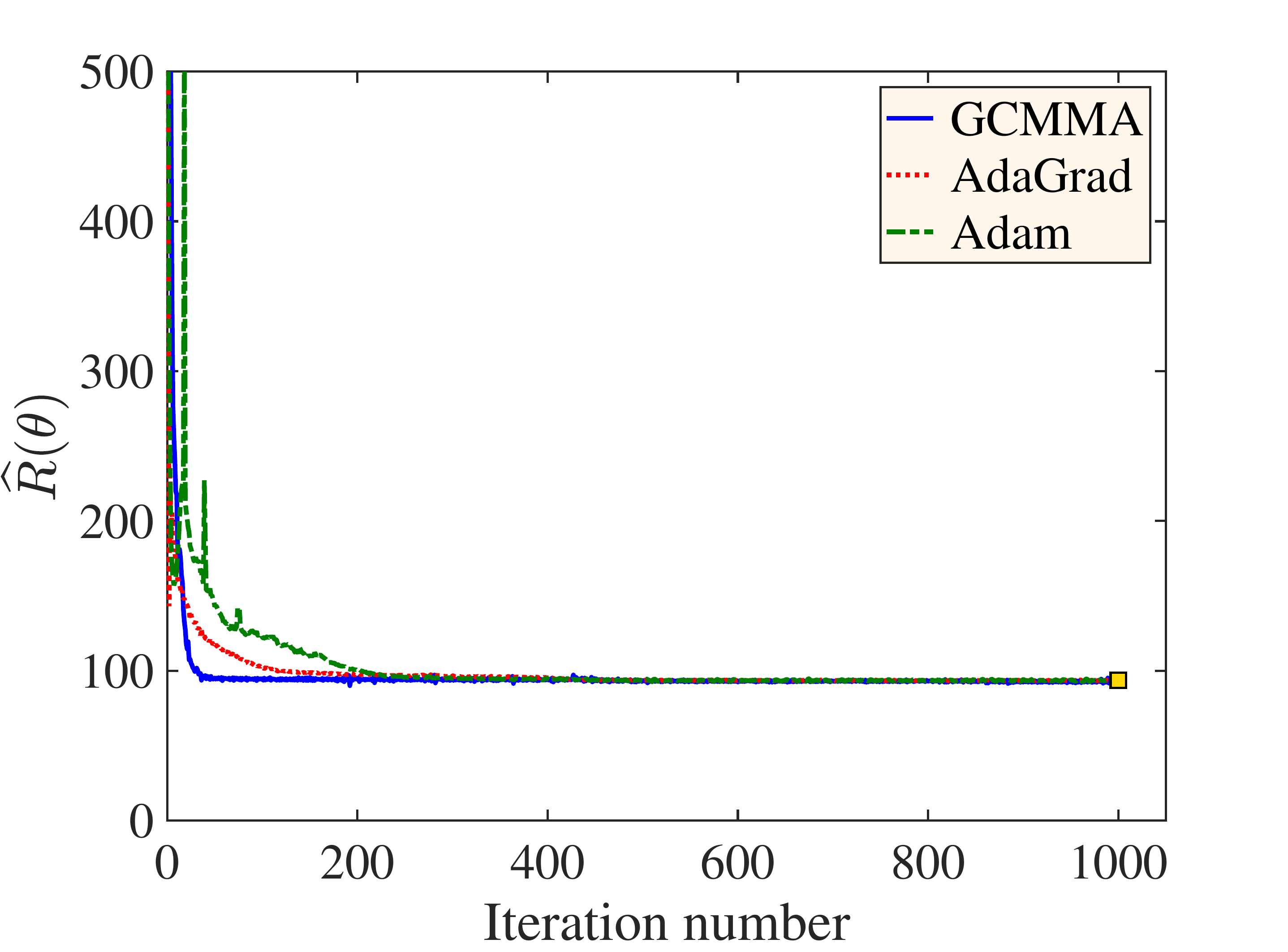}
		\caption{Objective estimates (the square at the right end shows $\widehat{R}(\ppm)$ for the optimal design from Adam evaluated using 1000 random samples)}
        \label{fig:obj_exI}
	\end{subfigure}%
	\\
	\begin{subfigure}[t]{\textwidth}
		\centering
		\includegraphics[scale = 0.35]{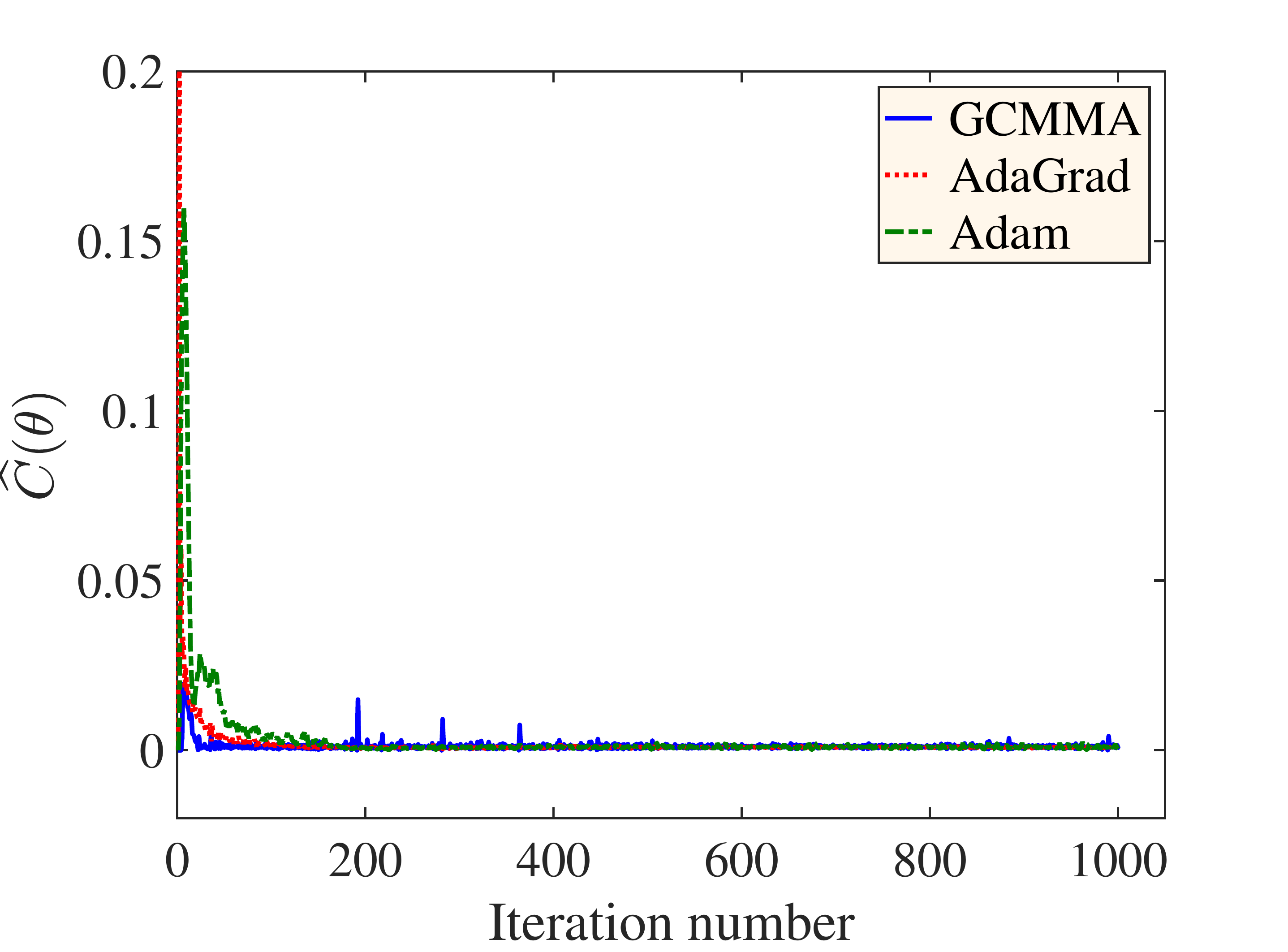}
		\caption{Constraint estimates}
	\end{subfigure}
	\caption{Objectives and constraints for the non-failing methods with $\lambda = 0$ for Example I.}
    \label{fig:con_exI}
	\label{fig:prims_mean_graphs}
\end{figure}

\begin{figure}[htb!]
	\centering
	\begin{subfigure}[t]{0.5\textwidth}
		\centering
		\begin{tikzpicture}
		\node[inner sep=0pt] () at (0,0){\includegraphics[scale =0.18]{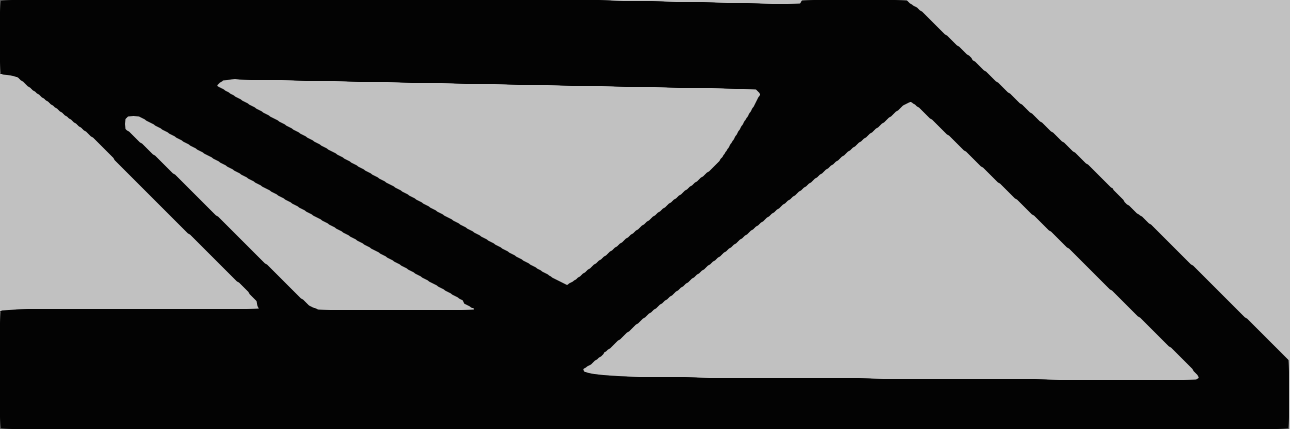}};
        \node[draw=none,white] at (-3, -0.85)   (2) {2};
        \node[draw=none,white] at (0, -0.25)   (9) {9};
        \node[draw=none,white] at (0.8, 0.5)   (10) {10};
        \node[draw=none,white] at (0.2, -1.15)   (1) {1};
        \node[draw=none,white] at (-2, 0.2)   (7) {7};
        \node[draw=none,white] at (-2, 1.05)   (5) {5};
        \node[draw=none,white] at (-3.5, 0.9)   (8) {8};
        \node[draw=none,white] at (-2.85, 0)   (3) {3};
        \node[draw=none,white] at (2.4, 0.25)   (4) {4};
        \node[draw=none,white] at (3.3, -0.5)   (11) {11};
        \node[draw=none,white] at (1, 1)   (6) {6};
		\end{tikzpicture}
		\caption{GCMMA design}
	\end{subfigure}%
	\\
	\begin{subfigure}[t]{0.5\textwidth}
		\centering
		\begin{tikzpicture}
		\node[inner sep=0pt] () at (0,0){\includegraphics[scale = 0.18]{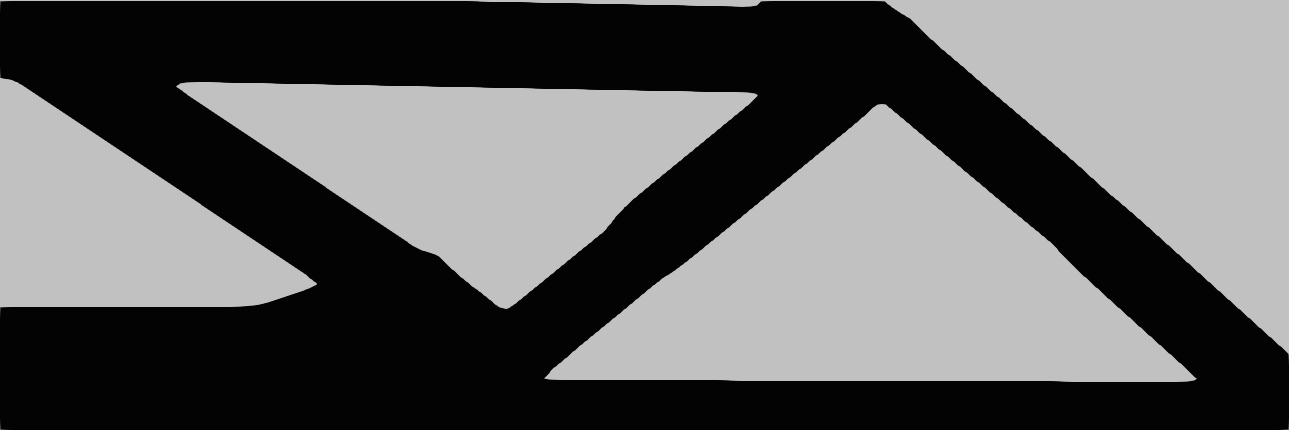}};
		
		\node[draw=none,white] at (-2, 1.05)   (5) {5};
		\node[draw=none,white] at (0.5, -1.2)   (2) {2};
		\node[draw=none,white] at (2.35, 0.25)   (6) {6};
        \node[draw=none,white] at (3.5, -0.75)   (11) {11};
        \node[draw=none,white] at (-0.5, -0.6)   (9) {9};
        \node[draw=none,white] at (0.4, 0.25)   (10) {10};
        \node[draw=none,white] at (-2.75, 0.35)   (8) {8};
        \node[draw=none,white] at (-3.5, -1)   (1) {1};
        \node[draw=none,white] at (-1.9, -0.75)   (7) {7};
		\end{tikzpicture}
		\caption{AdaGrad design}\label{fig:prims_mean_adagrad}
	\end{subfigure}
	\\
	\begin{subfigure}[t]{0.5\textwidth}
		\centering
		\begin{tikzpicture}
		\node[inner sep=0pt] () at (0,0){\includegraphics[scale = 0.18]{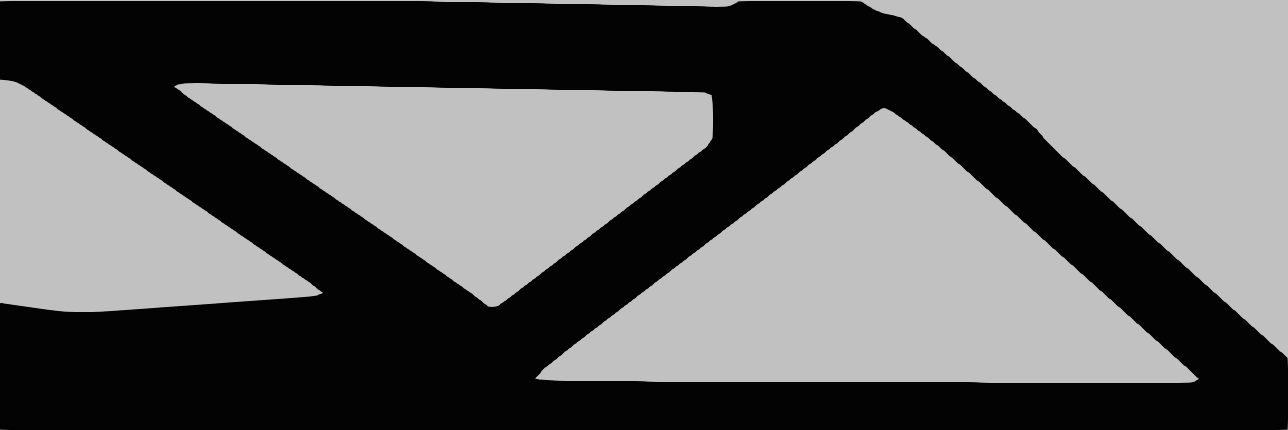}};
		\node[draw=none,white] at (-2, 1.05)   (5) {5};
        \node[draw=none,white] at (-3.5, -1)   (1) {1};
        \node[draw=none,white] at (-2.7, -0.8)   (7) {7};
        \node[draw=none,white] at (0.5, -1.2)   (2) {2};
        \node[draw=none,white] at (1.75, 0.75)   (6) {6};
        \node[draw=none,white] at (0.8, 0.75)   (10) {10};
        \node[draw=none,white] at (-2.75, 0.35)   (8) {8};
        \node[draw=none,white] at (-0.2, -0.25)   (9) {9};
        \node[draw=none,white] at (3.2, -0.5)   (11) {11};
		\end{tikzpicture}
		\caption{Adam design}\label{fig:prims_mean_adam}
	\end{subfigure}
	\caption{Computed Designs for the non-failing methods with $\lambda = 0$ for Example I. {The bar labels are also shown with respect to their corresponding positions in the final design.}}
	\label{fig:prims_mean_g}
\end{figure}

Figure~\ref{fig:prims_mean_graphs} shows the mean objective and constraint values for designs constructed by the non-failing methods when using $n=4$ realizations of uncertain inputs at each iteration to form the estimates of the objective and constraint functions, as well as their gradients. Figure \ref{fig:obj_exI} shows that among the SGD methods, Adam outperforms AdaGrad. However, GCMMA supplied with stochastic gradients turns out to be the best among all the methods compared here. The optimal $\ppm$ obtained are validated by evaluating the objective using 1000 random samples. For example, the square box in Figure \ref{fig:obj_exI} shows the objective for the optimized Adam design evaluated using 1000 random samples. Figure~\ref{fig:prims_mean_g} shows the final designs associated with the non-failing methods. 
%

\begin{figure}[htb!]
	\centering
	\begin{subfigure}[t]{\textwidth}
		\centering
		\includegraphics[scale =0.3]{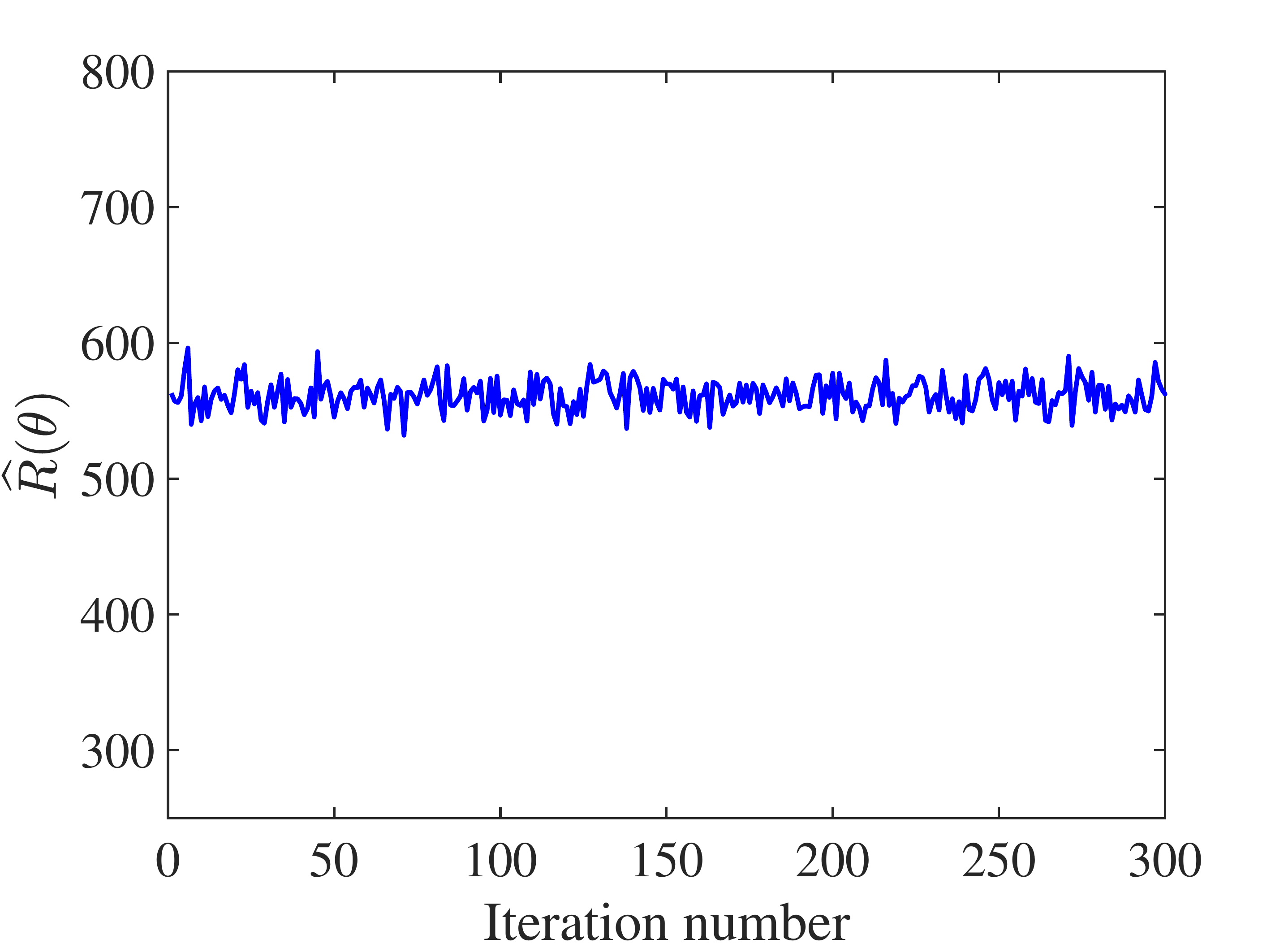}
		\caption{Adadelta (with $\zeta = 0.95$)}\label{fig:failing_a}
	\end{subfigure}\\
	\begin{subfigure}[t]{\textwidth}
		\centering
		\includegraphics[scale = 0.3]{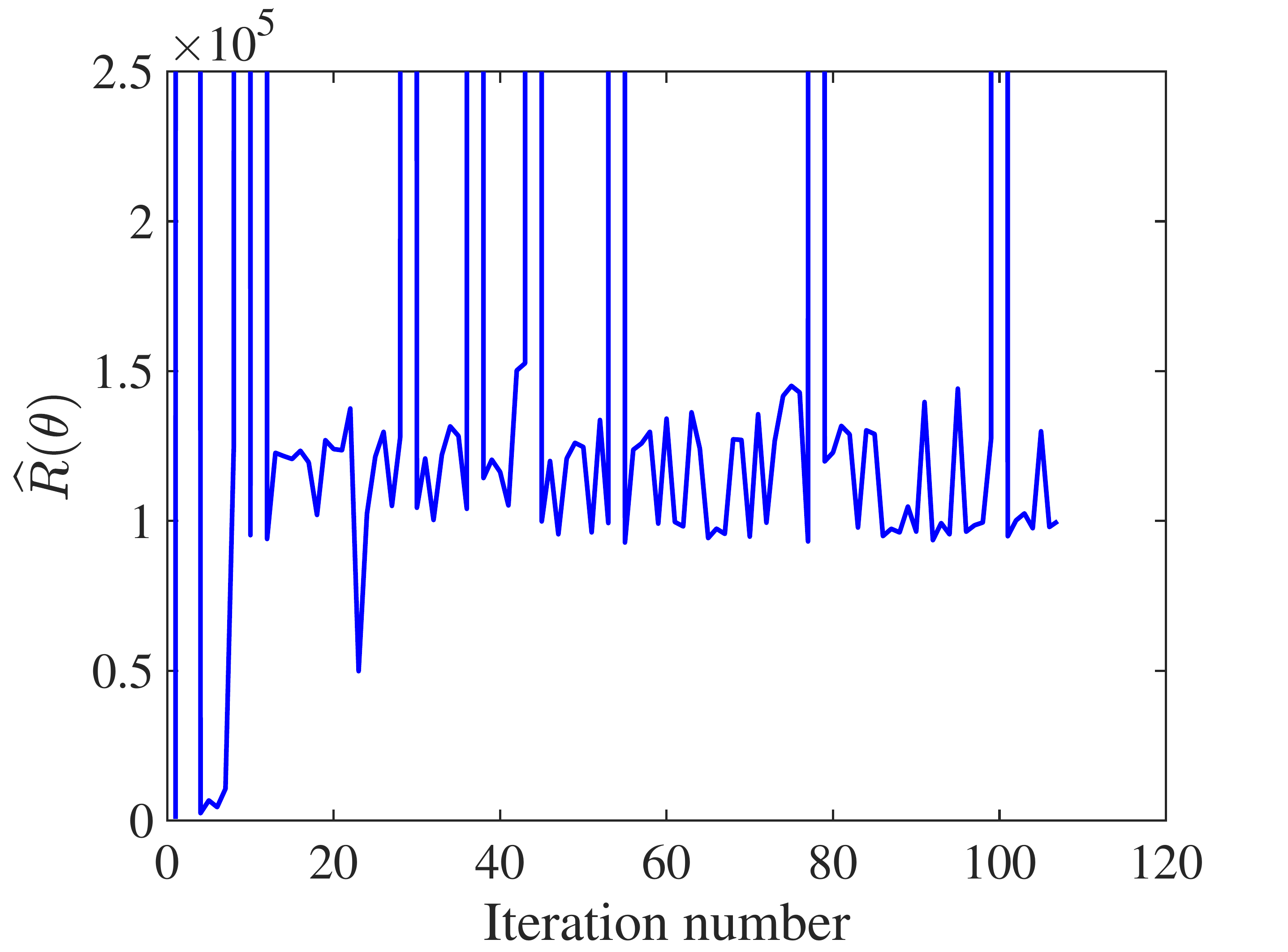}
		\caption{SAG (with $N=100$ and $n = 4$)}\label{fig:failing_b}
	\end{subfigure}
	\\
	\begin{subfigure}[t]{\textwidth}
		\centering
		\includegraphics[scale = 0.3]{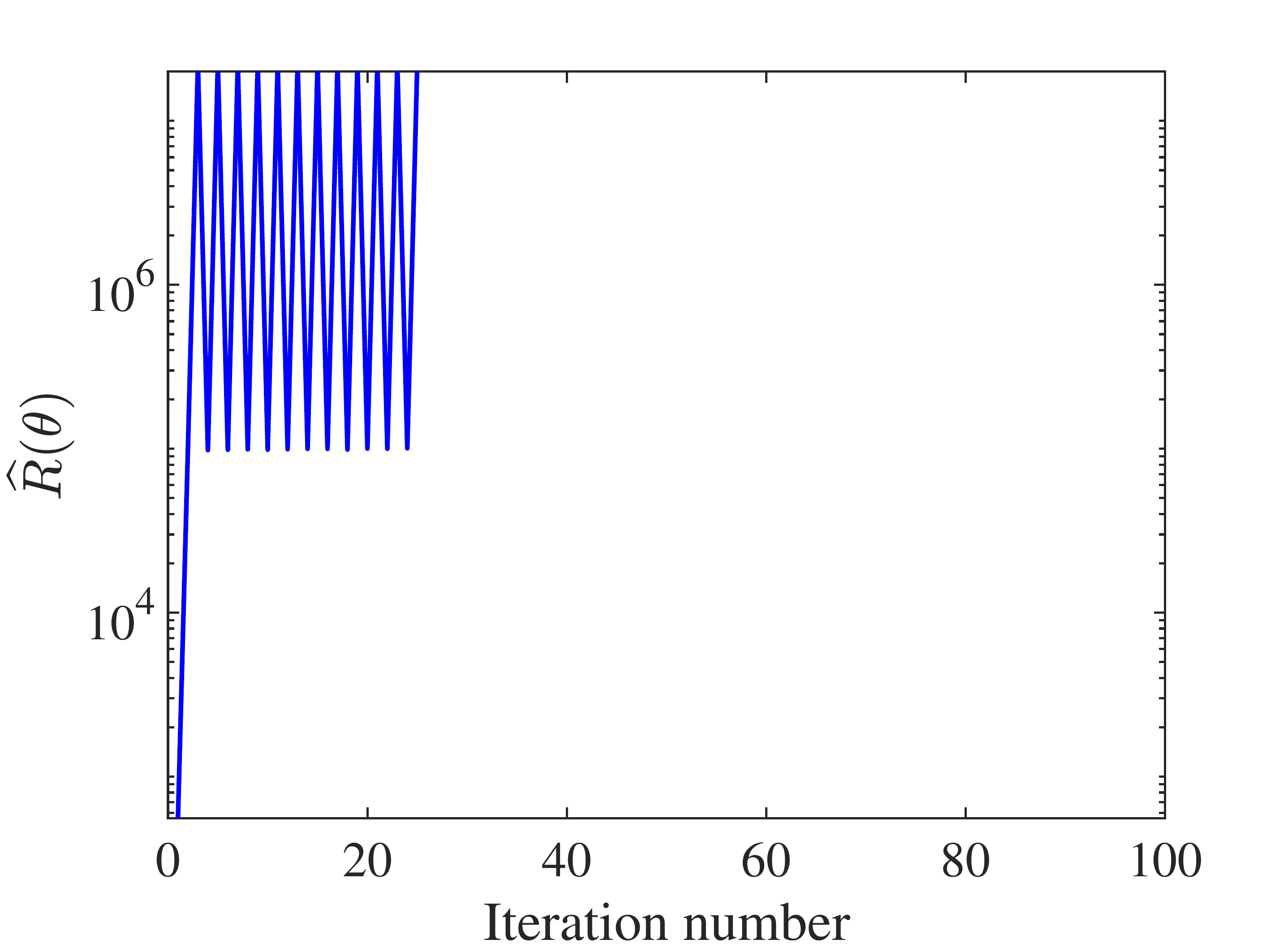}
		\caption{SVRG (with $N=20$, $m=4$, and $\eta=0.001$)}\label{fig:failing_c}
	\end{subfigure}
	\caption{Objectives for the failing methods with $\lambda = 0$ for Example I.}
	\label{fig:prims_failing}
\end{figure}

In this example, Adadelta with an exponential decay rate $\zeta = 0.95$ does not converge or converges very slowly as shown in Figure \ref{fig:failing_a}. We can explain this behavior by considering Algorithm \ref{alg:adadelta} in Section \ref{sec:methodology}, where we initialize $\am_\mathbf{h}$ and $\am_{\thetaa}$ to zero. This, coupled with $\rho$ close to 1.0 and small initial gradients, produce very small updates.
We can use smaller values for $\zeta$, which results in Adadelta relying mostly on gradient values from previous iterations.
However, significantly smaller $\zeta$ does not alleviate this issue entirely. 
%

The algorithms SAG and SVRG also fail to converge. The objective history of SAG and SVRG are shown in Figures \ref{fig:failing_b} and \ref{fig:failing_c}, respectively. 
%
%
Here, SAG is implemented with an ensemble size of $\Ns=100$ and $n=4$ gradients are updated at every iteration. In TO, the design may go through large changes during the initial iterations. Hence, the failure of the algorithm is due to the significant use of gradient information from past iterations. For example, in SAG, only a handful of the gradients are updated at every iteration and most of the gradients from previous iterations are used. Updating many gradients (\textit{i.e.}, large $n$) at every iteration will produce converged results when using SAG  as it will essentially become a full-batch gradient descent. However, this increases the overall computational cost of the optimization considerably compared to Adam and AdaGrad. We, therefore, do not pursue SAG any further in this work.

Similarly, in SVRG, gradients and parameter estimates from past iterations are maintained (see Algorithm \ref{alg:svrg}). For example, in Figure \ref{fig:failing_c}, SVRG with a small learning rate $\eta= 0.005$ and a smaller number of inner iterations $m=4$, \textit{i.e.}, $\hb_{\mathrm{best}}$ can be from the last four iterations, shows the algorithm starts diverging when there is a significant change in the design (near iteration 75 in this case). For larger $\eta$ and larger $m$, the algorithm diverges earlier. A smaller number of inner iterations, $m$, along with a much smaller learning rate, $\eta$, can be used, but this will increase the overall computational cost as $\hb_{\mathrm{best}}$ will be evaluated more frequently as well as it will slow down the convergence. Hence, significant non-convexity of the objective function and heavy use of previous gradients cause SAG and SVRG algorithms to fail in this example.

\begin{figure}[htb!]
	\centering
	\includegraphics[scale=0.35]{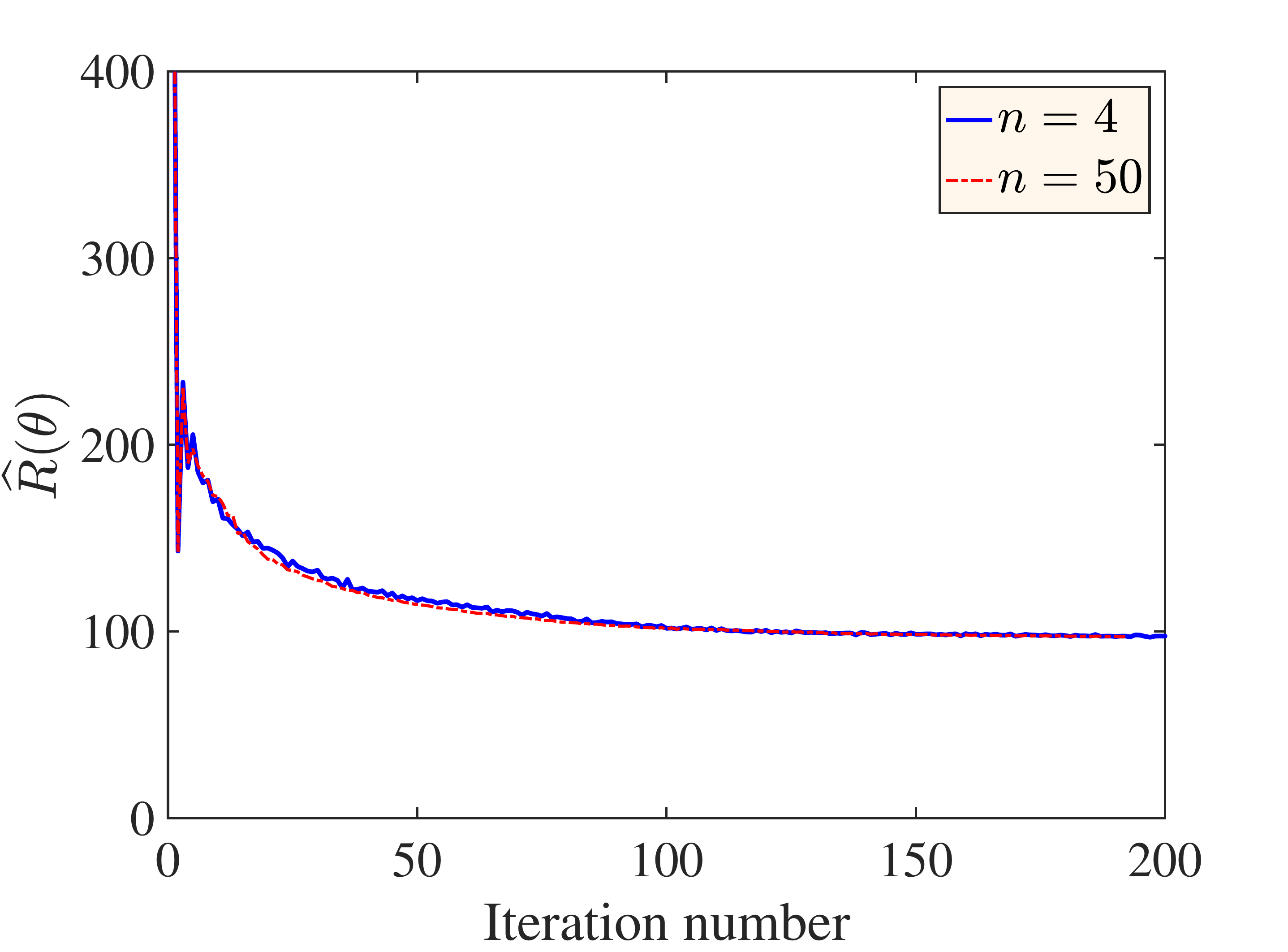}
	\caption{Plot of objectives using AdaGrad with $n=4$ and $n=50$ samples per iteration with $\lambda = 0$ for Example I.}
	\label{fig:prims_mean_by_sample}
\end{figure}

%
Next, we use $n=50$ random samples per iteration and compare the objective for AdaGrad with our previous results obtained using $n=4$ random samples. Figure~\ref{fig:prims_mean_by_sample} shows that there is no noticeable improvement in the objective even though performing $50$ forward solves per iteration is computationally expensive. 

%

\begin{figure}[htb!]
	\centering
	\begin{subfigure}[t]{\textwidth}
		\centering
		\includegraphics[scale =0.2625]{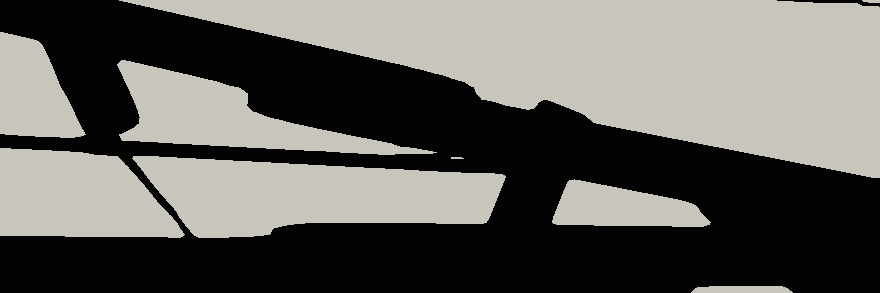}
		\caption{Untuned design}
	\end{subfigure}%
	\\
	\begin{subfigure}[t]{\textwidth}
		\centering
		\includegraphics[scale = 0.18]{figs/Adagrad_lambda0_iter1000}
		\caption{ Tuned design}
	\end{subfigure}
	\caption{Tuning of learning rate $\eta$ plays an important role in { the success of the SGD} methods. This figure shows final designs using AdaGrad { without and with $\eta$ tuned}.}
	\label{fig:prims_tuned_untuned}
\end{figure}

\begin{figure}[htb!]
    \begin{tikzpicture}
    	\centering
	\node (obj) at (0,0) {\includegraphics[scale=0.35]{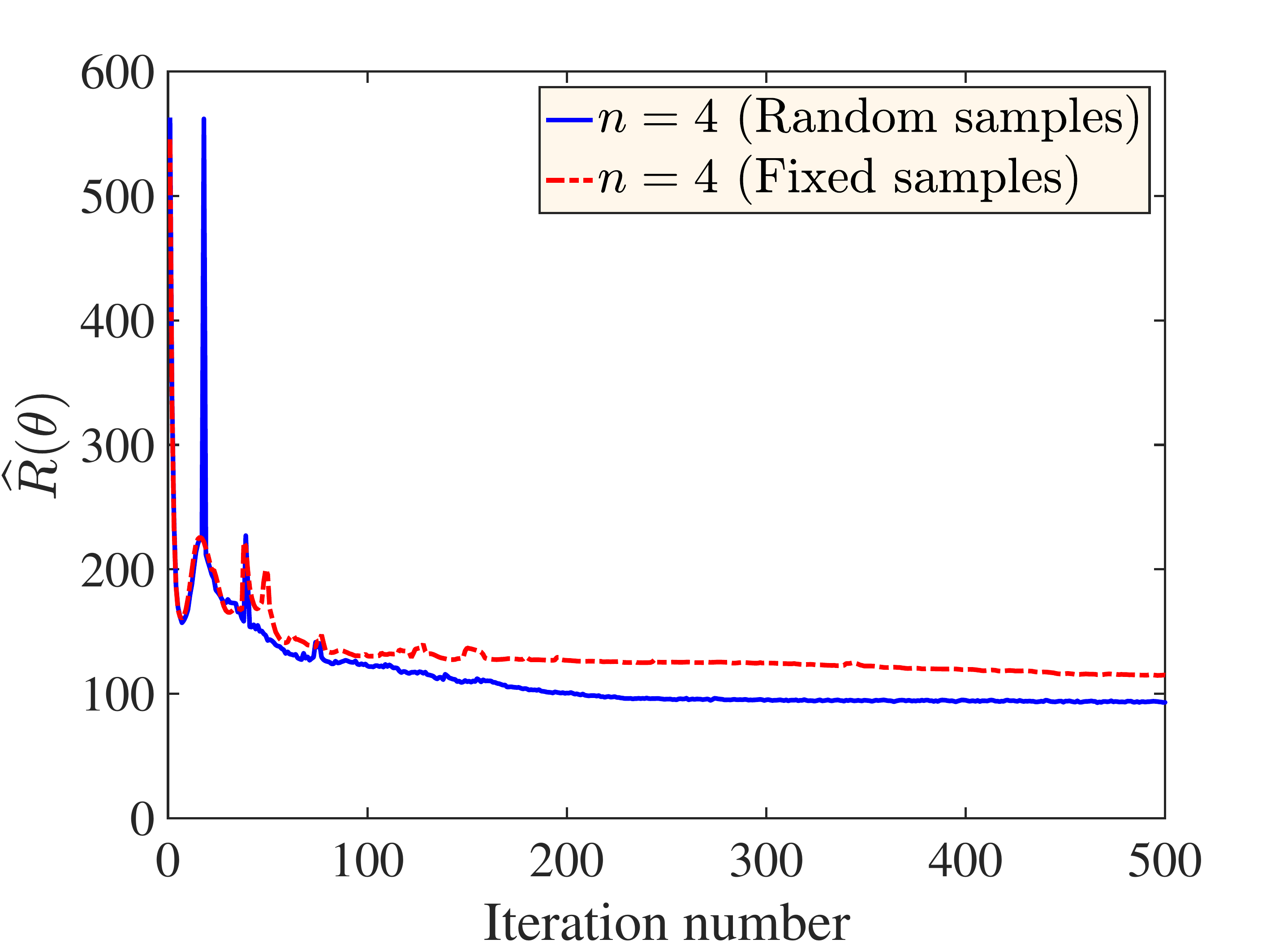}};
    \node (fixed_design) at (6.75,0) {\includegraphics[scale=0.08]{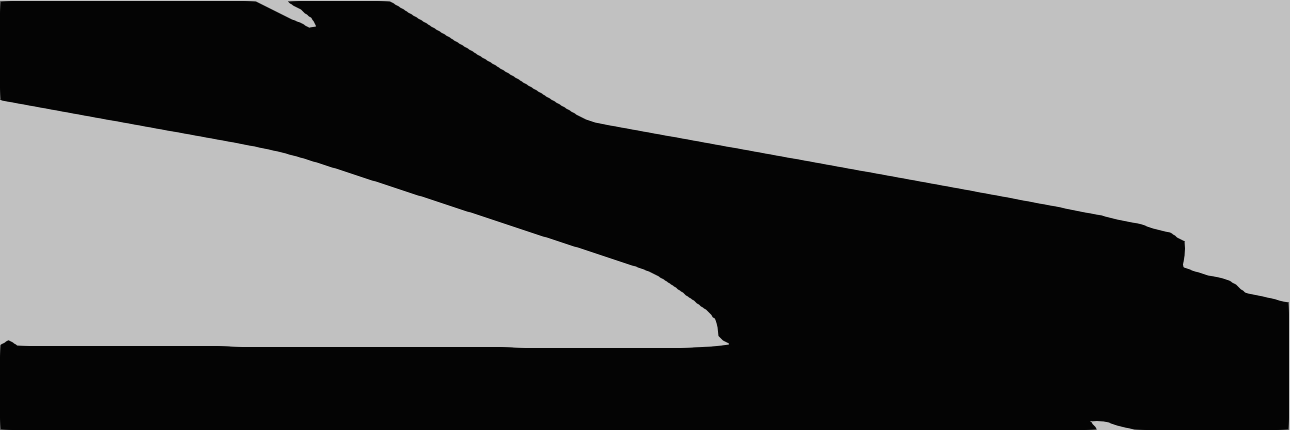}};
    \node (adam_design) at (6.75,-2) {\includegraphics[scale=0.08]{figs/Adam_lambda0_iter1000}};
    \draw[-latex,thick] (4.75,-2) -- (4.2,-1.8);
    \draw[-latex,thick] (4.75,-0.5) -- (4.2,-1.5);
    \end{tikzpicture}
	\caption{Plot of objectives using AdaGrad with fixed and random samples per iteration with $\lambda = 0$ for Example I.}
	\label{fig:prims_fixedsample}
\end{figure}

The success of the SGD methods (except Adadelta) depends on the choice of the learning rate $\eta$. Figure \ref{fig:prims_tuned_untuned} shows results from AdaGrad with untuned and tuned $\eta$. Tuning $\eta$ can be performed in a computationally inexpensive manner by performing only a few design optimization iterations with different $\eta$ and selecting the most promising $\eta$. Furthermore, it is important to note that the success of these SGD methods for TO strongly depends on the use of an independent set of random samples $\{\Ym_i\}_{i=1}^n$ at every iteration.  Figure \ref{fig:prims_fixedsample} shows that if the same sample set $\{\Ym_i\}_{i=1}^4$ is used at every iteration the optimization method AdaGrad may get trapped in a local minimum similar to the untuned learning rate setting.

\subsubsection{Robust Design}
\label{subsec:prims_var_opt}

Next, we solve the optimization problem defined in \eqref{eq:opt_def2}, for $\lambda = 0.1$. Here, we optimize a combination of mean and the variance of the objective subject to a similar combination of mean and variance of constraint values. Note that in this case the failures of some of the SGD methods are even more pronounced since the objective scales with the square of the compliance now. Again, we utilize a learning rate $\eta = 0.05$ for the algorithms in Section~\ref{sec:methodology} and a sample size of $n=4$. Again, we confine the presentation of the results to algorithms that converge to more or less meaningful designs. For example, Adadelta, SVRG, and SAG fail to produce a meaningful design due to the reasons discussed before.
%

\begin{figure}[htb!]
	\centering 
	\begin{subfigure}[t]{0.5\textwidth}
		\centering
		\includegraphics[scale = 0.18]{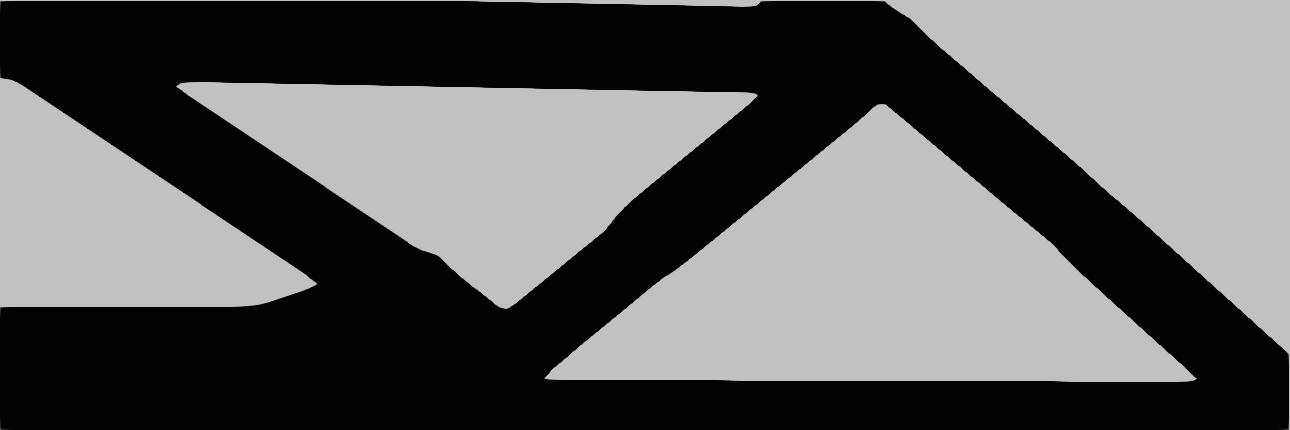}
		\caption{GCMMA design}
	\end{subfigure}
	\\
	\begin{subfigure}[t]{0.5\textwidth}
		\centering
		\includegraphics[scale = 0.18]{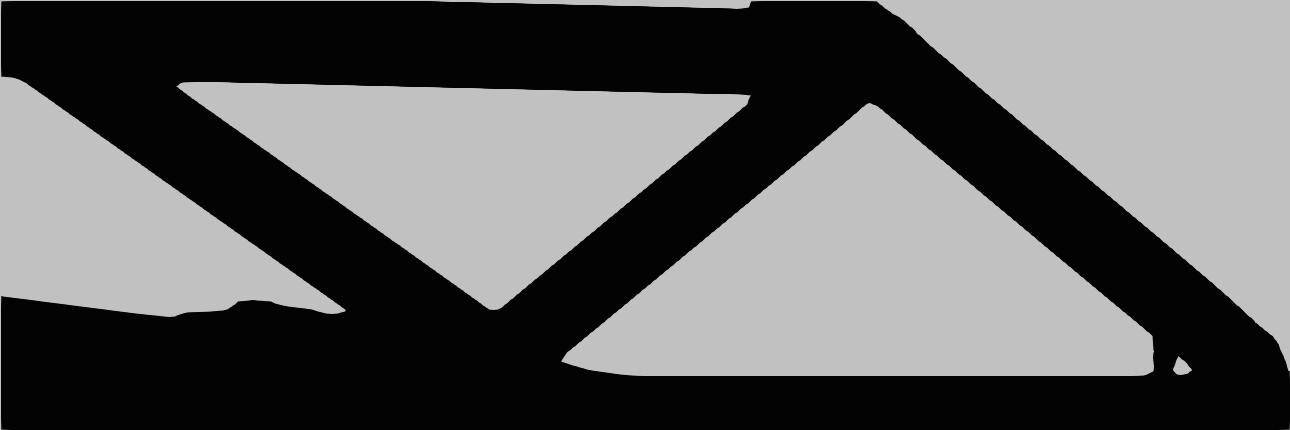}
		\caption{Adam design}
	\end{subfigure}
	\caption{Optimal robust designs for the non-failing methods with $\lambda = 0.1$ for Example I.}
	\label{fig:prims_var_lambda0_1}
\end{figure}
\FloatBarrier

\begin{figure}[htb!]
	\centering
	\begin{subfigure}[t]{\textwidth}
		\centering
		\includegraphics[scale =0.35]{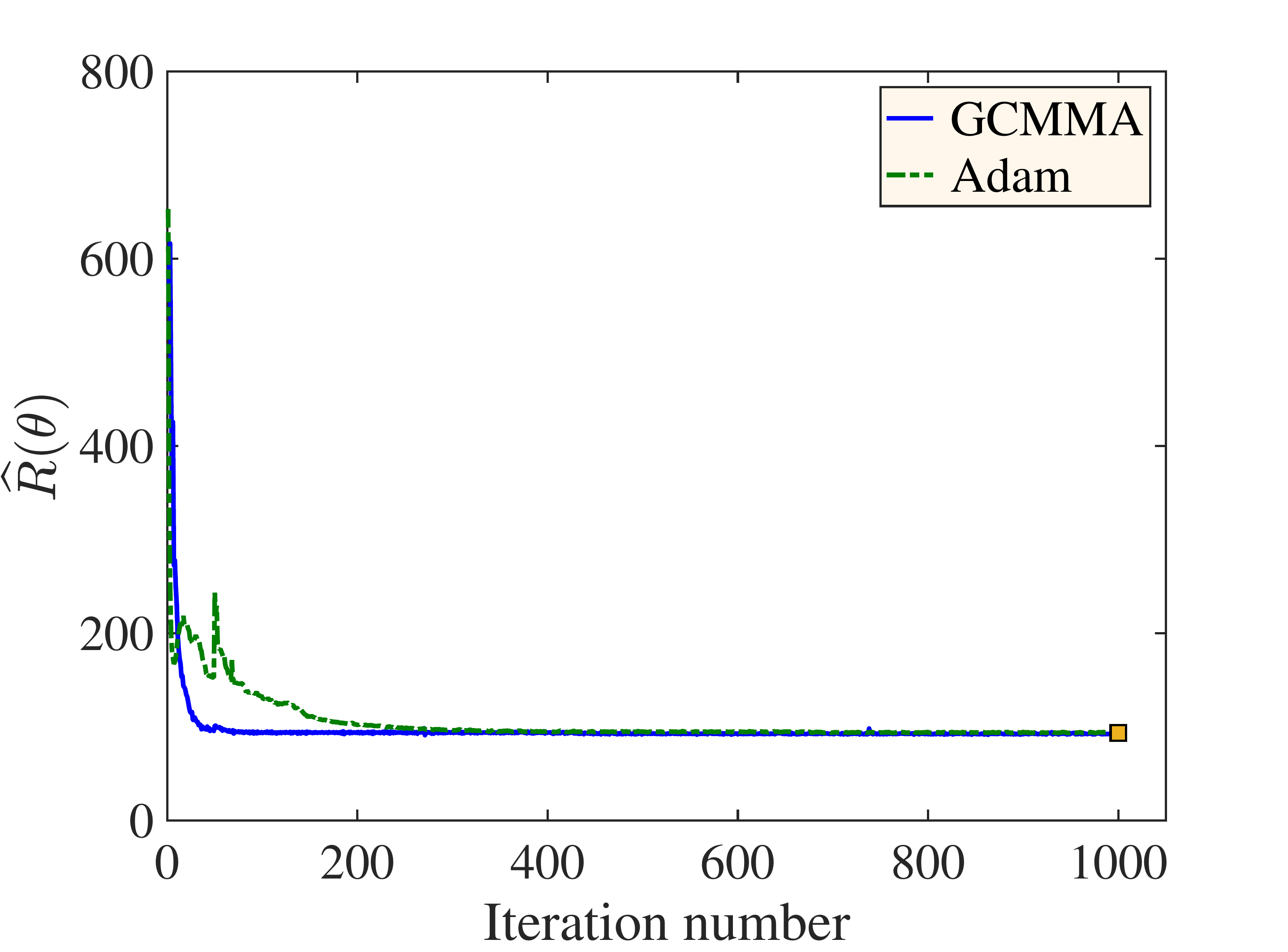}
		\caption{Objective estimates (the square at the right end shows $\widehat{R}(\ppm)$ for the optimal design from Adam evaluated using 1000 samples)}\label{fig:ExI_ObjVar}
	\end{subfigure}%
	\\
	\begin{subfigure}[t]{\textwidth}
		\centering
		\includegraphics[scale = 0.35]{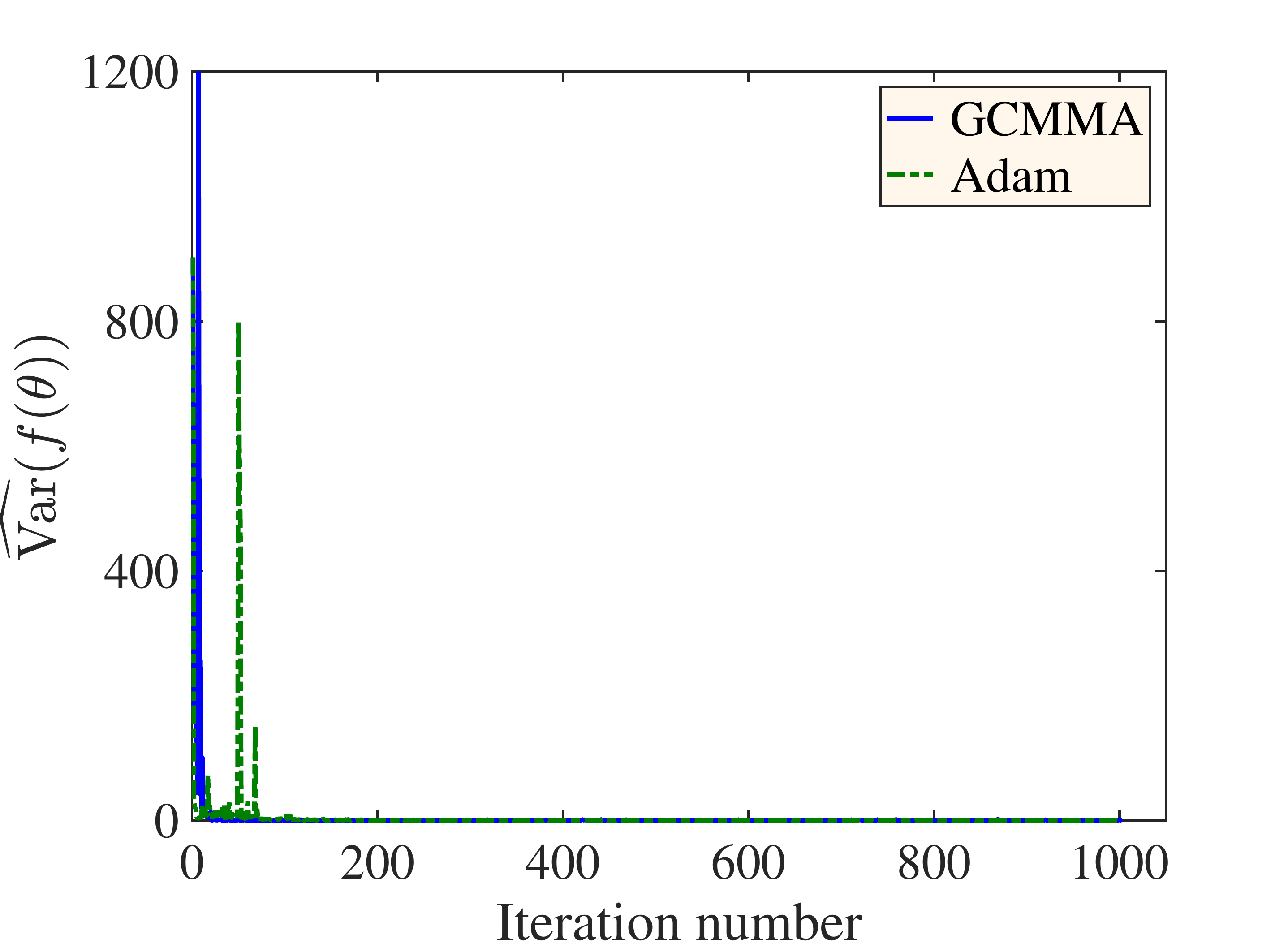}
		\caption{Variance estimates}
	\end{subfigure}
	\caption{Objectives and variances for the non-failing methods with $\lambda = 0.1$ for Example I.}\label{fig:ExI_VarVar}
	\label{fig:prims_var_graphs}
\end{figure}

Among the SGD methods, only Adam maintains its stability while all other methods fail to converge to meaningful designs. Figure~\ref{fig:prims_var_lambda0_1} shows the design produced by GCMMA and Adam. Figure~\ref{fig:prims_var_graphs} shows the objective function and the variance part of it for the designs constructed by GCMMA and Adam when using $n=4$ realizations of the uncertain inputs.  

{We observe from Figure~\ref{fig:prims_var_graphs} that the objective exhibits a low variance through the later optimization iterations. This, together with the use of independent samples of the uncertain input for each new design, implies that indeed the variance of the objective, estimated otherwise by a large number of samples, is small. We also confirm the small variability of the objective by evaluating the optimal solution obtained from Adam using 1000 random samples (the square at the right end in Figure \ref{fig:ExI_ObjVar}), which is very similar to the one obtained using only 4 random samples.}

\begin{figure}[htb!]
	\centering
	\includegraphics[scale = 0.35]{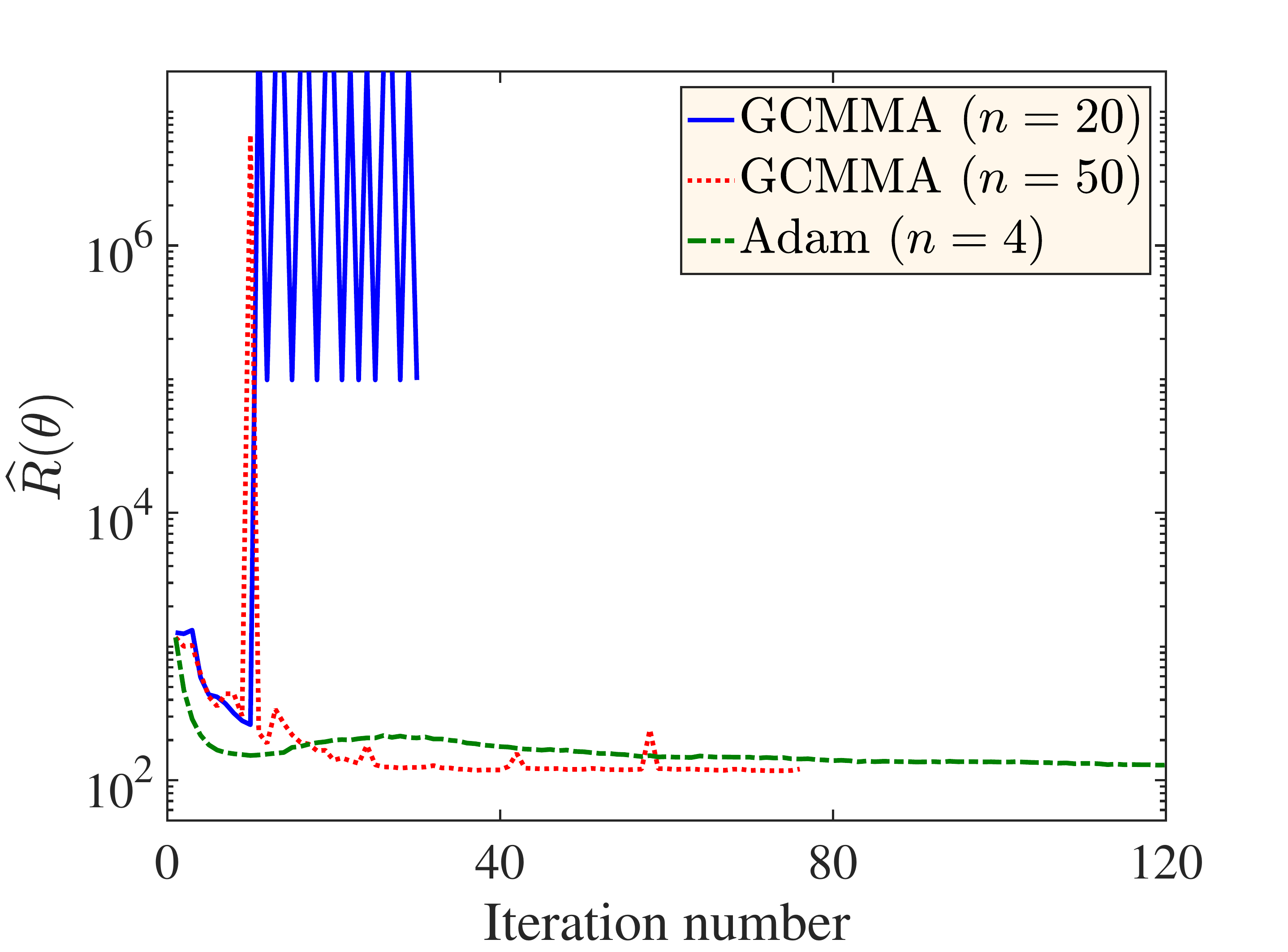}
	\caption{Objective estimates with $\lambda =1$ for Example I.}
	\label{fig:gcmma_fails}
\end{figure}

To observe the effect of $\lambda$ in the objective, next, we increase {$\lambda$ to 1}. In Figure \ref{fig:gcmma_fails}, we compare Adam and GCMMA for various mini-batch sizes. We run Adam with $n=4$, which shows the convergence of the objective. However, GCMMA with $n=4$ and $20$ does not converge. In Figure \ref{fig:gcmma_fails}, we show the objective for GCMMA using $n=20$ with a solid blue line. We then increase $n=50$ for GCMMA and note that GCMMA converges (shown using the red dotted line in Figure~\ref{fig:gcmma_fails}). 
{Note that, during the early optimization iterations, the variance of the objective remains large and with larger $\lambda$ the variance of the stochastic gradient estimates will be consequently larger. Therefore, to achieve a similar level of variance in the gradient estimates, more random samples are needed in each iteration.
The results in Figure~\ref{fig:gcmma_fails} indicate the sensitivity of GCCMA's performance with respect to the variance of the stochastic gradients. A theoretical study must be performed on the convergence of GCMMA when supplied with stochastic gradients to quantify this effect of $\lambda$. However, this is beyond the scope of the current paper. With $n=50$ random samples per iteration, GCMMA performs similarly to Adam using $n=4$ samples.}

%
%
This example uses a beam design problem under uncertainty to show the use of stochastic gradients for TO. We employ different stochastic gradient descent algorithms, namely, AdaGrad, Adadelta, Adam, SVRG, and SAG using only a handful of random samples per iteration. Among these algorithms, AdaGrad and Adam perform well after some tuning of their learning parameter and provides meaningful designs for $\lambda=0$ in the objective (see \eqref{eq:exp_risk}). However, AdaGrad fails to provide a meaningful design when we increase $\lambda$ to 0.1. GCMMA, which is typically not used for stochastic design problems, however, performs well and produces meaningful designs for $\lambda=0$ and $\lambda= 0.1$. Hence, among all the algorithms studied for this example, Adam and GCMMA turn out to be the most useful to optimize the objective defined in \eqref{eq:exp_risk} and \eqref{eq:glam_def}.
Note that, only first order gradients are used in all of these methods and hence their costs are comparable. 


\subsection{Example II: Load Support over Elastic Bedding}
\label{sec:elastic_bedding}

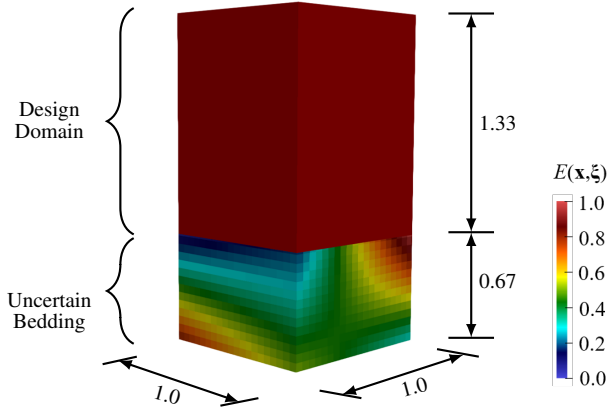
\begin{figure}[htb!]
	\centering
\begin{tikzpicture}
\node[inner sep=0pt] (block) at (0,0)	{\includegraphics[scale = 0.35]{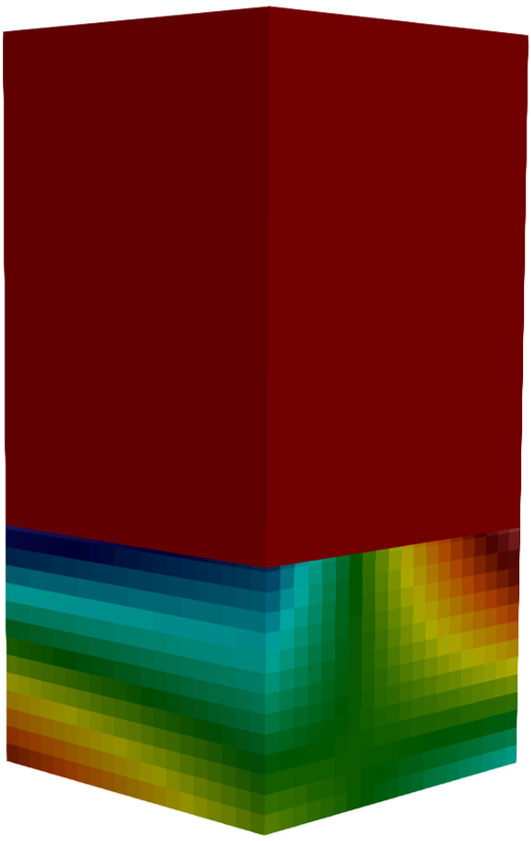}};
\draw [thick,decorate,decoration={brace,amplitude=10pt},xshift=-4pt,yshift=-25pt] (-2,0.25) -- (-2,3.25) node [black,midway,xshift=-0.6cm] {};
\draw [thick,decorate,decoration={brace,amplitude=10pt},xshift=-4pt,yshift=-25pt] (-2,-1.2) -- (-2,0.2) node [black,midway,xshift=-0.6cm] {};
\node[draw=none] at (-3.25, 1)     (a2)    {Design};
\node[draw=none] at (-3.25, 0.75)  (b2)    {Domain};
\node[draw=none] at (-3.25, -1.5)  (a2)    {Uncertain};
\node[draw=none] at (-3.25, -1.8)  (b2)    {Bedding};

\draw[thick, latex-latex] (2.3,-0.6) -- (2.3,2.3); 
\draw[thick] (2,-0.6) -- (2.6,-0.6); 
\draw[thick] (2,2.3) -- (2.6,2.3);
\node[draw=none] at (2.65,0.85)    {1.33};

\draw[thick, latex-latex] (2.3,-0.6) -- (2.3,-2);
\draw[thick] (2,-2) -- (2.6,-2);
\node[draw=none] at (2.65,-1.25)    {0.67};

\draw[thick, latex-latex] (2.2,-2.225) -- (0.6,-2.75);
\draw[thick] (1.85,-2.15) -- (2.4,-2.25);
\draw[thick] (0.305,-2.7) -- (0.855,-2.8);
\node[draw=none,rotate = 15] at (1.55,-2.65)    {1.0};

\draw[thick, latex-latex] (-2.35,-2.3) -- (-0.75,-2.85);
\draw[thick] (-0.5,-2.8) -- (-1,-2.9);
\draw[thick] (-2.1,-2.25) -- (-2.55,-2.35);
\node[draw=none,rotate = -15] at (-1.7,-2.75)    {1.0};
\node[draw=none] at (4,-1.2) {\includegraphics[scale = 0.5]{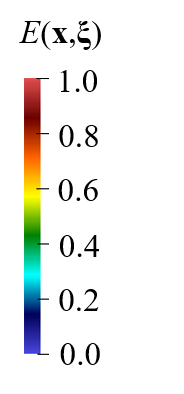}};
    \end{tikzpicture} 
	\caption{Design area (top) and bedding (bottom) for Example II; the color contours show the Young's modulus distribution for specific realization of the uncertain bedding.}
    \label{fig:block_schematic}
\end{figure}

In the second example, we consider a 3D design problem with a structure being subject to uncertain loading and resting on an uncertain bedding.
A design domain of size $1.0 \times 1.0 \times 1.33$ rest on top of a non-design domain of size $1.0 \times 1.0 \times 0.67$; see Figure~\ref{fig:block_schematic}. 
%
%
The non-design domain is occupied by a material with random stiffness and represents an uncertain bedding. The non-design domain is clamped at the bottom face. At the center of the top face of the design domain, a point load with random  direction is applied. The performance measure of the objective function is the strain energy, and the only constraint is to ensure that the mass-ratio of the structure is no more than 15\% of the maximum design mass, \ie~the entire design-domain filled with bulk material. In this example, to enforce the constraint $\kappa = 1000$ is used.

%
{In this example, the structural geometry is described by the SIMP method with $\beta_\mathrm{SIMP}=3$ described in Section \ref{sec:des_model}. For projection in \eqref{eq:projection}, we use $\nu_\mathrm{pr}=0.0001$ and set $\beta_\mathrm{pr}$ to 5 for first 400 iterations with GCMMA and AdaGrad and then it is increased to 20. With Adam, we use $\beta_\mathrm{pr}=2$ for first 400 iterations and then increase it to 20.}
The random stiffness distribution in the non-design domain, \ie~the Young's modulus $E(\xm,\Ym)$, is defined by the following tri-linear interpolation
\begin{equation}
    E(\xm,{\Ym}) = E_0 \sum_{i=1}^8 N_i(\xm) \ \xi_i,
\end{equation}
where $E_0$ is a constant and assumed as unity, $N_i$ are the tri-linear shape functions, and $\xi_i$ the random variables associated with the Young's moduli at the corner points of the non-design domain. The random variables $\xi_i$ are sampled from uniform distribution in $(0,1)$.
%
Uncertainty in the direction of the point load is described by two additional uniform random variables that define the direction cosines of the load vector as follows
\begin{equation}
\begin{split}
        \cos \gamma_x &= \sin(\pi\xi_9)\sin(2\pi\xi_{10});\\
        \cos \gamma_y &= \sin(\pi\xi_9)\cos(2\pi\xi_{10});\\
        \cos \gamma_z &= \cos(\pi\xi_{9}),\\
\end{split}
\end{equation}
where $\gamma_x$, $\gamma_y$, and $\gamma_z$ are the angles between the load vector and the $x$, $y$, and $z$ axes, respectively. In this example, we assume the Poisson's ratio as 0.3, the magnitude of the point load as unity. 
The finite element model consists of tri-linear hexahedron elements with edge length as 0.05. This results in a mesh with 16,000 design parameters.


\begin{figure}[htb!]
	\centering
	\includegraphics[scale=0.35]{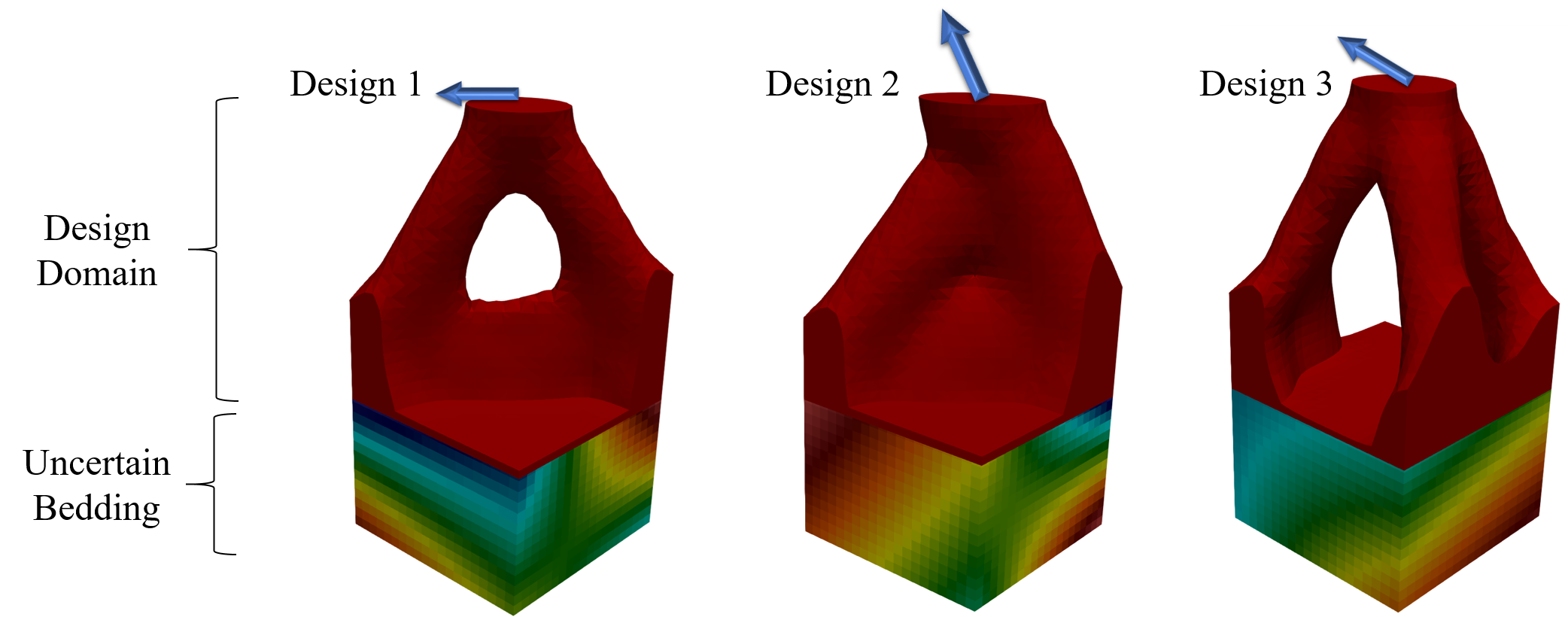}~~~~\includegraphics[scale = 0.5]{./figs/colorbar_ExII}
	\caption{Deterministic optimization results for three realizations of the random variables, {namely, tri-linearly interpolated Young's modulus and direction of the load}. }
	\label{fig:detdes}
\end{figure}

As the following numerical studies show, this optimization problem demonstrate key design differences that can arise from different stochastic gradient methods and GCMMA. We only use AdaGrad and Adam among the SGD methods in this example as the results for the previous example suggest that these algorithms are more useful to TO problems under uncertainty. For illustration purposes, we present first results for the deterministic formulation of the optimization problem. Figure \ref{fig:detdes} shows optimized geometries with a density threshold of 0.5 for certain realizations of the random variables, highlighting the strong influence of the uncertain parameters on the design solutions.

\subsubsection{Average Design}
\label{subsec:block_mean_opt}

\begin{figure}[htb!]
	\centering
	\begin{subfigure}[t]{\textwidth}
		\centering
		\includegraphics[scale = 0.35]{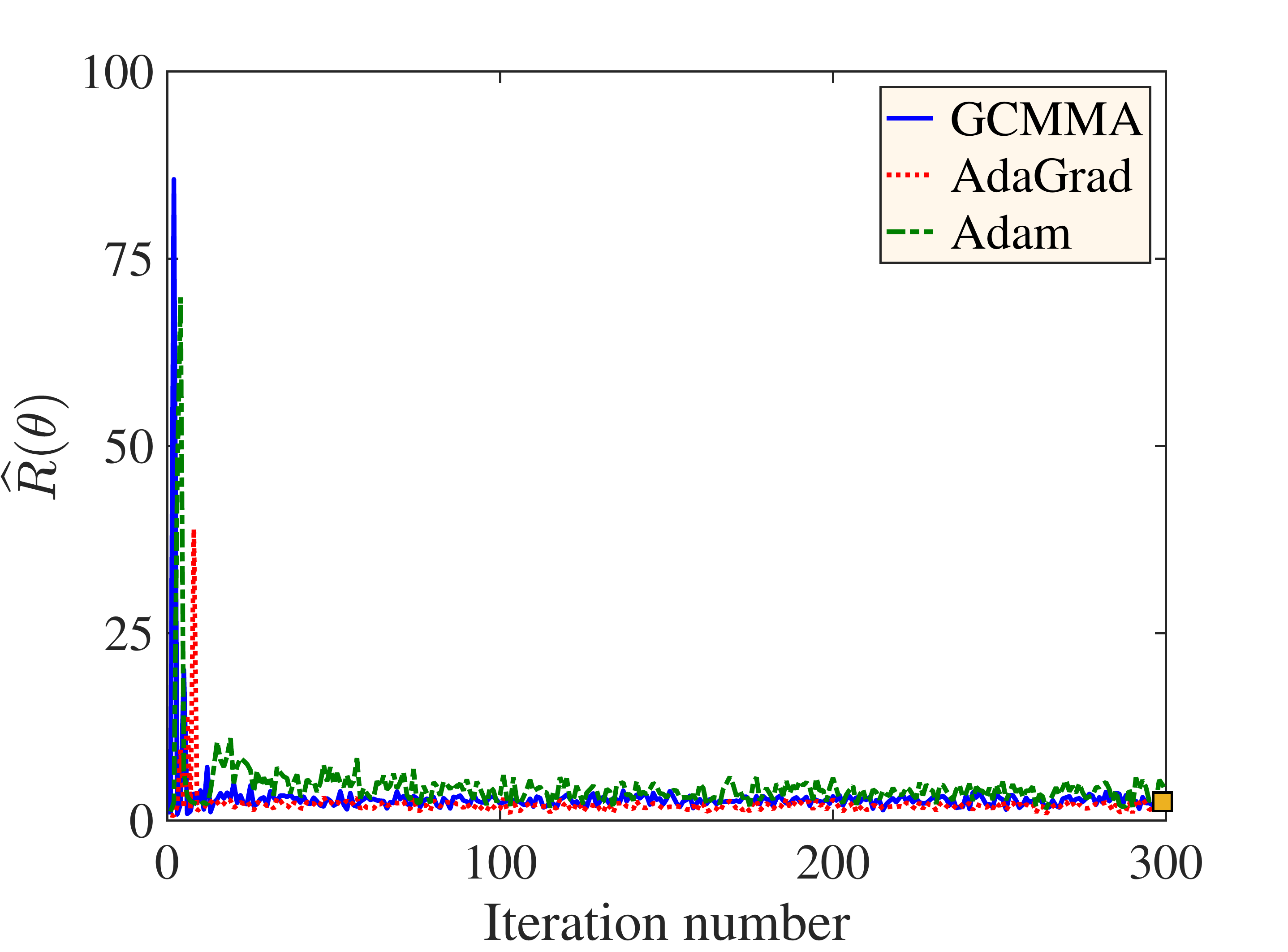}
		\caption{Objective estimates (the square at the right end shows $\widehat{R}(\ppm)$ for the optimal design from Adam evaluated using 1000 random samples)}\label{fig:obj_exII}
	\end{subfigure}
	\\
	\begin{subfigure}[t]{\textwidth}
		\centering
		\includegraphics[scale = 0.35]{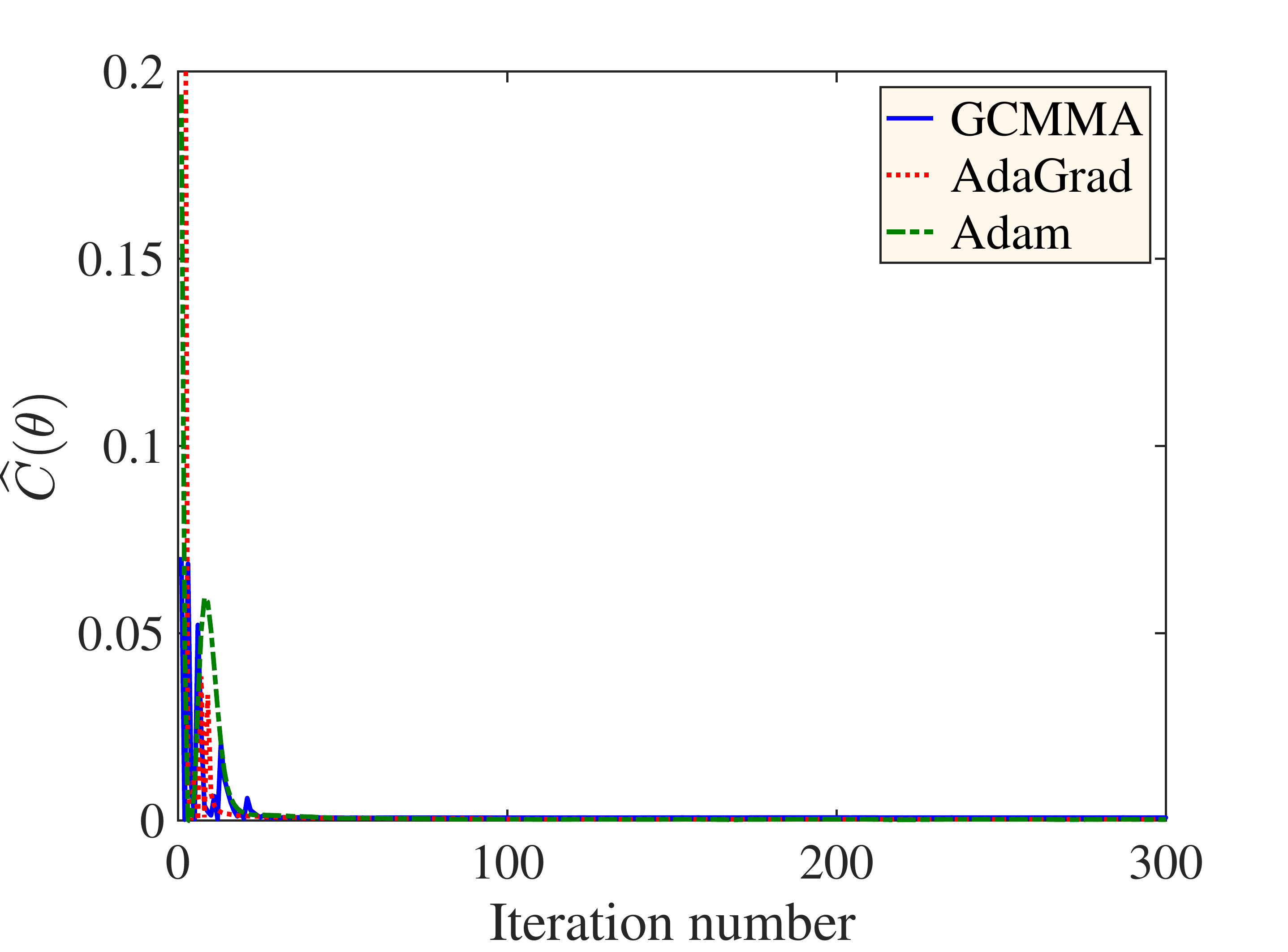}
		\caption{Constraint estimates}
	\end{subfigure}
	\caption{Objectives and constraints for three methods, {namely, GCMMA, AdaGrad, and Adam} with $\lambda = 0$ for Example II.}
	\label{fig:ExII_obj_con}
\end{figure}

First, we consider solving the optimization problem (\ref{eq:exp_risk}) and (\ref{eq:glam_def}) using the formulation in \eqref{eq:opt_def2}, for $\lambda = 0$, \textit{i.e.}, we do not account for the variance terms in the objective and constraint. We utilize a learning rate $\eta= 0.25$ for the algorithms in Section~\ref{sec:methodology} and a sample size of $n=4$. The value for the learning rate was determined by a few preliminary runs.  Figure~\ref{fig:ExII_obj_con} shows the evolution of the mean objective and constraint values. {Further, the optimal $\ppm$ obtained are validated by evaluating the objective using 1000 random samples. For example, the square box in Figure \ref{fig:obj_exII} shows the objective for the optimized Adam design evaluated using 1000 random samples.}

\begin{figure}[htb!]
	\centering
	\begin{subfigure}[b]{\textwidth}
		\centering
		\includegraphics[scale =0.33]{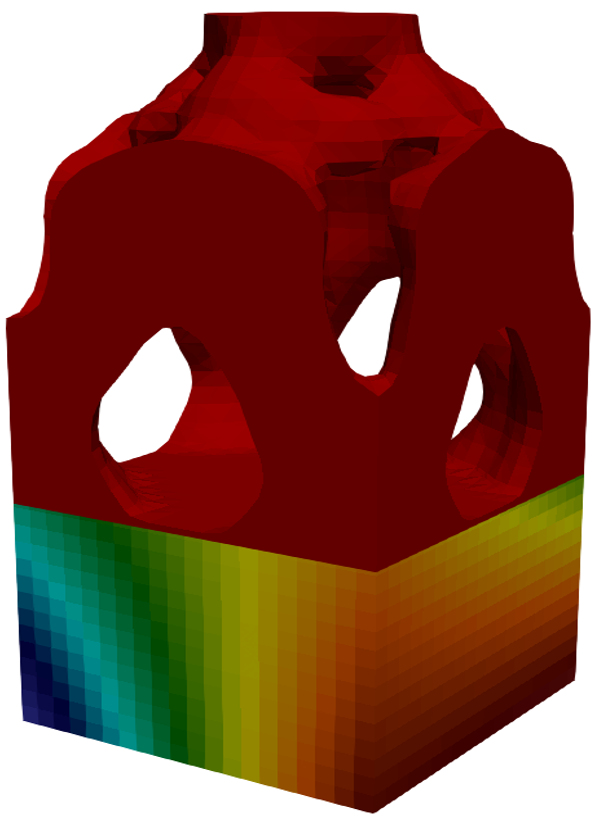}~~~~\includegraphics[scale =0.34]{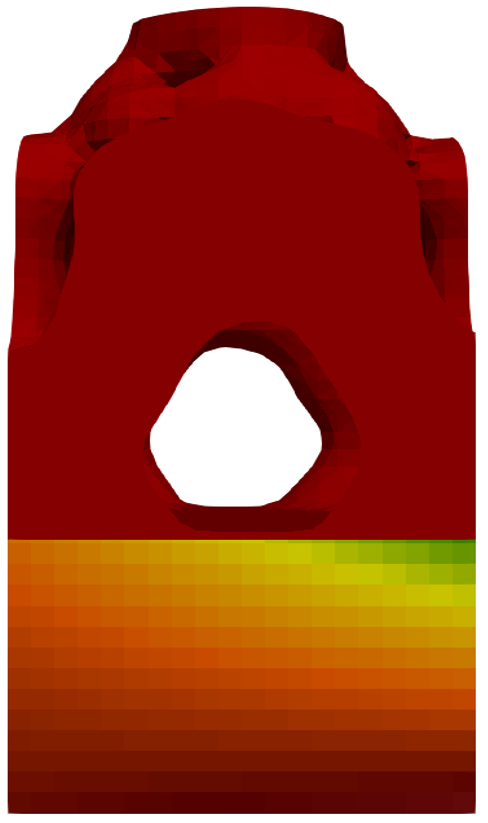}~~~~
		\includegraphics[scale =0.8]{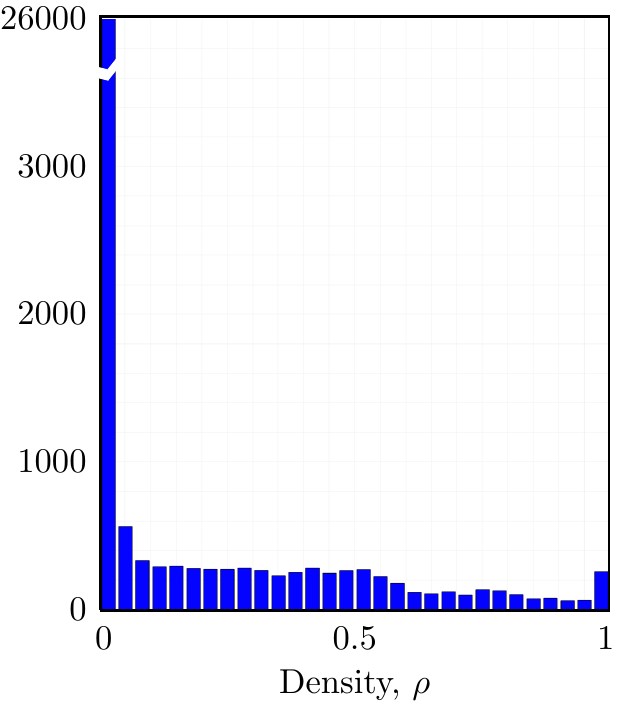}
		\\\includegraphics[scale = 0.5]{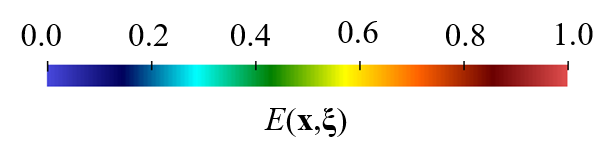}
		\caption{GCMMA design}
		\label{fig:exII_gcmma1}
	\end{subfigure}%
	\\
	\begin{subfigure}[b]{\textwidth}
		\centering
		\includegraphics[scale =0.32]{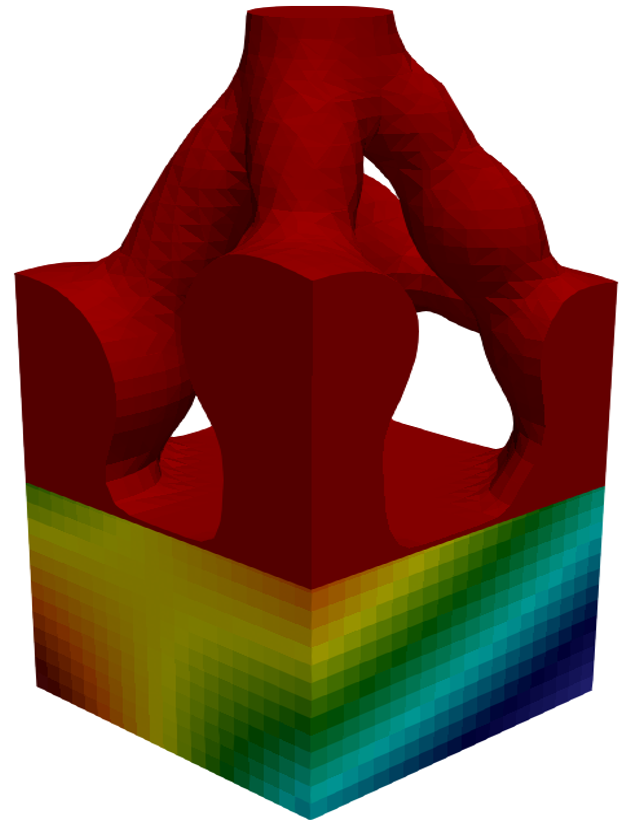}~~~~\includegraphics[scale =0.34]{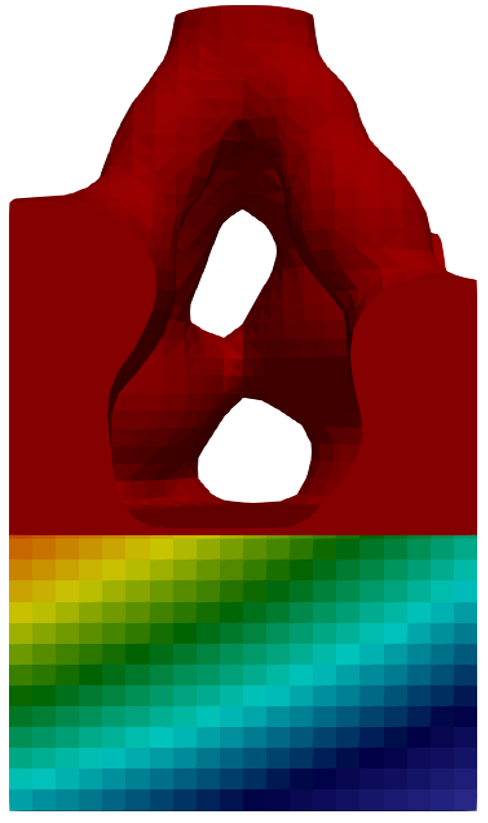}
		~~~~\includegraphics[scale =0.8]{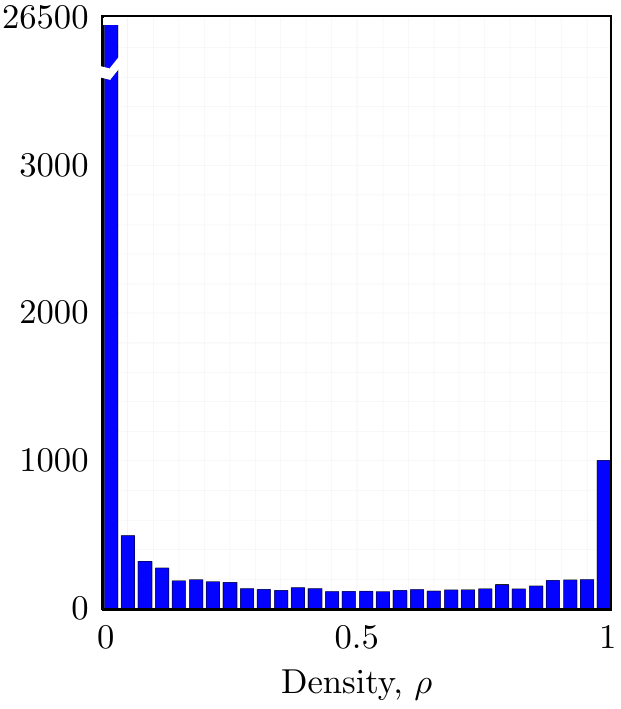}
		\\\includegraphics[scale = 0.5]{./figs/colorbar_ExII_horz}
		\caption{AdaGrad design}
		\label{fig:exII_adagrad1}
	\end{subfigure}%
	\\
	\begin{subfigure}[b]{\textwidth}
		\centering
		\includegraphics[scale =0.34]{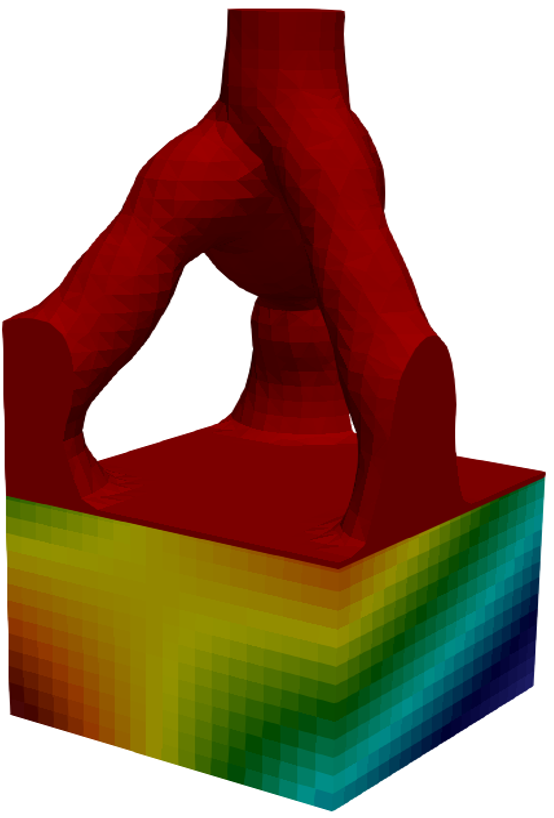}~~~~\includegraphics[scale =0.34]{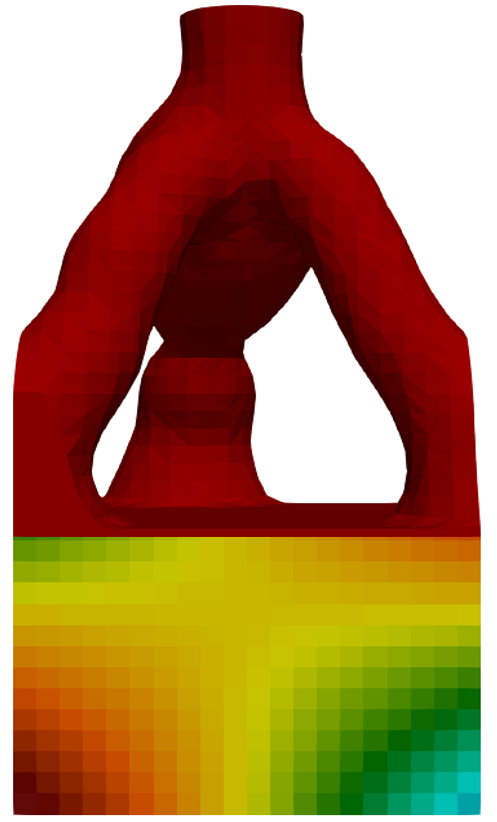}~~~~
		\includegraphics[scale =0.8]{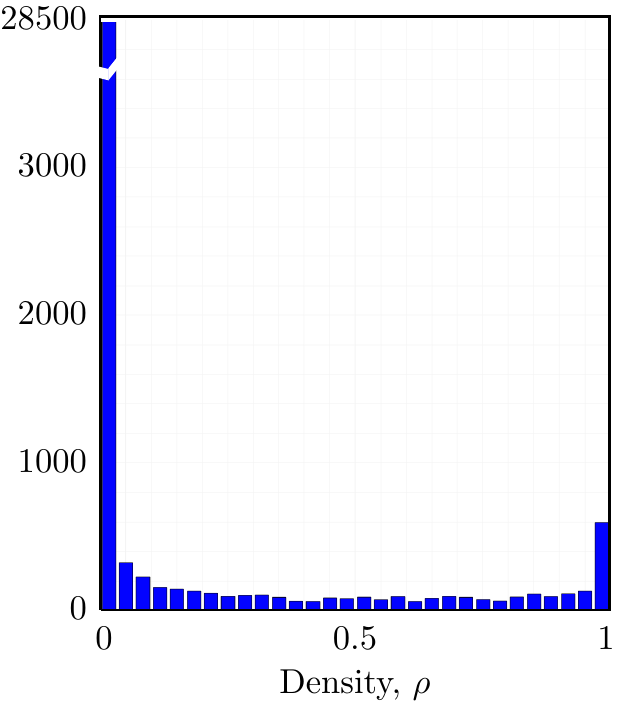}\\\includegraphics[scale = 0.5]{./figs/colorbar_ExII_horz}
		\caption{Adam design}
		\label{fig:exII_adam1}
	\end{subfigure}%
	\caption{Computed designs {with a density threshold of 0.5 and histogram of the density} for three methods, {namely, GCMMA, AdaGrad, and Adam} with $\lambda = 0$ for Example II {(side and front views {of the design} are shown for each method, respectively)}.}
	\label{fig:ExII_designs}
\end{figure}

\begin{figure}[htb!]
	\centering
	\begin{subfigure}[b]{\textwidth}
		\centering
		\includegraphics[scale =0.15]{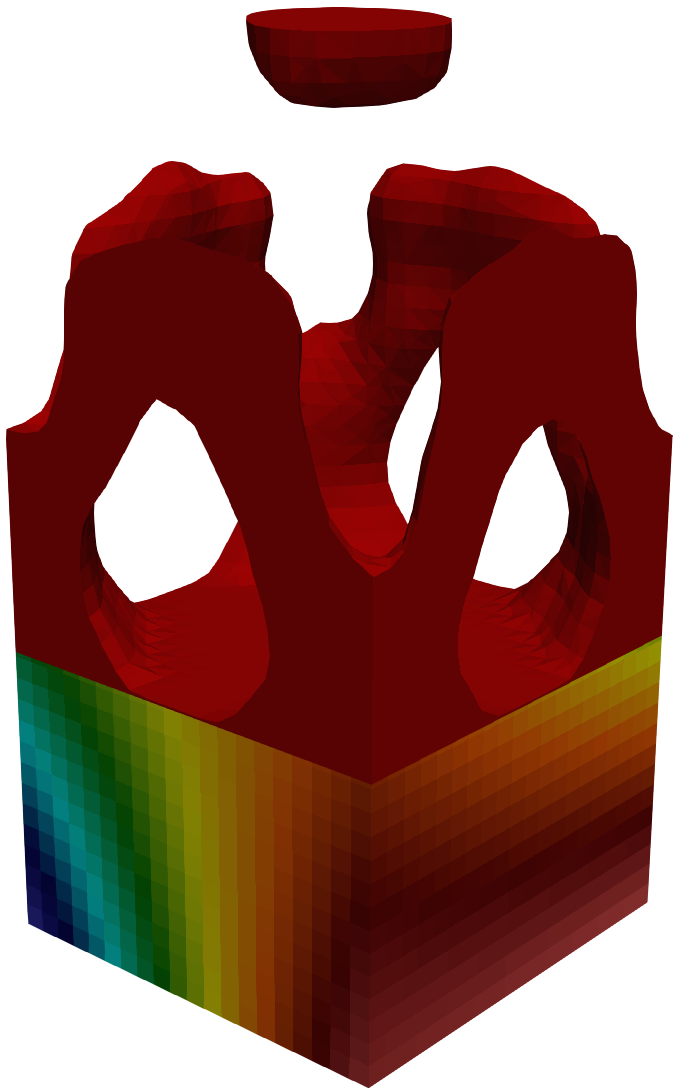}~~~~\includegraphics[scale =0.16]{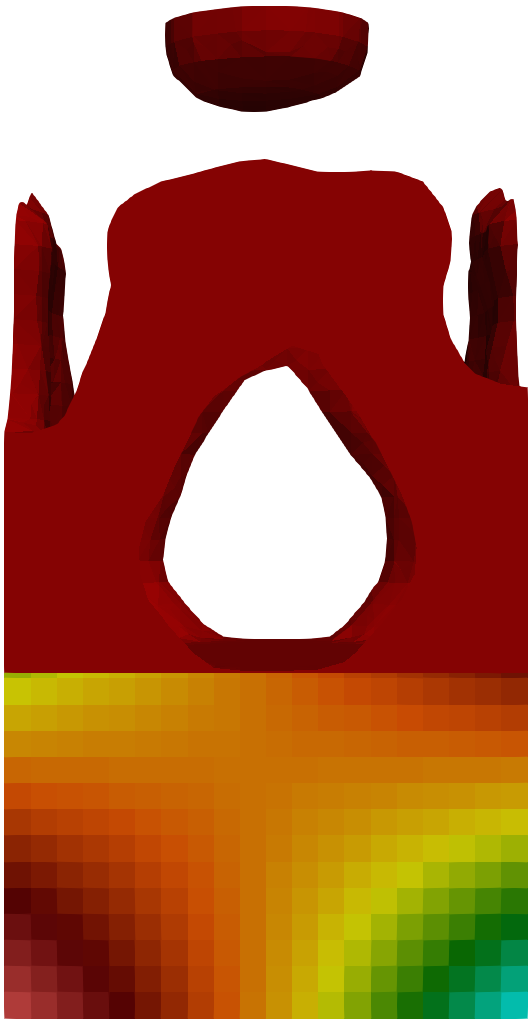}
		\\\includegraphics[scale = 0.5]{./figs/colorbar_ExII_horz}
		\caption{GCMMA design}
		\label{fig:exII_gcmma12}
	\end{subfigure}%
	\\
	\begin{subfigure}[b]{\textwidth}
		\centering
		\includegraphics[scale =0.15]{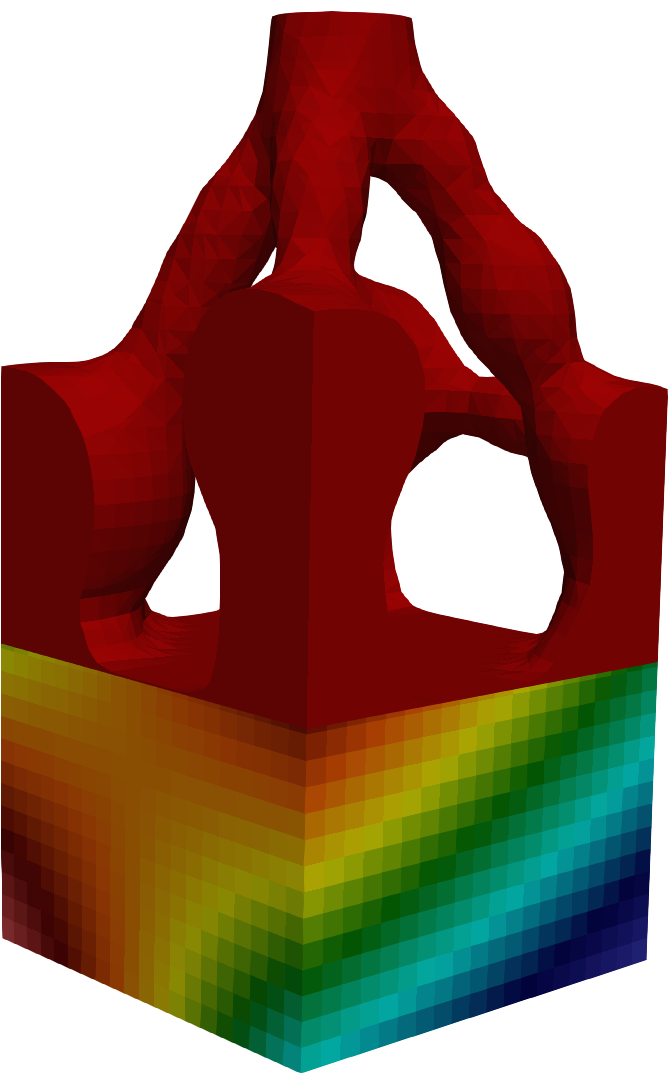}~~~~\includegraphics[scale =0.16]{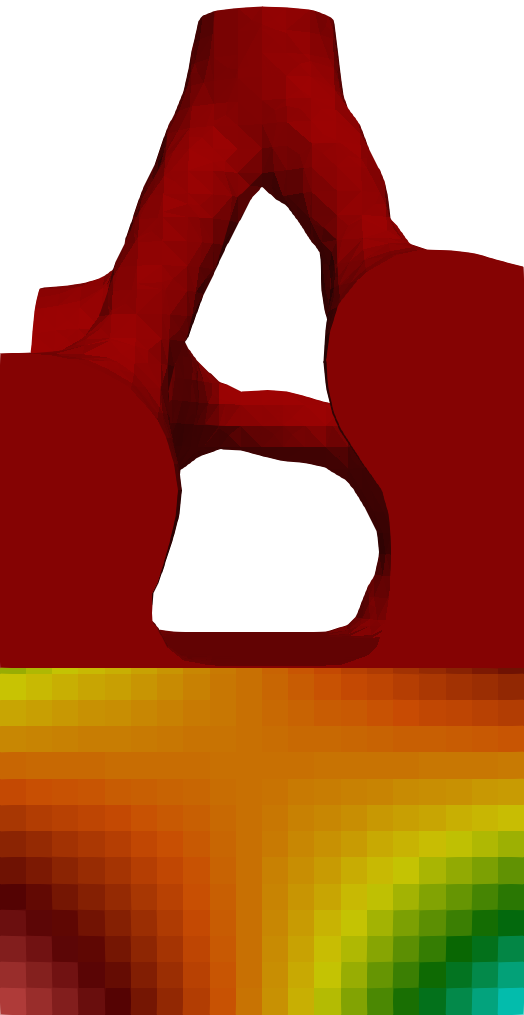}
		\\\includegraphics[scale = 0.5]{./figs/colorbar_ExII_horz}
		\caption{AdaGrad design}
		\label{fig:exII_adagrad12}
	\end{subfigure}%
	\\
	\begin{subfigure}[b]{\textwidth}
		\centering
		\includegraphics[scale =0.15]{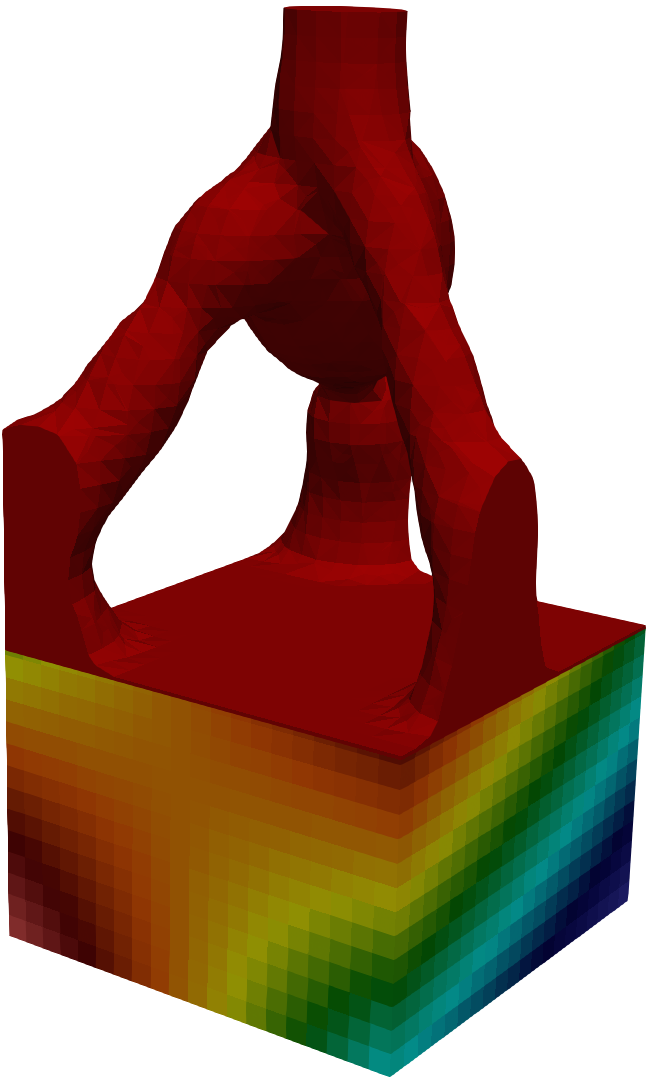}~~~~\includegraphics[scale =0.16]{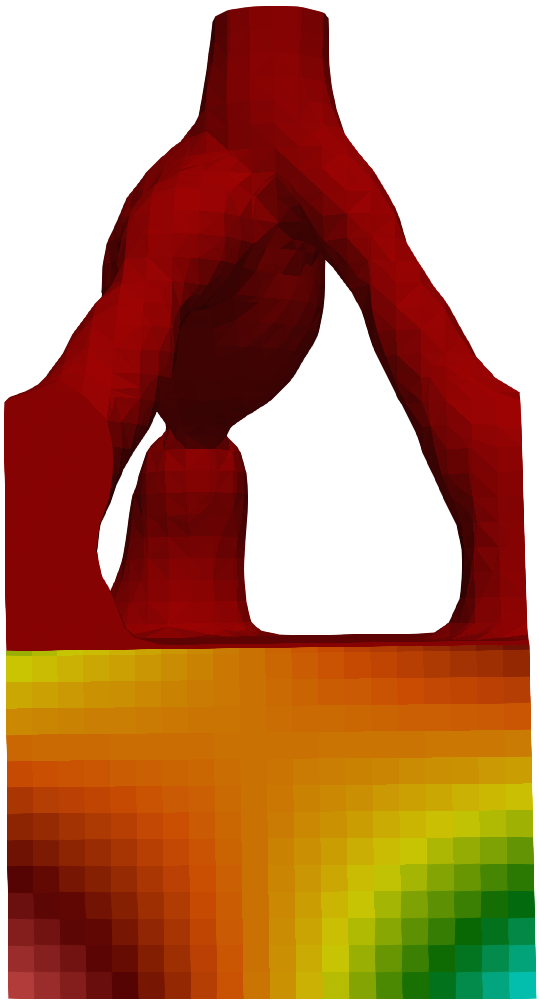}\\\includegraphics[scale = 0.5]{./figs/colorbar_ExII_horz}
		\caption{Adam design}
		\label{fig:exII_adam12}
	\end{subfigure}%
	\caption{{Computed designs with a density threshold of 0.7 for three methods, {namely, GCMMA, AdaGrad, and Adam} with $\lambda = 0$ for Example II {(side and front views of the design} are shown for each method, respectively)}.}
	\label{fig:ExII_designs2}
\end{figure}

\begin{figure}[htb!]
	\centering
	\includegraphics[scale=0.35]{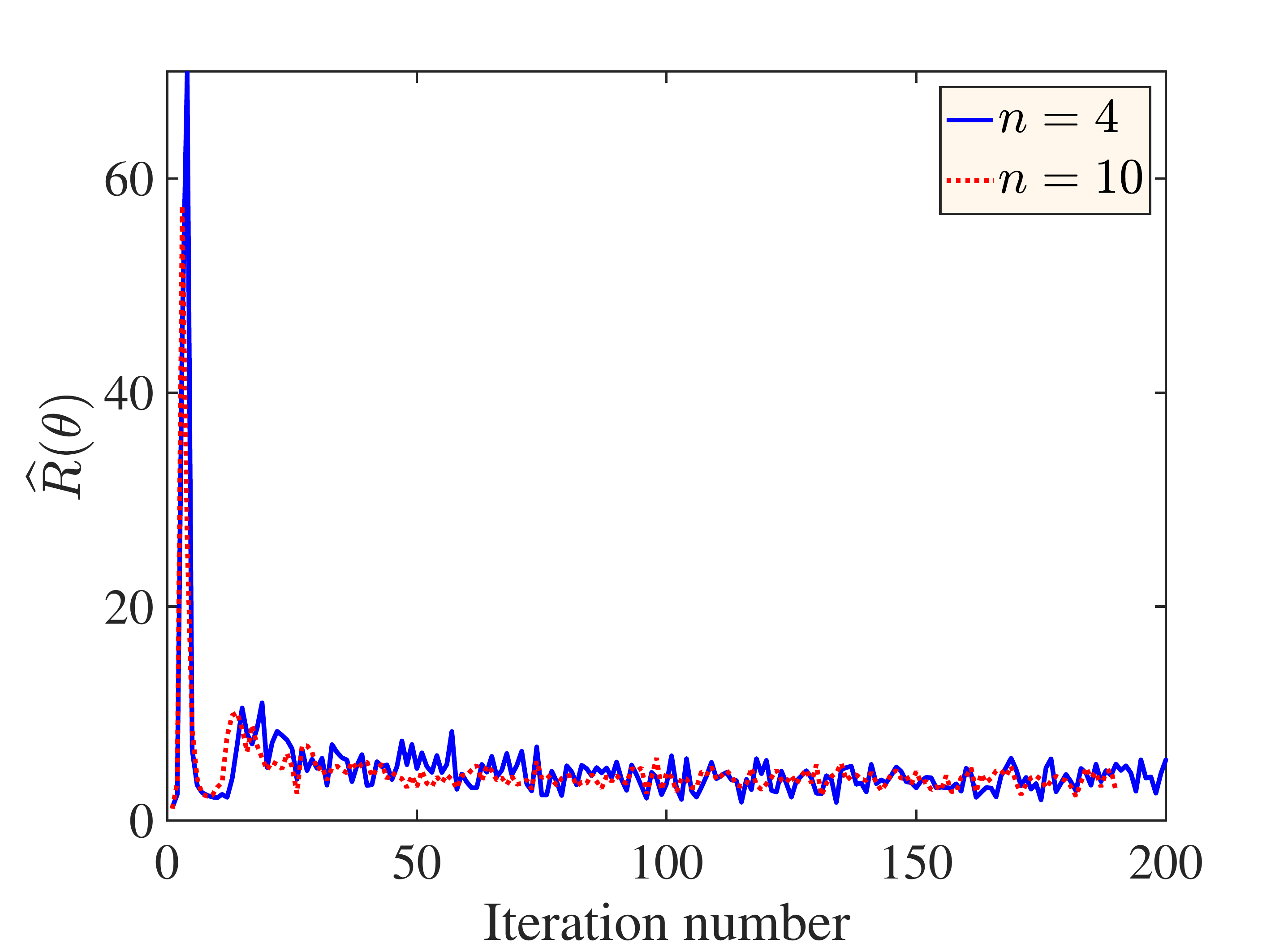}
	\caption{Plot of objectives using Adam with $n=4$ and $n=10$ samples per iteration with $\lambda = 0$ for Example II.}
	\label{fig:ExII_Adam_by_sample}
\end{figure}

{Figure~\ref{fig:ExII_designs} shows the optimal designs obtained using the three methods with a threshold of 0.5 on the material density. The final design obtained from these three methods are quite different -- GCMMA provides a canopy type design; AdaGrad leads to a design with four-legs and two horizontal members, and Adam delivers a three-legged design. When the material density threshold is increased to 0.7 GCMMA design disconnects as shown in Figure \ref{fig:ExII_designs2}. Adam and AdaGrad designs in Figure \ref{fig:ExII_designs2}, on the other hand, remain similar to their respective designs in Figure \ref{fig:ExII_designs}. 
The best designs obtained from all these three methods result in objective values that shows $98\sim99\%$
drop from the initial value, and all of them produce rather intuitive designs.
} 
%
%
When increasing the number of samples to $n=10$, we notice that using more samples at each iteration does not improve the optimization significantly; see Figure \ref{fig:ExII_Adam_by_sample} for Adam. The same observation of the convergence was made in Example I.

{The histograms of the density values in Figure~\ref{fig:ExII_designs}, however, show that the GCMMA performs poorly in producing a distinct $0-1$ design compared to AdaGrad and Adam. As GCMMA solves a subproblem constructed using gradients, presence of randomness in these gradients may results in poor convergence of the density values to distinct $0-1$ in the SIMP approach. Further analysis of GCMMA and potential strategies to improve its performance within SIMP is an interesting future work.}


\FloatBarrier
\subsubsection{Robust Design}
\label{subsec:block_var_opt}

Here, we consider solving the optimization problem (\ref{eq:exp_risk}) {subject to} (\ref{eq:glam_def}), for $\lambda = 0.01$, \textit{i.e.}, we seek to optimize a combination of mean and the variance of the objective subject to a similar combination of mean and variance of constraint values. We again utilize a learning rate parameter $\eta = 0.25$ for the algorithms in Section~\ref{sec:method_details} and a sample size of $n=4$.

\begin{figure}[htb!]
	\centering
	\begin{subfigure}[t]{\textwidth}
		\centering
		\includegraphics[scale = 0.35]{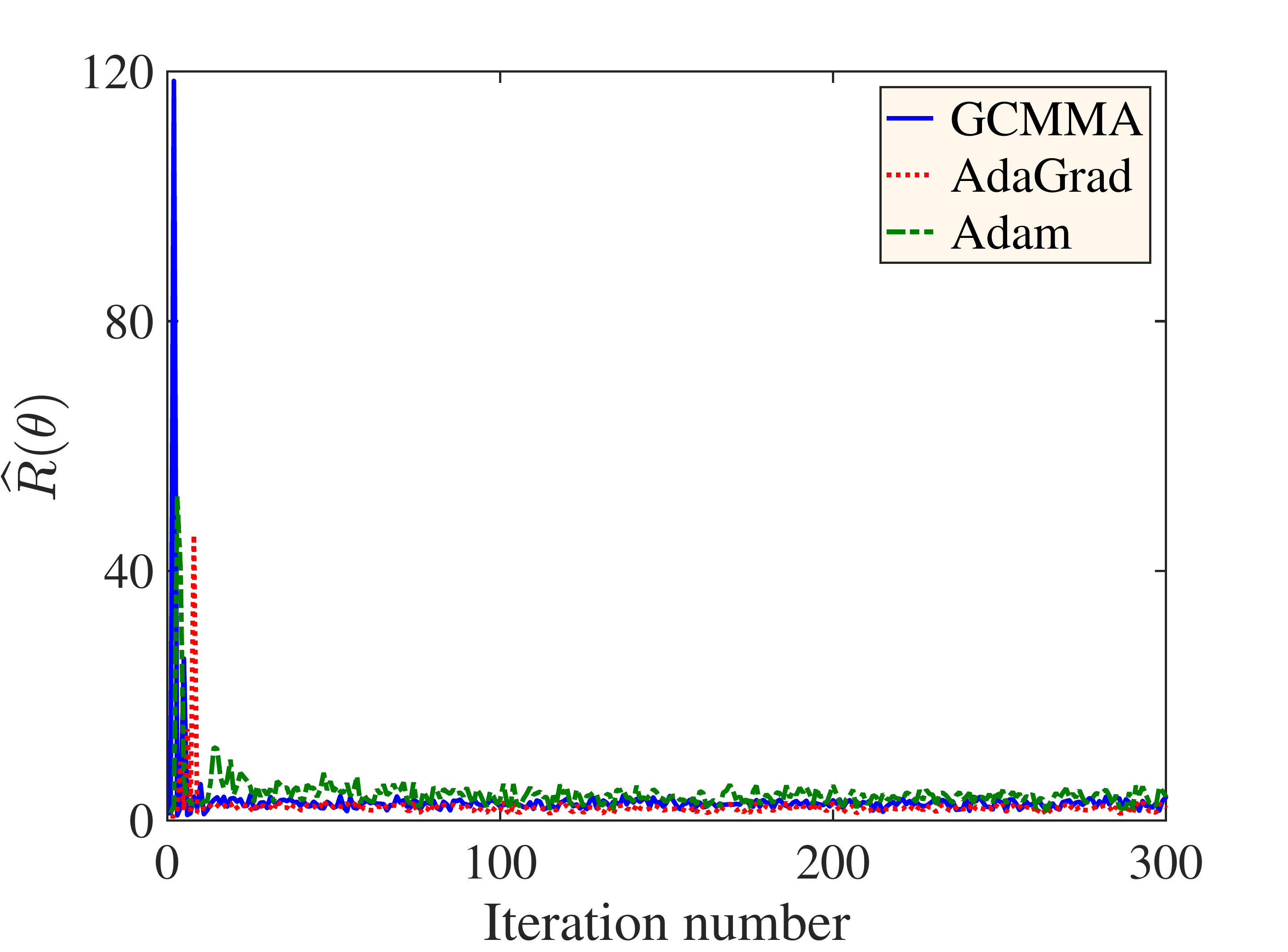}
		\caption{Objective estimates}
	\end{subfigure}
	\\
	\begin{subfigure}[t]{\textwidth}
		\centering
		\includegraphics[scale = 0.35]{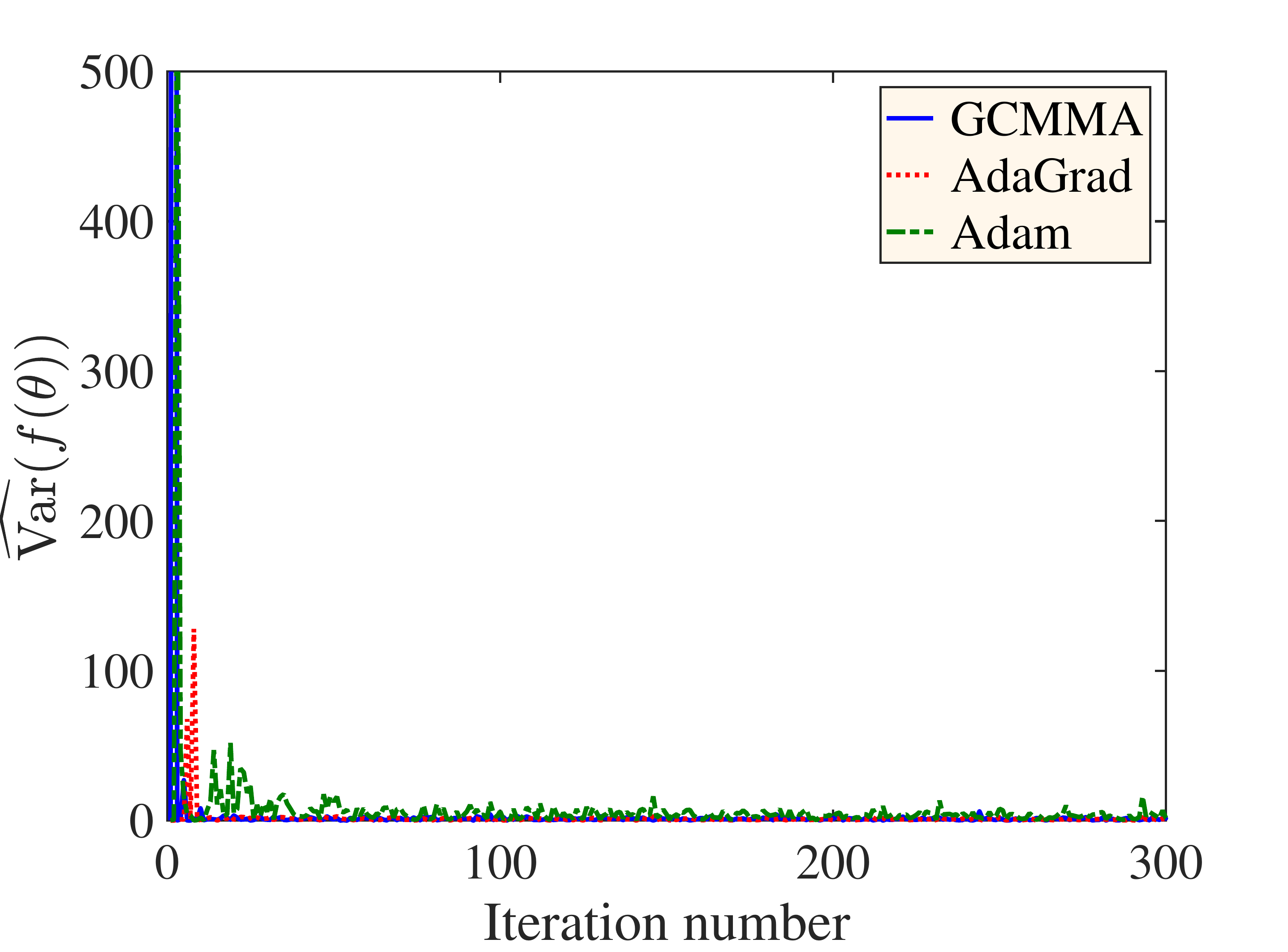}
		\caption{Variance estimates}\label{fig:ExII_var}
	\end{subfigure}
	\caption{Objectives and variances for the three methods used with $\lambda = 0.01$ for Example II.}
	\label{fig:ExII_var_graphs}
\end{figure}

\begin{figure}[htb!]
	\centering
	\begin{subfigure}[t]{\textwidth}
		\centering
		\includegraphics[scale =0.37]{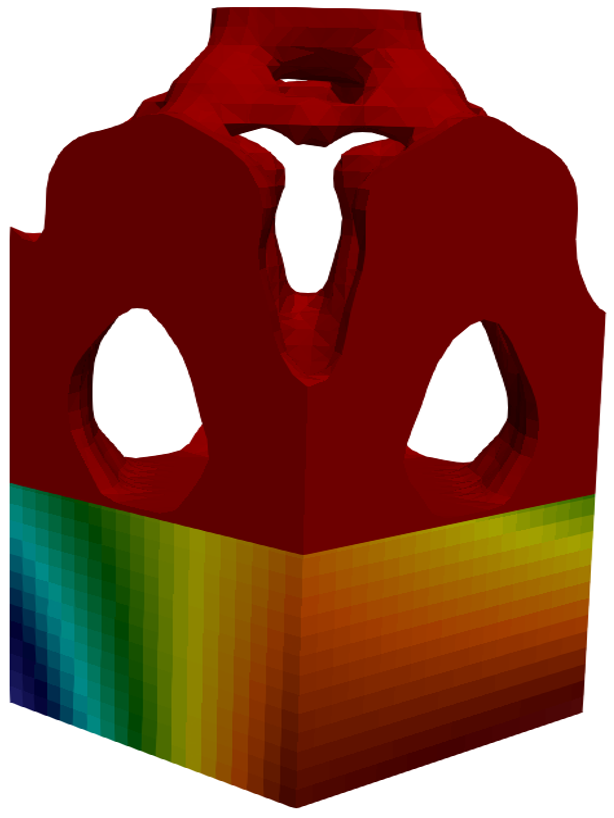}~~~~\includegraphics[scale =0.38]{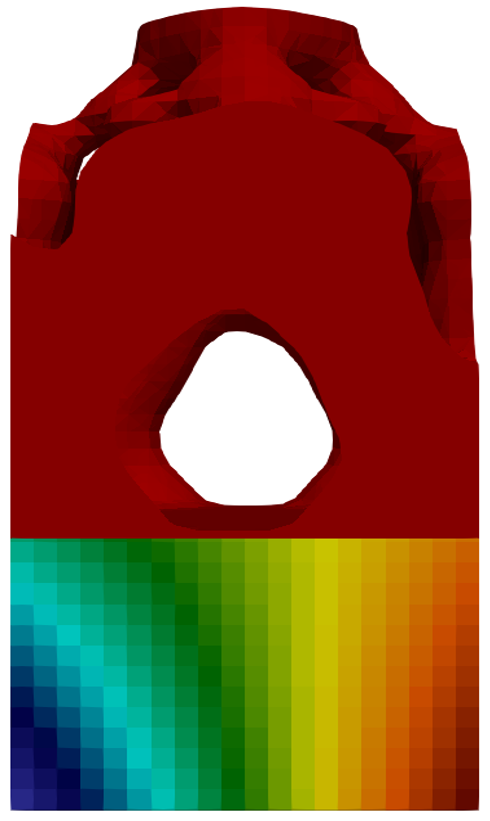}~~~~\includegraphics[scale =0.8]{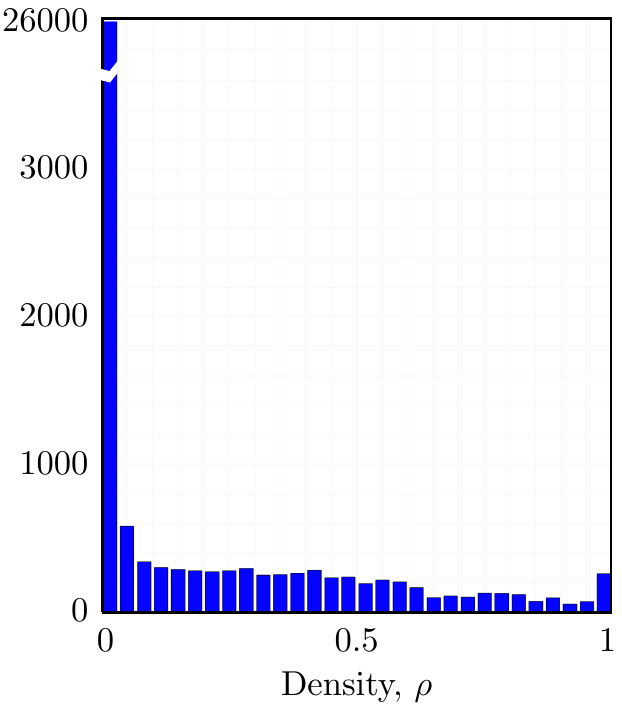}\\\includegraphics[scale = 0.5]{./figs/colorbar_ExII_horz}
		\caption{GCMMA design}
		\label{fig:exII_gcmma2}
	\end{subfigure}%
	\\
	\begin{subfigure}[t]{\textwidth}
		\centering
		\includegraphics[scale =0.37]{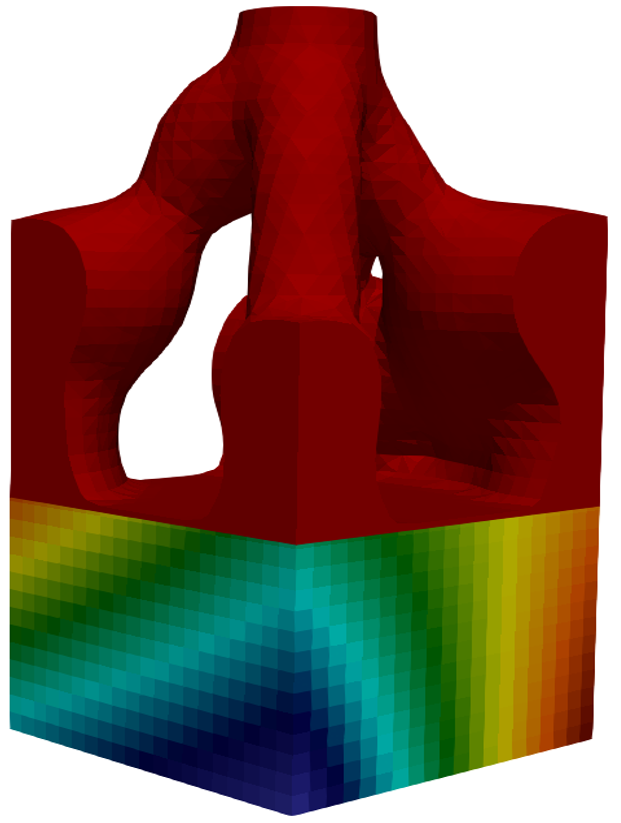}~~~~\includegraphics[scale =0.38]{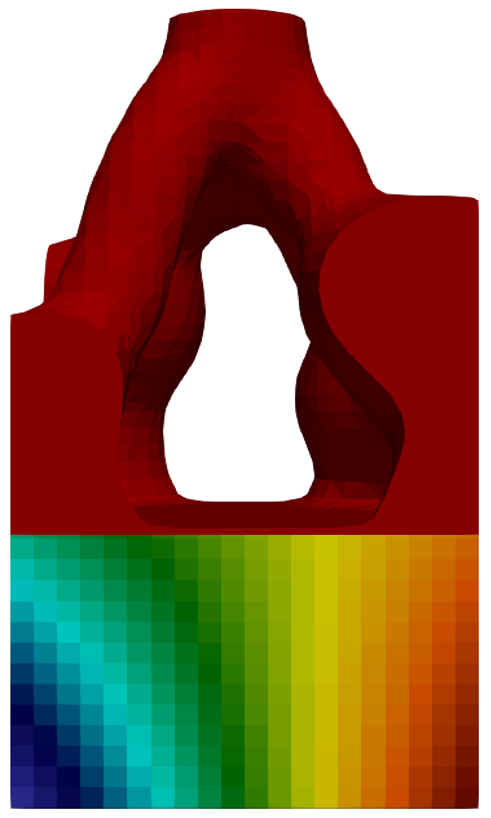}~~~~\includegraphics[scale =0.8]{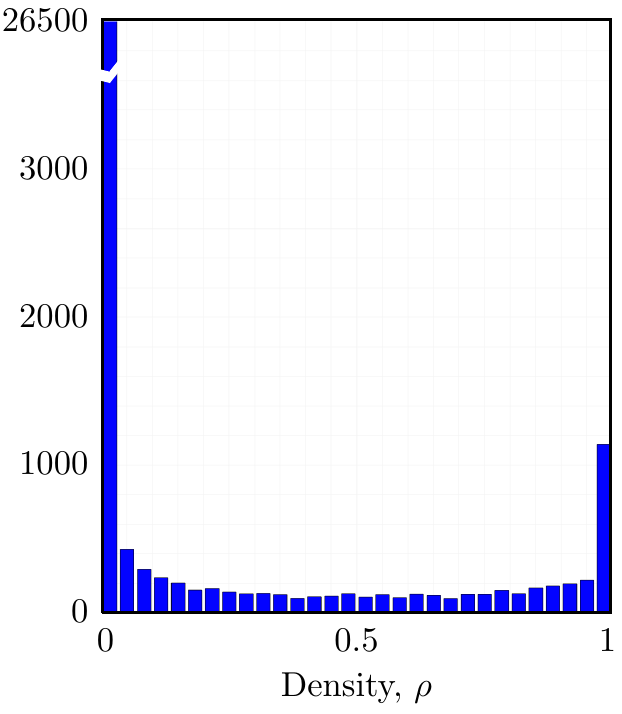}\\\includegraphics[scale = 0.5]{./figs/colorbar_ExII_horz}
		\caption{AdaGrad design}
		\label{fig:exII_adagrad2}
	\end{subfigure}%
	\\
	\begin{subfigure}[t]{\textwidth}
		\centering
		\includegraphics[scale =0.39]{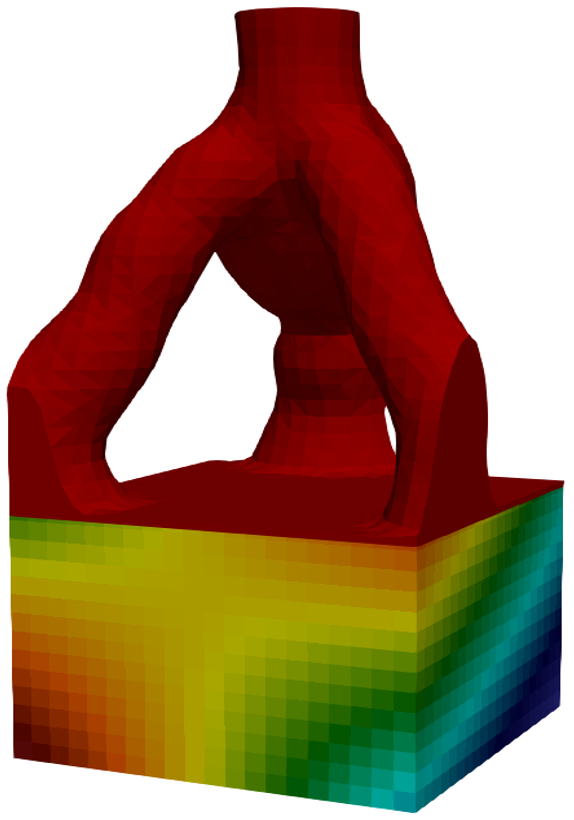}~~~~\includegraphics[scale =0.42]{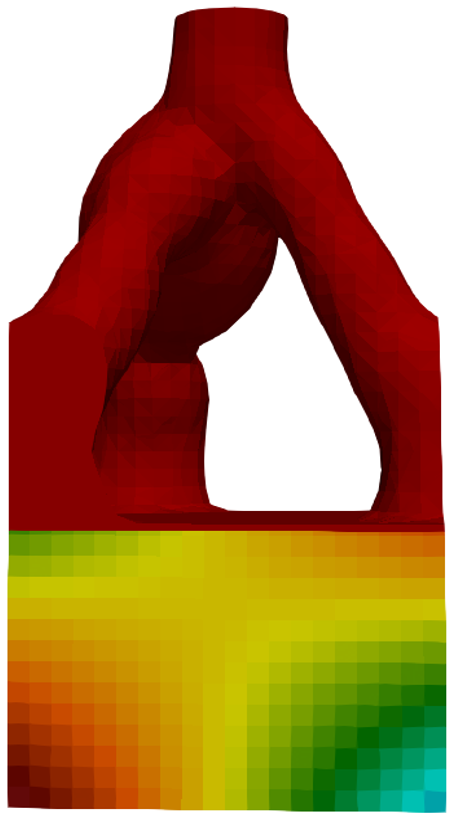}~~~~\includegraphics[scale =0.8]{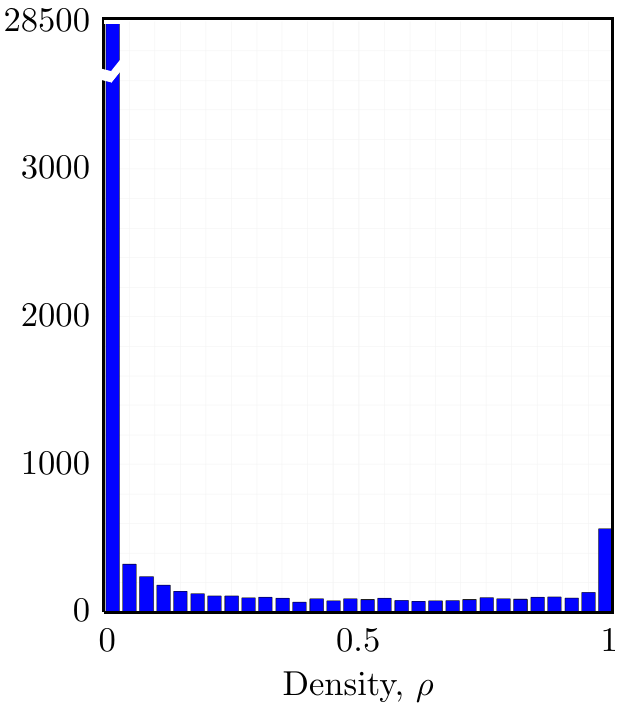}\\\includegraphics[scale = 0.5]{./figs/colorbar_ExII_horz}
		\caption{Adam design}
		\label{fig:exII_adam2}
	\end{subfigure}%
	\caption{Computed designs {with a density threshold of 0.5 and histogram of the density} for three methods, {namely, GCMMA, AdaGrad, and Adam} with $\lambda = 0.01$ for Example II {(side and front views {of the design} are shown for each method, respectively)}.}
	\label{fig:ExII_designs_var}
\end{figure}

Figure~\ref{fig:ExII_var_graphs} shows the evolution of the objective $\widehat{R}(\ppm)$ and variance $\widehat{\mathrm{Var}}(f(\ppm))$ for designs determined by GCMMA, AdaGrad, and Adam. The proposed methods show rapid decrease in the objective $\widehat{R}(\ppm)$ and variance $\widehat{\mathrm{Var}}(f(\ppm))$ in the first few iterations. Figure \ref{fig:ExII_var} indicates that the proposed use of stochastic gradients for TOuU is useful to produce a robust design with a smaller variance of $f(\ppm)$. Figure~\ref{fig:ExII_designs_var} shows the final designs obtained using the three methods. {While GCMMA provides a canopy type design as before, AdaGrad and Adam provide four- and three-legged designs.} However, if we compare their objectives, {
the best designs obtained from all these three methods result in objective values that shows $97\sim99\%$ drop from the initial value. 
All of these three methods produce rather intuitive designs as well. 
}
%


\begin{figure}[htb!]
	\centering
	\begin{subfigure}[b]{\textwidth}
		\centering
		\includegraphics[scale =0.15]{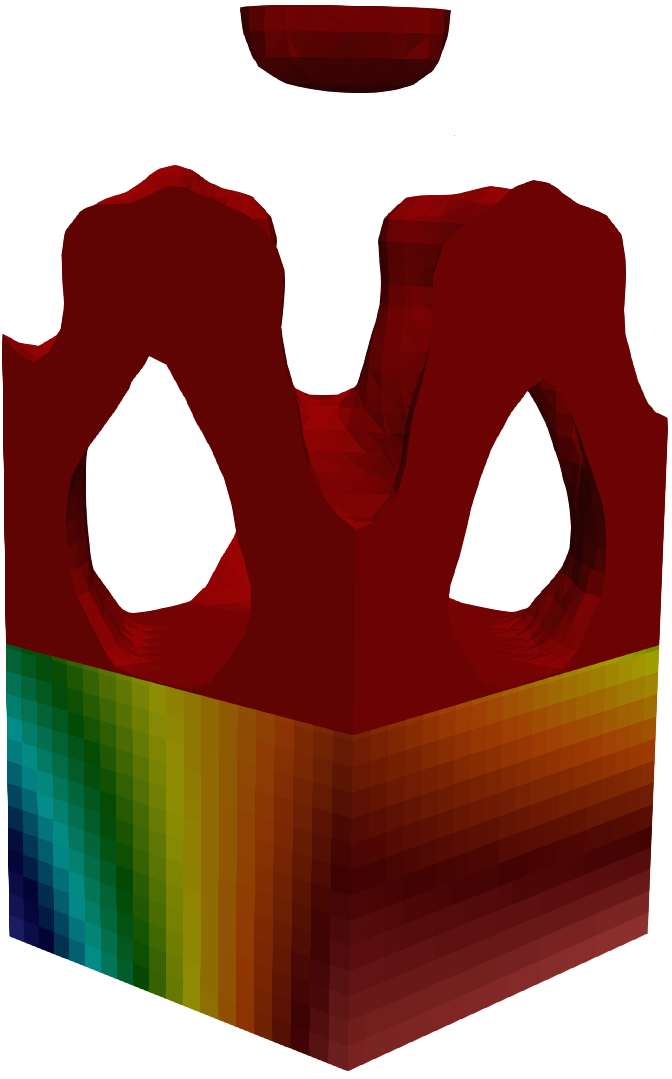}~~~~\includegraphics[scale =0.16]{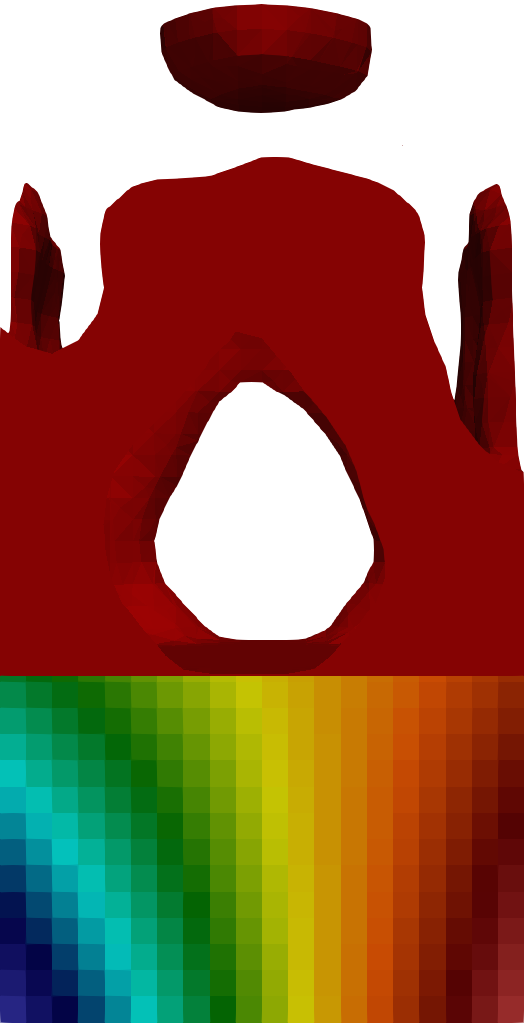}
		\\\includegraphics[scale = 0.5]{./figs/colorbar_ExII_horz}
		\caption{GCMMA design}
		\label{fig:exII_gcmma22}
	\end{subfigure}%
	\\
	\begin{subfigure}[b]{\textwidth}
		\centering
		\includegraphics[scale =0.15]{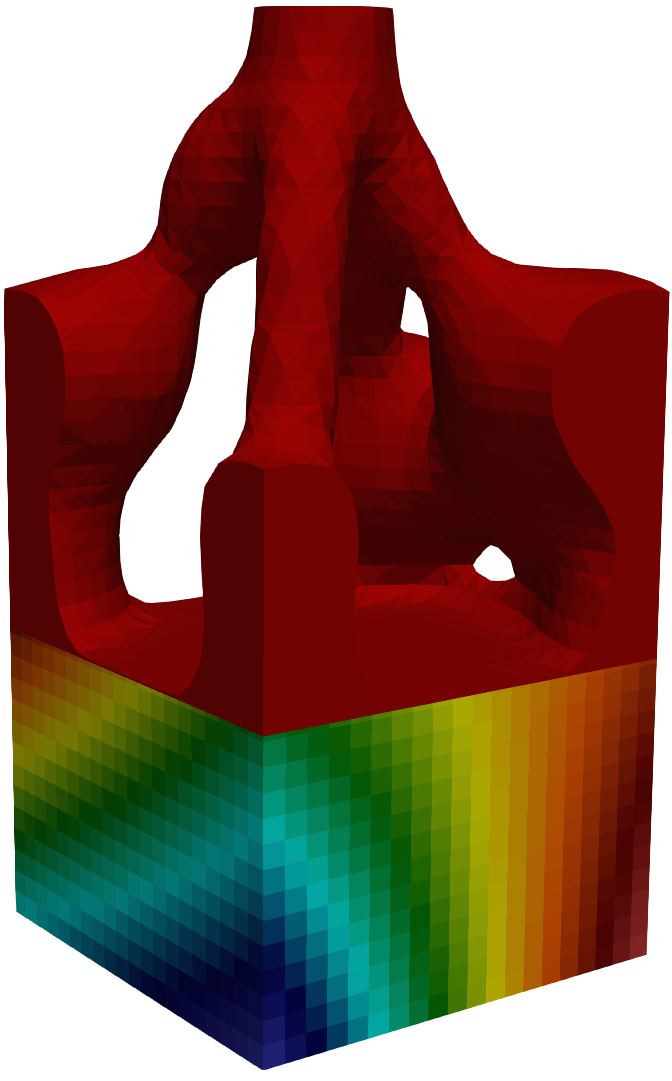}~~~~\includegraphics[scale =0.16]{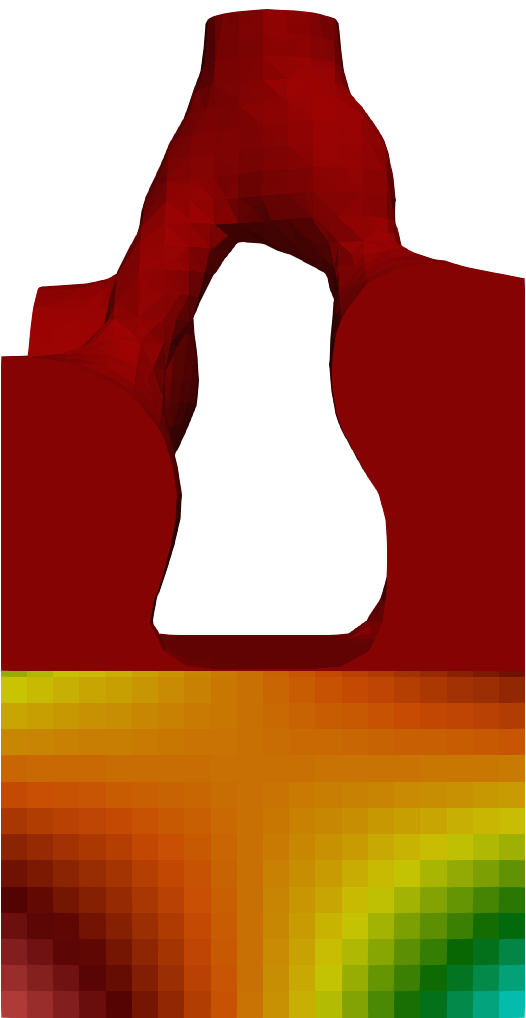}
		\\\includegraphics[scale = 0.5]{./figs/colorbar_ExII_horz}
		\caption{AdaGrad design}
		\label{fig:exII_adagrad22}
	\end{subfigure}%
	\\
	\begin{subfigure}[b]{\textwidth}
		\centering
		\includegraphics[scale =0.15]{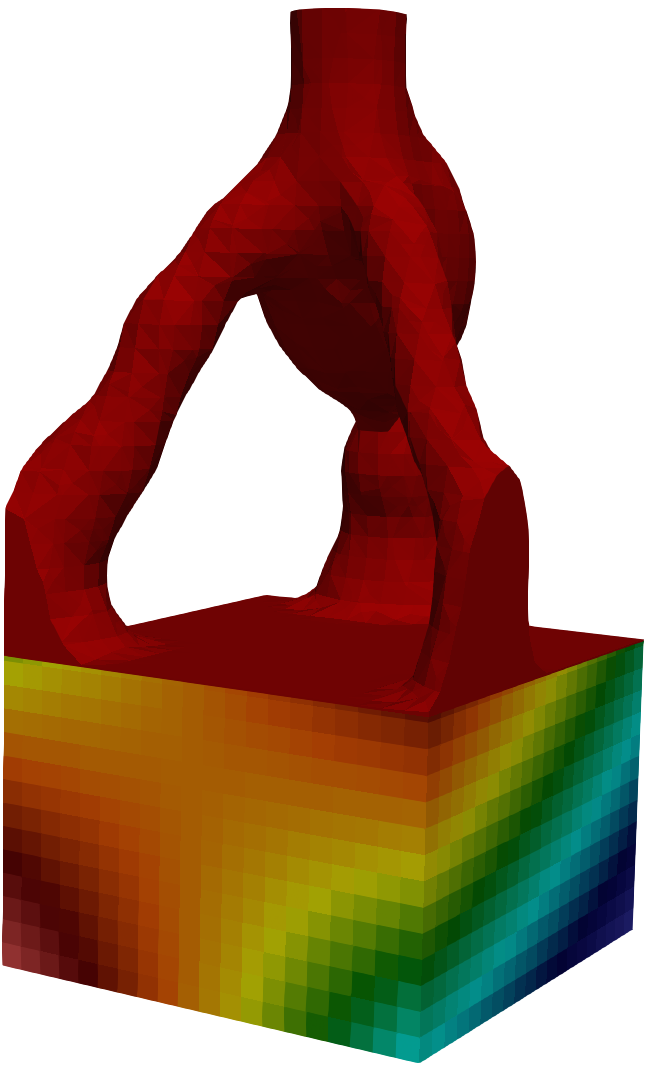}~~~~\includegraphics[scale =0.16]{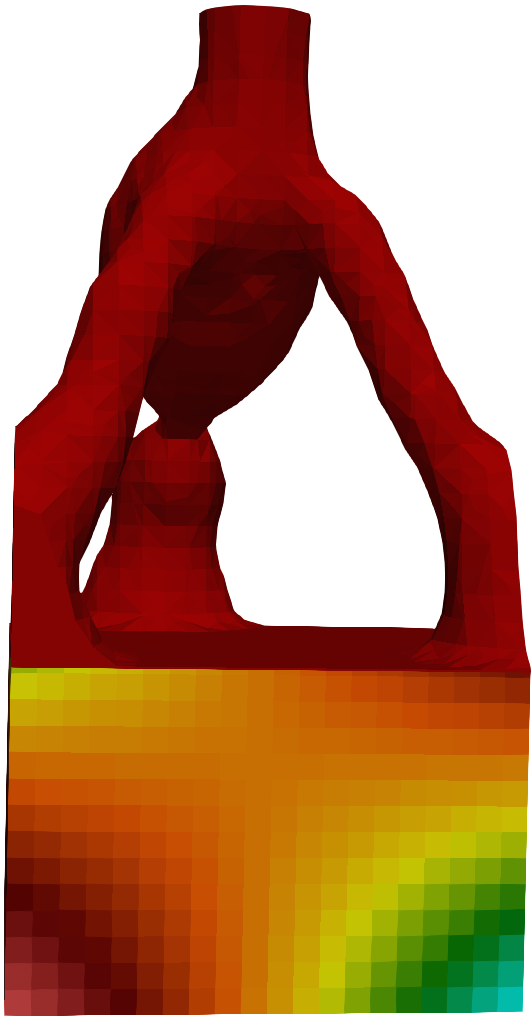}\\\includegraphics[scale = 0.5]{./figs/colorbar_ExII_horz}
		\caption{Adam design}
		\label{fig:exII_adam22}
	\end{subfigure}%
	\caption{{Computed designs with a density threshold of 0.7 for three methods, {namely, GCMMA, AdaGrad, and Adam} with $\lambda = 0.01$ for Example II {(side and front views of the design} are shown for each method, respectively)}.}
	\label{fig:ExII_designs_var2}
\end{figure}

{In this example, Adam and AdaGrad perform well in the presence of uncertainty even though final designs obtained from them are different. The GCMMA design in Figure \ref{fig:exII_gcmma2} has more pronounced and thicker bars at the top of the canopy that connects the four sides to the top of the structure compared to the average design  in Figure \ref{fig:exII_gcmma1}. For the AdaGrad design the robust structure removes the two horizontal bars connecting the legs on two sides as seen in Figure \ref{fig:exII_adagrad1} but produces thicker legs. The Adam design, however, remains similar. 
Figure \ref{fig:ExII_designs_var2} shows the designs obtained using these three methods when the material density threshold is increased to 0.7. The GCMMA design again disconnects showing its poor convergence to a $0-1$ design.
The histograms of the density values of the final designs in Figure \ref{fig:ExII_designs_var} show that the Adam and AdaGrad produce better $0-1$ designs in terms of density than GCMMA as before.
}
\FloatBarrier
%


\section{\texorpdfstring{Conclusions}{Conclusions}}
\label{sec:conc}

{ Topology optimization under high-dimensional uncertainty that utilizes gradient-based approaches requires multiple evaluations of the objective, constraints, and the gradients. 
This incurs a significant computational cost and TOuU often becomes practically infeasible.
In this paper, we reduce the computational cost of TOuU by using a stochastic approximation of the gradients. These approximations are calculated with only a few random samples per iteration (\textit{e.g.}, $n=4$ in the numerical examples). This approach helps us reduce the cost of TOuU to only a modest multiple of the one of deterministic TO. To implement our proposed approach, we study some of the 
SGD methods popular in machine learning while supplying them with stochastic gradients. 
We also consider the GCMMA in our numerical examples while supplying stochastic gradients to it. 
The first TOuU problem, a highly non-convex beam design problem, shows that while GCMMA can often perform well, it requires fair estimates of the gradients when the weight on the variance term in the objective is large. Among the five SGD methods employed here, AdaGrad and Adam perform well in this example. The damping of gradients in Adagrad and Adam helps to navigate the design space. 
In the second TO problem,  a load is supported over an uncertain elastic bedding with uncertain loading direction. In this example, GCMMA, AdaGrad, and Adam again produce designs with comparable objectives. {However, GCMMA performs poorly in achieving a $0-1$ design.} Both of these examples show that the use of stochastic gradients estimated with a small number of random samples at each iteration produces meaningful designs in a computationally inexpensive manner. Further, the optimized structure obtained using SGD methods has objective values similar to GCMMA for topology optimization under uncertainty. 
}

\section{\texorpdfstring{Acknowledgements}{Acknowledgements}}
\label{sec:ack}

The authors acknowledge the support of the Defense Advanced Research Projects Agency (DARPA) TRADES project {under agreement HR0011-17-2-0022}.
Computations for the results presented in Section~\ref{sec:examples} are performed on Lonestar5, a high-performance computing resource operated by the Texas Advanced Computing Center (TACC) at University of Texas, Austin.

\section{\texorpdfstring{Conflict of Interest}{Conflict of Interest}}
\label{sec:coi}

On behalf of all authors, the corresponding author states that there is no conflict of interest. 

\section{\texorpdfstring{Replication of Results}{Replication of Results}}
\label{sec:replication}

The optimization algorithms of Sections \ref{sec:sgd} and \ref{sec:gcmma} have been implemented in MATLAB and are accessible via the GitHub page \url{https://github.com/CU-UQ/TOuU}. The link also includes the necessary input files for the elastic bedding example of Section \ref{sec:elastic_bedding}. The TO and XFEM calculations were performed using an in-house solver that is not at the stage of being publicly available. However, the MATLAB codes used in this study can interface with other TO and XFEM codes.

\bibliographystyle{apalike} 
\bibliography{citations}

\end{document}